\documentclass[11pt]{article}
\usepackage{amssymb,mathrsfs,latexsym,bm,cite}
\usepackage[boxed]{algorithm2e}

\begin{document}
\newcommand{\qed}{\hphantom{.}\hfill $\Box$\medbreak}
\newcommand{\proof}{\noindent{\bf Proof \ }}
\newtheorem{Theorem}{Theorem}[section]
\newtheorem{Lemma}[Theorem]{Lemma}
\newtheorem{Corollary}[Theorem]{Corollary}
\newtheorem{Remark}[Theorem]{Remark}
\newtheorem{Example}[Theorem]{Example}
\newtheorem{Definition}[Theorem]{Definition}
\newtheorem{Construction}[Theorem]{Construction}
\newtheorem{Proposition}[Theorem]{Proposition}
\renewcommand{\algorithmcfname}{Figure}

\begin{center}
{\Large\bf Combinatorial constructions for optimal two-dimensional
optical orthogonal codes with $\lambda=2$} \footnote{Supported by
NSFC grant No. $61071221$, $10831002$ (Y. Chang), and NSFC grant No.
$10901016$ (T. Feng)}

\vskip24pt

Tao Feng, Yanxun Chang \\ Institute of Mathematics\\ Beijing
Jiaotong University\\ Beijing 100044, P. R. China \\ {\tt
tfeng@bjtu.edu.cn}\\ {\tt yxchang@bjtu.edu.cn}

\vskip12pt

\end{center}

\noindent {\bf Abstract:} In this paper, we are concerned about optimal
two-dimensional optical orthogonal codes with $\lambda=2$. Some
combinatorial constructions are presented and many infinite families
of optimal two-dimensional optical orthogonal codes with weight $4$
and $\lambda=2$ are obtained. Especially, we shall see that in many
cases an optimal two-dimensional optical orthogonal code can not
achieve the Johnson bound.



\section{Introduction}

An optical orthogonal code is a family of sequences with good auto-
and cross-correlation properties. Its study has been motivated by an
application in an optical code-division multiple access (OCDMA)
system. In a bursty traffic environment of a multiple access local
area network (LAN), asynchronous multiplexing schemes are more
efficient than synchronous schemes. OCDMA is such one asynchronous
multiplexing scheme suitable for high speed LANs. For more
information, the interested reader may refer to
\cite{kpp,ml,mmc,s,sb}.

In an OCDMA system different users share both time and frequency,
and are distinguished by using a unique spreading sequence. Each
user's data is multiplied by its spreading sequence, and then all
the users are coupled into the shared channel. Optical orthogonal
codes can be taken as the spreading sequences used in an OCDMA
system.

Let $u$, $v$, $k$ and $\lambda$ be positive
integers. A  {\em two-dimensional $(u\times v,k,\lambda)$ optical
orthogonal code} (briefly $2$-D $(u\times v,k,\lambda)$-OOC),
$\cal{C}$, is a family of $u\times v$ $(0, 1)$-matrices (called {\em
codewords}) of Hamming weight $k$ satisfying: for any matrix
${\mathbf{A}}=(a_{ij})_{u\times v}\in\cal{C}$,
${\mathbf{B}}=(b_{ij})_{u\times v}\in\cal{C}$ and any integer $r$:
$$\sum_{i=0}^{u-1}\sum_{j=0}^{v-1}a_{ij}b_{i,j+r}\leq\lambda,$$
where either $\mathbf{A}\neq \mathbf{B}$ or $r\neq 0$, and the
arithmetic $j+r$ is reduced modulo $v$. Especially, when $u=1$, a two-dimensional $(1\times v,k,\lambda)$ optical orthogonal code is said to be a \emph{one-dimensional $(v,k,\lambda)$-optical orthogonal code}, denoted by $1$-D $(v,k,\lambda)$-OOC.

$1$-D OOC was first suggested in $1989$ \cite{csw}. Since then much work has been done on $1$-D OOCs. The interested reader may refer to \cite{ab,am,am2,be,b,cfm,cj,cm,cy,cc,cc1,ck,fcj,fm,gms,gy,mc1,mc,mms,mokl,mzkz,ngm,wc,y}.
One limitation in applying $1$-D OOCs is that the length of the
sequences increases rapidly when the number of users or the weight
of codes is increased, which means a large bandwidth expansion is
required. Thus the bandwidth utilization is reduced. And a large code length
causes the chip rate of the OCDMA system to exceed the maximum chip
rate currently attainable in practice.

$1$-D OOCs spread the input data bits only in the time domain. By spreading in both time and wavelength domain, the chip rate can be
reduced considerably. Technologies such as
wavelength-division-multiplexing (WDM) and dense-WDM have made it
possible to spread codes in both time and wavelength domain \cite{yk}.
These codes are referred to wavelength-time hopping codes,
multiple-wavelength codes, two-dimensional optical orthogonal
codes, etc., which tend to require smaller code length and hence
lower chip rate. Here we always refer to these codes as
two-dimensional optical orthogonal codes.

The number of codewords of a $2$-D OOC is called the {\em size} of the $2$-D OOC. From a practical point of view, a code with a large size is required \cite{sb}. For fixed values of $u$, $v$, $k$ and $\lambda$, the largest possible size of a $2$-D
$(u\times v,k,\lambda)$-OOC is denoted by $\Phi(u\times v,k,\lambda)$. A $2$-D $(u\times v,k,\lambda)$-OOC
with $\Phi(u\times v,k,\lambda)$ codewords is said to be {\em
optimal}. Generally speaking, it is difficult to determine the exact value of $\Phi(u\times v,k,\lambda)$. Based on the Johnson bound \cite{John} for constant weight
codes, the size of a $2$-D $(u\times v,k,\lambda)$-OOC is upper bounded
\cite{yk} by
$$
 \Phi(u\times v,k,\lambda)\leq J(u\times v,k,\lambda),
$$
where
$$J(u\times v,k,\lambda)=\lfloor\frac{u}{k}\lfloor\frac{uv-1}{k-1}
\lfloor\frac{uv-2}{k-2}\lfloor\cdots\lfloor\frac{uv-\lambda}{k-\lambda}
\rfloor\cdots\rfloor\rfloor\rfloor\rfloor.$$

In optical code-division multiple-access applications, performance analysis shows that codes with $\lambda\leq 3$ are the most desirable. As pointed out by \cite{ms}, from a multiple-access and synchronization point of view, the most desirable on-off signature sequences are OOCs with $\lambda=1$. However, these families of codes may suffer from low
cardinality in some applications. It was hinted that in \cite{asl} OOCs with $\lambda=2$ could, under certain conditions, have better performance than that of OOCs with $\lambda=1$. In this paper, we are concerned about OPTIMAL $2$-D $(u\times v,k,\lambda)$-OOCs with $\lambda=2$ and $k=4$.

We will neither try to explore the applications of $2$-D OOCs, nor try to provide the performance analysis of a code-division multiple-access system which uses $2$-D OOCs. Mathematically, combinatorial design theory, projective geometry and finite field theory are three main tools to investigate the constructions for $2$-D OOCs. In this paper, we focus our attention on the combinatorial structures of $2$-D OOCs. Many terminologies and results related to combinatorial design theory will be used. To ensure smooth reading of the paper, most of the proofs related to design theory have been moved to the Appendices. For more information on design theory, the interested reader may refer to \cite{bjl}.

\subsection{Literature review}

There is a considerable literature  on $2$-D OOC constructions. Yang and Kwong \cite{yk} used a $1$-D OOC to achieve spreading in the wavelength and time domains to construct a $2$-D OOC. The construction by Lee and Seo \cite{ls} spreads in the wavelength
and the time domain by using two different $1$-D OOCs.
Sun et al. \cite{sywx} constructed a $2$-D OOC by employing a frequency
hopping code and a $1$-D OOC to spread in the wavelength domain and
the time axis, respectively. The construction by Alderson and Mellinger \cite{am1} are based on certain point sets in finite projective spaces of dimension $k$ over GF$(q)$.
Omrani et al. \cite{ogkeb} constructed some $2$-D OOCs using polynomials over finite fields and rational functions. Cao and Wei \cite{cws} first gave a combinatorial description of $2$-D OOCs. Wang et al. \cite{wsy} discussed the existence of optimal $2$-D OOCs with weight $3$ and index $\lambda=1$ using combinatorial design theory.

For more information on $2$-D OOC constructions, the interested reader may refer to \cite{am1,cws,gw,ky,ls,lyqwx,mghbl,ogkeb,sowk,sss,sss1,sywx,wsy,wy,yk96,yk,ycb} and the references therein. However, in this paper, we only focus our attention on OPTIMAL $2$-D OOC constructions. In applications, optimal OOCs facilitate the largest possible number of asynchronous users to transmit information effectively and reliably. A quick review of the majority of constructions for optimal $2$-D OOCs in the literature is presented in Table I.

 {\footnotesize
\tabcolsep 0.07in
\begin{center}
{\bf Table I\\ Optimal $2$-D OOCs in the literature} \vspace{6pt}

\begin{tabular}{|c|c|c|c|}\hline
Parameters & Conditions & Code size & Reference\\\hline
$(u\times v,3,1)$ & $u,v\geq 1$ and $v\equiv 1\ ({\rm mod }\ 2)$ & if $u\equiv 5\ ({\rm mod }\ 6)$ and $v=1$ &\cite{wsy}\\
&& $J(u\times v,3,1)-1$; & \\
&&  otherwise, $J(u\times v,3,1)$ & \\\hline
$(u\times u,k,1)$ & $u\equiv 1\ ({\rm mod }\ k(k-1))$ and $u$ a prime & $J(u\times u,k,1)$ & \cite{yk}\\\hline

$(p^n\times p,p,1)$ & $p$ a prime and $n\geq 1$ & $J(p^n\times p,p,1)$ & \cite{cws}\\\hline

$(u\times v,q,1)$ & $uv=q^n-1$, $n\geq 1$, and $q$ a prime power& $J(u\times v,q,1)$ & \cite{am1}\\\hline
$(u\times v,q+1,1)$ & $uv=(q^{n+1}-1)/(q-1)$, & $J(u\times v,q+1,1)$ & \cite{am1}\\
& $q$ a prime power, $n\geq 1$,&&\\
& either $n\equiv 0\ ({\rm mod }\ 2)$, &&\\
& or $n\equiv 1\ ({\rm mod }\ 2)$ and gcd$(q+1,v)=1$ && \\\hline

$(u\times v,4,2)$ & $uv=2^n-1$ and $n\geq 3$ & $J(u\times v,4,2)$ &\cite{am1}\\\hline
$(u\times v,6,2)$ & $uv=(4^n-1)/3$, $n\geq 3$, & $J(u\times v,6,2)$ &\cite{am1}\\
& either $n\equiv 0,1\ ({\rm mod }\ 3)$, &&\\
& or $n\equiv 2\ ({\rm mod }\ 3)$ and gcd$(21,v)=1$ &&\\\hline
$(u\times v,q+1,2)$ & $uv=q^n+1$, $q$ a prime power, & $J(u\times v,q+1,2)$ &\cite{am1}\\
& $n\geq 1$, either $n\equiv 0\ ({\rm mod }\ 2)$, &&\\
& or $n\equiv 1\ ({\rm mod }\ 2)$ and gcd$(q+1,v)=1$ &&\\\hline
\end{tabular} \end{center} }

\subsection{Outline of the paper}

The rest of this paper is structured as follows. In Section $2$ based on the
relationship between $1$-D OOCs and $2$-D OOCs, many optimal $2$-D $(u\times v,4,2)$-OOCs are derived. Cao and Wei \cite{cws} showed that an optimal $2$-D $(u\times v ,k,t-1)$-OOC is equivalent to an optimal strictly $v$-cyclic $t$-$(u\times
v,k,1)$-packing, provided that $t\leq k$ holds. We restate this combinatorial equivalence in Section $3$. In this section perfect $2$-D OOCs are defined as a special case of optimal $2$-D OOCs. We point out in Remark \ref{approach-perfect} that the problem for the existence of perfect $2$-D $(u\times v,4,2)$-OOCs can be reduced to the problem for the existence of perfect $2$-D $(w\times v,4,2)$-OOCs, $w\in\{1,2\}$. When $w=1$, perfect $2$-D $(1\times v,4,2)$-OOCs have been widely investigated as a kind of combinatorial object called strictly cyclic Steiner quadruple system. Thus we pay our attention to the case of $w=2$ in Section $4$. We give a construction for perfect $2$-D $(2\times v,4,2)$-OOCs. In Section $5$ we improve the upper bound for optimal $2$-D $(u\times v,4,2)$-OOCs (not only focus on perfect), which is tighter than the well-known Johnson bound in many cases. In Sections $6$ and $7$ some auxiliary designs are introduced to establish recursive constructions for $2$-D $(u\times v,k,2)$-OOCs with general $k$. Using these recursive
constructions and some direct constructions, we obtain many infinite
families of optimal $2$-D $(u\times v,4,2)$-OOCs in Sections $8$ and
$9$. Finally, Section $10$ gives a brief conclusion.

Our main results are summarized in Table II (in Section $9$), Tables III and IV (in Section $10$).

\section{$2$-D OOCs from $1$-D OOCs}

$2$-D OOCs are very closely related to $1$-D OOCs. A $2$-D
$(1\times v,k,\lambda)$-OOC is just a $1$-D $(v,k,\lambda)$-OOC. A
$1$-D $(v,k,\lambda)$-OOC  with $\Phi(1\times v,k,\lambda)$
codewords is said to be {\em optimal}. In this section we shall derive some optimal $2$-D OOCs from the known results on optimal $1$-D OOCs. First we quote the following result from Alderson and Mellinger \cite{am1}.

\begin{Theorem}
\label{relation2} {\rm (\cite{am1})} Suppose that there exists an
optimal $2$-D $(u\times v,k,\lambda)$-OOC with $\Phi(u\times
v,k,\lambda)$ codewords. Then for any integer factorization
$v=v_1v_2$, there exists a $2$-D $(uv_1\times v_2,k,\lambda)$-OOC
with $v_1\Phi(u\times v,k,\lambda)$ codewords.
\end{Theorem}

As a corollary of Theorem \ref{relation2}, we have

\begin{Corollary}
\label{relation1} If there is an optimal $1$-D
$(uv,k,\lambda)$-OOC with $\Phi(1\times uv,k,\lambda)$ codewords,
then a $2$-D $(u\times v,k,\lambda)$-OOC with $u\Phi(1\times
uv,k,\lambda)$ codewords exists.
\end{Corollary}

Therefore by Corollary \ref{relation1}, if $u\Phi(1\times
uv,k,\lambda)$ is just equal to $\Phi(u\times v,k,\lambda)$, then
the resulting $2$-D $(u\times v,k,\lambda)$-OOC is optimal. In the following we shall give some analysis on $u\Phi(1\times
uv,k,\lambda)=\Phi(u\times v,k,\lambda)$. According to the Johnson
bound, $\Phi(u\times v,k,\lambda)\leq J(u\times v,k,\lambda)$ and
$\Phi(1\times uv,k,\lambda)\leq J(1\times uv,k,\lambda)$. Assume
that
$$J_1(1\times
uv,k,\lambda)=\frac{1}{k}\lfloor\frac{uv-1}{k-1}
\lfloor\frac{uv-2}{k-2}\lfloor\cdots\lfloor\frac{uv-\lambda}{k-\lambda}
\rfloor\cdots\rfloor\rfloor\rfloor.$$ We have the following lemma.

\begin{Lemma}
\label{iff-1D-2D} $J(u\times v,k,\lambda)=uJ(1\times uv,k,\lambda)$
if and only if $J_1(1\times uv,k,\lambda)-J(1\times
uv,k,\lambda)<1/u$.
\end{Lemma}
\proof Let $x=\lfloor\frac{uv-1}{k-1}
\lfloor\frac{uv-2}{k-2}\lfloor\cdots\lfloor\frac{uv-\lambda}{k-\lambda}
\rfloor\cdots\rfloor\rfloor\rfloor$ and $x=ak+b$, where $0\leq b<k$.
Then $J(u\times v,k,\lambda)=\lfloor\frac{u}{k}x\rfloor$ and
$J(1\times uv,k,\lambda)=\lfloor\frac{1}{k}x\rfloor$. It is easy to
verify that

\begin{tabular}{lllll}
$J(u\times v,k,\lambda)=uJ(1\times uv,k,\lambda)$ &
$\Longleftrightarrow$ &
$\lfloor\frac{u}{k}x\rfloor=u\lfloor\frac{1}{k}x\rfloor$ &
$\Longleftrightarrow$ & $ua+\lfloor \frac{ub}{k}\rfloor=ua$\\
&$\Longleftrightarrow$&$\lfloor
\frac{ub}{k}\rfloor=0$&$\Longleftrightarrow$&$ub<k$\\
&$\Longleftrightarrow$&$\frac{b}{k}<\frac{1}{u}.$
\end{tabular}

\noindent Note that $J_1(1\times uv,k,\lambda)-J(1\times
uv,k,\lambda)=\frac{b}{k}$.  \qed

\begin{Theorem}
\label{2D from 1D} If there exists an optimal $1$-D
$(uv,k,\lambda)$-OOC with $J(1\times uv,k,\lambda)$ codewords and
$J_1(1\times uv,k,\lambda)-J(1\times uv,k,\lambda)<1/u$, then there
exists an optimal $2$-D $(u\times v,k,\lambda)$-OOC with $J(u\times
v,k,\lambda)$ codewords.
\end{Theorem}
\proof By Corollary \ref{relation1}, if there is an optimal $1$-D
$(uv,k,\lambda)$-OOC with $J(1\times uv,k,\lambda)$ codewords, then
there is a $2$-D $(u\times v,k,\lambda)$-OOC with $uJ(1\times
uv,k,\lambda)$ codewords. Since $J_1(1\times uv,k,\lambda)-J(1\times
uv,k,\lambda)<1/u$, by Lemma \ref{iff-1D-2D}, $uJ(1\times
uv,k,\lambda)=J(u\times v,k,\lambda)$. Thus the resulting $2$-D OOC
is optimal. \qed

\begin{Corollary}
\label{2D from 1D k=4} If there exists an optimal $1$-D
$(n,4,2)$-OOC with $J(1\times n,4,2)$ codewords, then
\item{$(1)$} for any integer $n\equiv 1,3\
({\rm mod }\ 6)$ or $n\equiv 2,10\ ({\rm mod }\ 24)$, and for any
integer factorization $n=uv$, there exists an optimal $2$-D
$(u\times v,4,2)$-OOC with $J(u\times v,4,2)$ codewords;
\item{$(2)$} for any integer $n\equiv 4,20\ ({\rm
mod }\ 24)$ and $n=2n_1$, there exists an optimal $2$-D $(2\times
n_1,4,2)$-OOC with $J(2\times n_1,4,2)$ codewords.
\end{Corollary}

\proof When $n\equiv 1,3\ ({\rm mod }\ 6)$ or $n\equiv 2,10\ ({\rm
mod }\ 24)$, it is readily checked that $J_1(1\times
n,4,2)-J(1\times n,4,2)=0$. When $n\equiv 4,20\ ({\rm mod }\ 24)$,
it is readily checked that $J_1(1\times n,4,2)-J(1\times
n,4,2)=1/4$. The assertion then follows from Theorem \ref{2D from
1D}. \qed

As a special topic, $1$-D $(n,4,2)$-OOCs have been extensively
studied, for example Alderson and Mellinger \cite{am2}, Chu and
Colbourn \cite{cc,cc1}, Feng, Chang and Ji \cite{fcj,fcj1}. We only
quote partial known results on optimal $1$-D $(n,4,2)$-OOCs, which
are essential for our work.

\begin{Lemma}
\label{1D k=4}
\item{$(1)$} {\rm (\cite{fcj1})} There exists an
optimal $1$-D $(uv,4,2)$-$OOC$ with $J(1\times uv,4,2)$ codewords
for any $u\in \{4^n-1: {\rm integer} \ n\ge 1\}\cup \{1$, $27$,
$33$, $39$, $51$, $87$, $123$, $183\}$ and $v$  an integer taken
from the set $\{p\equiv 7 \ (mod \ 12): p$ is a prime$\}\cup
\{2^n-1: {\rm odd \ integer} \ n\geq 1\}\cup \{25$, $37$, $61$,
$73$, $109$, $157$, $181$, $229$, $277\}$, or a product of such
integers.

\item{$(2)$}  {\rm (\cite{fcj})} Let $n$ be a positive integer. If $n=p_1^{r_1}p_2^{r_2}\cdots
p_s^{r_s}$, where $p_i=13$ or $p_i$ is a prime, $p_i\equiv 5\ (mod\
12)$ and $p_i<1500000$, $r_i\geq 1$ for $1\leq i\leq s$, then there
is an optimal $1$-D $(4n,4,2)$-OOC with $J(1\times 4n,4,2)$
codewords.

\item{$(3)$} {\rm (\cite{fcj})} There exists an optimal $1$-D $(n,4,2)$-OOC
with $J(1\times n,4,2)$ codewords for all $7\leq n\leq 100$ with the
definite exceptions of $n \in \{9$, $12$, $13$, $24$, $48$, $72$,
$96\}$ and possible exceptions of $n\in \{45$, $47$, $53$, $55$,
$59$, $60$, $65$, $66$, $69$, $71$, $76$, $77$, $81$, $83$, $84$,
$85$, $89$, $91$, $92$, $95$, $97$, $99\}$.

\item{$(4)$} {\rm (\cite{cc,fcj})} There exists an optimal $1$-D $(n,4,2)$-OOC
with $J(1\times n,4,2)-1$ codewords for $n \in \{9$, $12$, $13$,
$24$, $48$, $72$, $96\}$.
\end{Lemma}

\begin{Theorem}
\label{2D k=4 from 1D}
\item{$(1)$} Let $m=uv$, where $u\in \{4^n-1:
\ n\ge 1\}\cup \{1$, $27$, $33$, $39$, $51$, $87$, $123$, $183\}$
and $v$ is an integer taken from the set $\{p\equiv 7 \ (mod \ 12):
p$ is a prime$\}\cup \{2^n-1: {\rm odd \ integer} \ n\geq 1\}\cup
\{25$, $37$, $61$, $73$, $109$, $157$, $181$, $229$, $277\}$, or a
product of such integers. Then for any integer factorization
$m=n_1n_2$, there exists an optimal $2$-D $(n_1\times
n_2,4,2)$-$OOC$ with $J(n_1\times n_2,4,2)$ codewords.

\item{$(2)$} Let $n\in\{10,15,21,25,26,27,33,34,39,49,50,51,57,58,63,74,75,82,87,93,98\}$.
Then for any integer factorization $n=n_1n_2$, there is an optimal
$2$-D $(n_1\times n_2,4,2)$-OOC with $J(n_1\times n_2,4,2)$
codewords.

\item{$(3)$} Let $n$ be a positive integer. If
$n=p_1^{r_1}p_2^{r_2}\cdots p_s^{r_s}$, where $p_i=13$ or $p_i$ is a
prime, $p_i\equiv 5\ (mod\ 12)$ and $p_i<1500000$, $r_i\geq 1$ for
$1\leq i\leq s$, then there is an optimal $2$-D $(2\times
2n,4,2)$-OOC with $J(2\times 2n,4,2)$ codewords.

\item{$(4)$} Let $2n\in\{20,28,44,52,68,100\}$.
Then there is an optimal $2$-D $(2\times n,4,2)$-OOC with $J(2\times
n,4,2)$ codewords.
\end{Theorem}
\proof It is readily checked that the number $n$ described in $(1)$
is congruent to $1$ or $3$ modulo $6$; the number $n$ described in
$(2)$ is congruent to $1,3$ modulo $6$, or $2,10$ modulo $24$; the
number $4n$ described in $(3)$ and the number $2n$ described in
$(4)$ are congruent to $4$ or $20$ modulo $24$. Combining the
results of Lemma \ref{1D k=4}, the assertion then follows from
Corollary \ref{2D from 1D k=4}. \qed

\section{Combinatorial descriptions}

Two-dimensional optical orthogonal codes are closely related to some
combinatorial configurations called strictly $v$-cyclic packings.
Throughout this paper we always assume that $I_u=\{0,1,\ldots,u-1\}$
and denote by $Z_v$ the additive group of integers modulo $v$.

\subsection{Combinatorial equivalence}

A $t$-$(v,k,1)$ {\em packing} is a pair $(X,{\cal B})$, where $X$ is
a set of $v$ points and ${\cal B}$ is a set of subsets of $X$
(called {\em blocks}), each of cardinality $k$, such that every
$t$-subset of $X$ occurs in at most one block. The set of all
uncovered $t$-subsets by $\cal B$ is said to be the {\em leave} of
the packing.

An {\em automorphism} $\alpha$ of a packing $(X,{\cal B})$ is a permutation
on $X$ such that
$$\{\{\alpha(x):x\in B\}:B\in{\cal B}\}={\cal B}.$$ In other words, a block of the packing is mapped to a block under an automorphism. A $t$-$(u\times v,k,1)$ packing
is said to be {\em $v$-cyclic} if it admits an automorphism $\pi$
consisting of $u$ cycles of length $v$. Without loss of generality
identify $X$ with $I_u\times Z_v$ and the automorphism $\pi$ can be
taken as $(i,x)\longmapsto(i,x+1)$ (mod $(-,v)$), $i\in I_u$ and
$x\in Z_v$. Then all blocks of this packing can be partitioned into
some orbits under $\pi$. Choose any fixed block from each orbit and
then call it a {\em base block}.

All automorphisms of a packing form a group, called the {\em full
automorphism group} of the packing. Any subgroup of the full
automorphism group is called an {\em automorphism group} of the
packing. Let $G$ be an automorphism group of a packing. For any
block $B$ of the packing, the subgroup $$\{\pi\in G: B^{\pi}=B\}$$ is
called the {\em stabilizer} of $B$ in $G$. If the stabilizer of each
block of a $v$-cyclic $t$-$(u\times v,k,1)$ packing is trivial in
$Z_v$, i.e., for each block $B$, $\{\delta\in Z_v: B+\delta=B\}=\{0\}$, where
$B+\delta=\{(i,x+\delta):(i,x)\in B\}$, then the packing is called {\em strictly $v$-cyclic}. When
$u=1$, a (strictly) $v$-cyclic $t$-$(1\times v,k,1)$ packing is
often simply referred to as a {\em $($strictly$)$ cyclic
$t$-$(v,k,1)$ packing}. When $v=1$, a (strictly) $1$-cyclic
$t$-$(u\times 1,k,1)$ packing is just a $t$-$(u,k,1)$ packing.

A strictly $v$-cyclic $t$-$(u\times v,k,1)$ packing is called
$optimal$ if it contains the largest possible block number. The main purpose of this paper is to construct optimal $2$-D OOCs. Cao and Wei \cite{cws} established the equivalence between optimal $2$-D OOCs and optimal strictly $v$-cyclic packings. Suppose $(X,{\cal B})$ is a strictly $v$-cyclic $t$-$(u\times
v,k,1)$ packing. Denote the family of base blocks of this packing by
$\cal F$. For each base block $B$ of $\cal F$, construct an $u\times
v$ $(0,1)$-matrix $M_B$ whose rows are indexed by $I_u$ and columns
are indexed by $Z_v$, such that its $(i,j)$ cell equals $1$ if and
only if $(i,j)\in B$. Since any two blocks intersect at most $t-1$
points and all the blocks can be generated by developing cyclically
the base blocks, $\{M_B:B\in {\cal F}\}$ forms a $2$-D $(u\times
v,k,t-1)$-OOC with $|{\cal F}|$ codewords. Conversely, given a $2$-D
$(u\times v,k,t-1)$-OOC, $\cal C$, for each $u\times v$
$(0,1)$-matrix $M\in\cal C$ whose rows are indexed by $I_u$ and
columns are indexed by $Z_v$, construct a $k$-subset $B_M$ of
$I_u\times Z_v$ such that $(i,j)\in B_M$ if and only if $M$'s
$(i,j)$ cell equals $1$. Then $\{B_M:M\in {\cal C}\}$ is the family
of base blocks of a strictly $v$-cyclic $t$-$(u\times v,k,1)$
packing.

\begin{Theorem}
\label{cws-optimal} {\rm (\cite{cws})} An optimal $2$-D $(u\times v ,k,t-1)$-OOC is
equivalent to an optimal strictly $v$-cyclic $t$-$(u\times
v,k,1)$-packing, provided that $t\leq k$ holds.
\end{Theorem}

Since a strictly $1$-cyclic $3$-$(u\times 1,4,1)$ packing is just a
$3$-$(u,4,1)$ packing, and the existence of an optimal $3$-$(u,4,1)$
packing has been investigated by Ji \cite{jip}, we can have the following result.

\begin{Theorem}
\label{ji-v=1} {\rm (\cite{jip})} There exists an optimal $2$-D
$(u\times 1,4,2)$-OOC $($i.e., an optimal $3$-$(u,4,1)$ packing$)$
with $\phi$ codewords, where
$$ \phi=\left\{
\begin{array}{lll}
\lfloor\frac{u}{4}\lfloor\frac{u-1}{3}\lfloor\frac{u-2}{2}
\rfloor\rfloor\rfloor & \hbox{ $u\not\equiv 0$ $({\rm mod}$ $6)$,} \\
    \\
\lfloor\frac{u}{4}(\lfloor\frac{u-1}{3}\lfloor\frac{u-2}{2}
\rfloor\rfloor-1)\rfloor & \hbox{ $u\equiv 0$ $({\rm mod}$ $6)$,} \\
\end{array}
\right.$$ with the exception of $21$ undecided values $u=6r+5$,
$r\in\{m:m$ is odd, $3\leq m\leq 35$, $m\neq 17,21\}\cup$
$\{45,47,75,77,79,159\}$.
\end{Theorem}

\begin{Example}
\label{example-6} There is a trivial optimal $2$-D $(1\times
6,4,2)$-OOC, whose number of codewords is
$\lfloor\frac{1}{4}\lfloor\frac{5}{3}\lfloor\frac{4}{2}
\rfloor\rfloor\rfloor=0$. By Theorem $\ref{ji-v=1}$, there is an
optimal $2$-D $(6\times 1,4,2)$-OOC with $3$ codewords. An optimal
$2$-D $(2\times 3,4,2)$-OOC has only $J(2\times 3,4,2)=1$ codeword
\begin{center}
$\left(
  \begin{array}{ccc}
    1 & 1 & 0 \\
    1 & 1 & 0 \\
  \end{array}
\right),$
\end{center}
\noindent whose corresponding base block of the optimal strictly
$3$-cyclic $3$-$(2\times 3,4,1)$-packing is
$\{(0,0),(1,0),(0,1),(1,1)\}$. A $2$-D $(3\times 2,4,2)$-OOC can not
contain $J(3\times 2,4,2)=2$ codewords. Otherwise, there were a
$2$-D $(6\times 1,4,2)$-OOC with $4$ codewords by Theorem
$\ref{relation2}$, which would be contradict to Theorem
$\ref{ji-v=1}$. Thus an optimal $2$-D $(3\times 2,4,2)$-OOC has only
one codeword
\begin{center}
$\left(
  \begin{array}{cc}
    1 & 1 \\
    1 & 0 \\
    1 & 0 \\
  \end{array}
\right),$
\end{center}
\noindent whose corresponding base block of the optimal strictly
$2$-cyclic $3$-$(3\times 2,4,1)$-packing is
$\{(0,0),(1,0),(2,0),(0,1)\}$.
\end{Example}

\subsection{Perfect $2$-D OOCs}

Let $K$ be a set of positive integers. A {\em $t$-wise balanced
design} (briefly $t$-design) is a pair $(X,{\cal B})$, where $X$ is
a set of $v$ points and ${\cal B}$ is a set of subsets of $X$
(called {\em blocks}), each of cardinality from $K$, such that every
$t$-subset of $X$ is contained in a unique block. Such a design is
denoted by $S(t,K,v)$. If $K=\{k\}$, we write $S(t,K,v)$ by
$S(t,k,v)$. An $S(2,3,v)$ is called a {\em Steiner triple system}
and denoted by $STS(v)$. An $S(3,4,v)$ is called a {\em Steiner
quadruple system} and denoted by $SQS(v)$.

Evidently an $S(t,k,v)$ is a special $t$-$(v,k,1)$-packing, whose
leave is an empty set. Thus one can similarly define (strictly)
$v$-cyclic $S(t,k,u\times v)$ as we have done for (strictly)
$v$-cyclic $t$-$(u\times v,k,1)$-packing. A strictly $v$-cyclic
$SQS(1\times v)$ is often simply written as an $sSQS(v)$ (cf.
\cite{fcj}).

If a $2$-D $(u\times v,k,t-1)$-OOC is equivalent to a strictly
$v$-cyclic $S(t,k,u\times v)$, then the OOC is said to be {\em
perfect}. It is easy to verify that a perfect OOC is an optimal OOC
that attains the Johnson bound without using the brackets (cf.
\cite{omk}).

\begin{Lemma}
\label{nece-perfect} The necessary conditions for the existence of a
strictly $v$-cyclic $SQS(u\times v)$ $($or equivalently, a perfect
$2$-D $(u\times v,4,2)$-OOC$)$ are $uv\equiv 2,4$ $({\rm mod}$ $6)$,
$u(uv-1)(uv-2)\equiv 0$ $({\rm mod}$ $24)$. Specifically, the
necessary conditions can be classified as follows:

\begin{enumerate}
\item[{\rm (1)}] $u\equiv 1,5$ $({\rm mod }$
$12)$ and $v\equiv 2,10$ $({\rm mod }$ $24)$;

\item[{\rm (2)}] $u\equiv 7,11$ $({\rm mod }$
$12)$ and $v\equiv 14,22$ $({\rm mod }$ $24)$;

\item[{\rm (3)}] $u\equiv 2,4$ $({\rm mod }$ $6)$ and $v\equiv 1,5$
$({\rm mod }$ $6)$;

\item[{\rm (4)}] $u\equiv 4,8$ $({\rm mod }$
$12)$ and $v\equiv 2,4$ $({\rm mod }$ $6)$.
\end{enumerate}
\end{Lemma}

\proof It is well known that an $SQS(uv)$ exists if and only if
$uv\equiv 2,4$ $({\rm mod}$ $6)$ \cite{hanani60}. Count the number
of base blocks of a strictly $v$-cyclic $SQS(u\times v)$. It follows
that $u(uv-1)(uv-2)\equiv 0$ $({\rm mod}$ $24)$. \qed

A natural question from Lemma \ref{nece-perfect} is whether the
necessary conditions for the existence of a perfect $2$-D $(u\times
v,4,2)$-OOC are sufficient. In Section $5$, by Corollary~\ref{cor no
perfect}, we show that for $u\equiv 4,8$ $({\rm mod }$ $12)$ and
$v\equiv 2,4$ $({\rm mod }$ $6)$, there is no perfect $2$-D
$(u\times v,4,2)$-OOC. In Section $9$, by Proposition \ref{infinite
family-u-v-1}, if there exists a perfect $2$-D $(2\times v,4,2)$-OOC
with $v\equiv 1,5$ $({\rm mod }$ $6)$, then a perfect $2$-D
$(u\times v,4,2)$-OOC exists for any $u\equiv 2,4$ $({\rm mod}$
$6)$. When $u$ and $v$ satisfy Conditions $(1)$ and $(2)$ in Lemma
\ref{nece-perfect}, $uv\equiv 2,10$ $({\rm mod }$ $24)$. Then by
Corollary $\ref{2D from 1D k=4}(1)$, if there exists an optimal
$1$-D $(uv,4,2)$-OOC with $J(1\times uv,4,2)$ codewords, a perfect
$2$-D $(u\times v,4,2)$-OOC exists. Note that when $uv\equiv 2,10$
$({\rm mod }$ $24)$, an optimal $1$-D $(uv,4,2)$-OOC with $J(1\times
uv,4,2)$ codewords is just a perfect $2$-D $(1\times uv,4,2)$-OOC.
Thus

\begin{Remark}
\label{approach-perfect}
The existence problem of perfect $2$-D $(u\times v,4,2)$-OOCs can be reduced to the existence problems of perfect $2$-D $(1\times v,4,2)$-OOCs and perfect $2$-D $(2\times v,4,2)$-OOCs.
\end{Remark}

\section{A construction for perfect $2$-D $(2\times v,4,2)$-OOCs}

According to Remark \ref{approach-perfect}, it is important to consider the existences of perfect $2$-D $(1\times v,4,2)$-OOCs and perfect $2$-D $(2\times v,4,2)$-OOCs. A perfect $2$-D $(1\times v,4,2)$-OOC is equivalent to an $sSQS(v)$. Much work has been done on $sSQS$s in the literature. The interested reader may refer to \cite{fcj} and the references therein. In this section, we shall present a construction for perfect $2$-D $(2\times v,4,2)$-OOCs.

The idea of this construction is originally from Hartman \cite{hartman}. In $1980$ Hartman \cite{hartman} gave a construction for an $SQS(2p)$, which can be obtained from an $SQS(p+1)$ with a cyclic derived Steiner triple system, where $p\equiv 1$ $({\rm mod}$ $6)$ is a prime. Here, we shall generalize Hartman's method to obtain a construction for strictly $p$-cyclic $SQS(2\times p)$s. The existence of a strictly $p$-cyclic
$SQS(2\times p)$ implies the existence of a perfect $2$-D $(2\times
p,4,2)$-OOC.

Our construction are based on the concept of rotational $SQSs$. A {\em rotational $SQS(n)$} is an $SQS(n)$ with an automorphism
consisting of one fixed point and a cycle of length $n-1$. Such a
design is denoted by $RoSQS(n)$.

Assume that $(X,\cal B)$ is an $RoSQS(n)$. We can identify $X$ with $Z_{n-1}\cup \{\infty\}$, and let the permutation $\alpha$ fixing $\infty$ and mapping $i$ to $i+1$ (mod $n-1$), $i\in Z_{n-1}$, be an automorphism of the $RoSQS$. Let $G$ be a cyclic group
generated by $\alpha$ under the compositions of permutations. Then
all blocks of the $RoSQS$ can be partitioned into some orbits
under $G$. Choose any fixed block from each orbit and then call it
a {\em base block}.

\begin{Example}
\label{RoSQS-8} An $RoSQS(8)$ $(X,{\cal B})$ is constructed on $X=Z_7\cup\{\infty\}$. All blocks of $\cal B$ are listed below:
\begin{center}
$\{i,i+1,i+2,i+5\}$, $\{i,i+1,i+3,\infty\}$, $0\leq i\leq 6$.
\end{center}
\noindent Obviously, all blocks of $\cal B$ can be obtained by developing the two base blocks $\{0,1,2,5\}$, $\{0,1,3,\infty\}$ by $+1$ modulo $7$, where $\infty+1=\infty$.
\end{Example}

\begin{Construction}
\label{RoSQS-recur} Let $p\equiv 1$ $({\rm mod}$ $6)$ be a prime. If
there exists an $RoSQS(p+1)$, then there exists a strictly
$p$-cyclic $SQS(2\times p)$.
\end{Construction}

\proof Here we only exhibit the algorithm in Figure $1$. The detailed proof of this construction has been moved to Appendix I. \qed
 {\footnotesize
\begin{algorithm}[H]
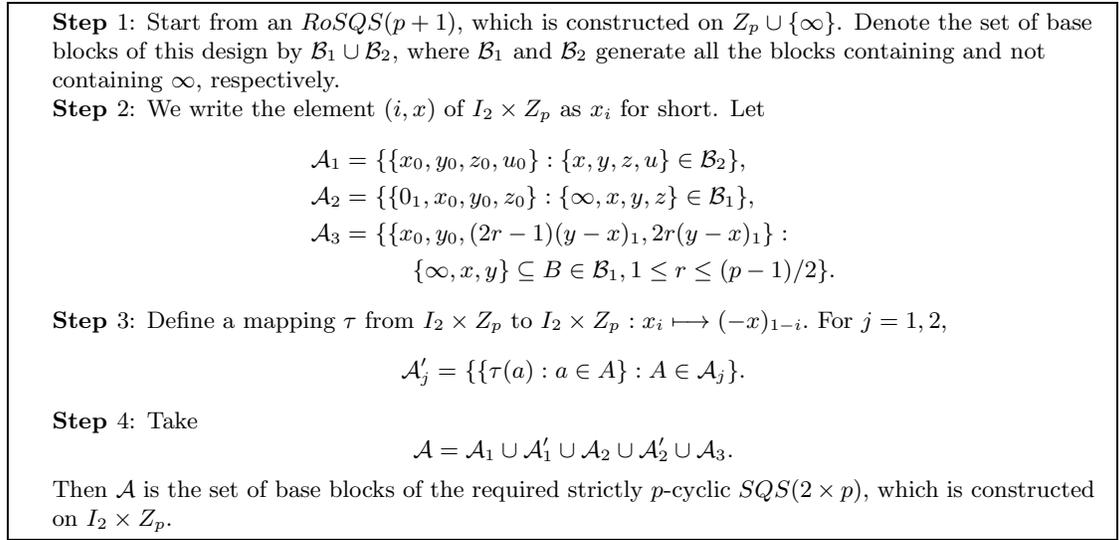
\caption{Algorithm to construct a strictly $p$-cyclic $SQS(2\times p)$}

{\rm \textbf{Step $1$}}: Start from an $RoSQS(p+1)$, which is constructed on $Z_p\cup
\{\infty\}$. Denote the set of base blocks of this design by ${\cal
B}_1\cup {\cal B}_2$, where ${\cal B}_1$ and ${\cal B}_2$ generate
all the blocks containing and not containing $\infty$, respectively.

{\rm \textbf{Step $2$}}: We write the
element $(i,x)$ of $I_2\times Z_p$ as $x_i$ for short. Let
\begin{eqnarray*}
{\cal A}_1=\hspace{-3mm}&\hspace{-3mm}&\{\{x_0,y_0,z_0,u_0\}:\{x,y,z,u\}\in {\cal B}_2\},\\
{\cal A}_2=\hspace{-3mm}&\hspace{-3mm}&\{\{0_1,x_0,y_0,z_0\}:\{\infty,x,y,z\}\in {\cal B}_1\},\\
{\cal A}_3=\hspace{-3mm}&\hspace{-3mm}&
\{\{x_0,y_0,(2r-1)(y-x)_1,2r(y-x)_1\}:\\
\hspace{-3mm}&\hspace{-3mm}&\hspace{0.5cm}
\{\infty,x,y\}\subseteq B\in {\cal B}_1,1\leq r\leq (p-1)/2\}.
\end{eqnarray*}

{\rm \textbf{Step $3$}}: Define a mapping $\tau$ from $I_2\times Z_p$ to $I_2\times Z_p$ $:  x_i \longmapsto (-x)_{1-i}.$ For $j=1,2$,
\begin{eqnarray*}
{\cal A}'_j=\hspace{-3mm}&\hspace{-3mm}&\{\{\tau(a):a\in A\}:{A\in
{\cal A}_j}\}.
\end{eqnarray*}

{\rm \textbf{Step $4$}}: Take $${\cal A}={\cal A}_1\cup{\cal A}'_1\cup{\cal A}_2\cup{\cal A}'_2\cup{\cal A}_3.$$ Then $\cal A$ is the set of base blocks of the required strictly $p$-cyclic $SQS(2\times p)$, which is constructed on $I_2\times Z_p$.
\end{algorithm}}

\vspace{0.5cm}

The following example illustrates the algorithm presented in Figure $1$.

\begin{Example}
\label{RoSQS-2*7} In this example we shall show how to construct a strictly
$7$-cyclic $SQS(2\times 7)$ from an $RoSQS(8)$. It is equivalent to a perfect $2$-D $(2\times 7,4,2)$-OOC by Theorem \ref{cws-optimal}.

\begin{itemize}
\item Step $1$: Start from an $RoSQS(8)$, which is given by Example $\ref{RoSQS-8}$. Take {\footnotesize$${\cal B}_1=\{\{0,1,3,\infty\}\},\hspace{0.5cm}{\cal B}_2=\{\{0,1,2,5\}\}.$$}
\item Step $2$: Construct the required strictly
$7$-cyclic $SQS(2\times 7)$ on $I_2\times Z_7$. Let
{\footnotesize\begin{eqnarray*}
{\cal A}_1=\hspace{-3mm}&\hspace{-3mm}&\{\{0_0,1_0,2_0,5_0\}\},\hspace{10mm}
{\cal A}_2=\{\{0_1,0_0,1_0,3_0\}\},\\
{\cal A}_3=\hspace{-3mm}&\hspace{-3mm}&
\{\{0_0,1_0,1_1,2_1\}\cup\{0_0,1_0,3_1,4_1\}\cup\{0_0,1_0,5_1,6_1\}\\
\hspace{-3mm}&\hspace{-3mm}&
\cup\{0_0,3_0,3_1,6_1\}\cup\{0_0,3_0,2_1,5_1\}\cup\{0_0,3_0,1_1,4_1\}\\
\hspace{-3mm}&\hspace{-3mm}&
\cup\{1_0,3_0,2_1,4_1\}\cup\{1_0,3_0,6_1,1_1\}\cup\{1_0,3_0,3_1,5_1\}.
\end{eqnarray*}}
\item Step $3$: Under the action of the mapping $\tau:  x_i \longmapsto (-x)_{1-i}$, we have {\footnotesize$$ {\cal A}'_1=\{\{0_1,6_1,5_1,2_1\}\},\hspace{0.5cm}{\cal A}'_2=\{\{0_0,0_1,6_1,4_1\}\}.$$}
\item Step $4$: Let ${\cal A}={\cal A}_1\cup{\cal A}'_1\cup{\cal A}_2\cup{\cal A}'_2\cup{\cal A}_3$. Then $|{\cal A}|=13$ and $\cal A$ is the set of base blocks of the required strictly $7$-cyclic $SQS(2\times 7)$.
\end{itemize}
\end{Example}

By Theorem \ref{cws-optimal}, a strictly $p$-cyclic $SQS(2\times p)$ is equivalent to a perfect $2$-D $(2\times p,4,2)$-OOC. Thus by Construction \ref{RoSQS-recur}, for obtaining some perfect $2$-D $(2\times p,4,2)$-OOCs, we need some results on $RoSQSs$.

\begin{Theorem}
\label{RoSQS-update} {\rm (\cite{fcj1})} There exists an
$RoSQS(uv+1)$ for any $u\in \{4^n-1: {\rm integer} \ n\ge 1\}\cup
\{1$, $27$, $33$, $39$, $51$, $87$, $123$, $183\}$ and $v$ is an
integer taken from the set $\{p\equiv 7 \ (mod \ 12): p$ is a
prime$\}\cup \{2^n-1: {\rm odd \ integer} \ n\geq 1\}\cup \{25$,
$37$, $61$, $73$, $109$, $157$, $181$, $229$, $277\}$, or a product
of such integers.
\end{Theorem}

Combining the results of Theorem \ref{cws-optimal}, Construction
\ref{RoSQS-recur} and Theorem \ref{RoSQS-update}, we have

\begin{Theorem}
\label{RoSQS-result}  There exist a strictly $p$-cyclic $SQS(2\times
p)$ and a perfect $2$-D $(2\times p,4,2)$-OOC for any prime $p\equiv
7 \ (mod \ 12)$ or $p\in \{37$, $61$, $73$, $109$, $157$, $181$,
$229$, $277\}$.
\end{Theorem}

\section{Tighter upper bound for $2$-D $(u\times v,4,2)$-OOCs}

In most cases an optimal $2$-D $(u\times v,4,2)$-OOC is not a perfect $2$-D $(u\times v,4,2)$-OOC. Thus the determination of the largest possible size $\Phi(u\times v,4,2)$ of a optimal $2$-D $(u\times v,4,2)$-OOC is of interest. Recall that in Section $1$, we mention that $\Phi(u\times v,4,2)\leq J(u\times v,4,2)$, where $J(u\times v,4,2)=\lfloor\frac{u}{4}\lfloor\frac{uv-1}{3}
\lfloor\frac{uv-2}{2}\rfloor\rfloor\rfloor$ is the famous Johnson bound.
Here we shall give a tighter upper bound for $2$-D $(u\times v,4,2)$-OOCs than the Johnson bound.

\begin{Lemma}
\label{bound-6 mod 0} Let $uv\equiv 0$ $({\rm mod}$ $6)$. Then
$\Phi(u\times v,4,2)\leq
\lfloor\frac{u}{4}(\lfloor\frac{uv-1}{3}\lfloor\frac{uv-2}{2}
\rfloor\rfloor-1)\rfloor$.
\end{Lemma}

\proof By Theorem \ref{relation2}, an optimal $2$-D $(u\times
v,4,2)$-OOC with $\Phi(u\times v,4,2)$ codewords implies a $2$-D
$(uv\times 1,4,2)$-OOC with $v\Phi(u\times v,4,2)$ codewords. Since
a $2$-D $(uv\times 1,4,2)$-OOC is equivalent to a strictly
$1$-cyclic $3$-$(uv\times 1,4,1)$-packing, by Theorem \ref{ji-v=1},
when $uv\equiv 0$ $({\rm mod}$ $6)$, it has at most
$\lfloor\frac{uv}{4}(\lfloor\frac{uv-1}{3}\lfloor\frac{uv-2}{2}
\rfloor\rfloor-1)\rfloor$ blocks. Thus we have $v\Phi(u\times
v,4,2)\leq
\lfloor\frac{uv}{4}(\lfloor\frac{uv-1}{3}\lfloor\frac{uv-2}{2}
\rfloor\rfloor-1)\rfloor$. It is readily checked that $\Phi(u\times
v,4,2)\leq\lfloor\frac{1}{v}\lfloor\frac{uv}{4}(\lfloor\frac{uv-1}{3}\lfloor\frac{uv-2}{2}
\rfloor\rfloor-1)\rfloor\rfloor=\lfloor\frac{u}{24}(u^2v^2-3uv-6)\rfloor=\lfloor\frac{u}{4}(\lfloor\frac{uv-1}{3}\lfloor\frac{uv-2}{2}
\rfloor\rfloor-1)\rfloor$.  \qed

The following two lemmas shows that in some cases of $uv\equiv 0$ $({\rm mod}$ $6)$, the bound for $\Phi(u\times v,4,2)$ in Lemma \ref{bound-6 mod 0} is not tight enough. Their proofs are lengthy. To ensure smooth reading of the paper, their proofs have been moved to Appendix II.

\begin{Lemma}
\label{bound-u-v-12-2} Let $u\equiv 0$ $({\rm mod}$ $12)$ and
$v\equiv 2,4$ $({\rm mod }$ $6)$. Then $\Phi(u\times v,4,2)\leq
\lfloor\frac{u}{4}(\lfloor\frac{uv-1}{3}\lfloor\frac{uv-2}{2}
\rfloor\rfloor-1)\rfloor-1$.
\end{Lemma}

\begin{Lemma}
\label{bound-u-v-4-6} Let $uv\equiv 0$ $({\rm mod}$ $12)$ and
$v\equiv 0$ $({\rm mod}$ $6)$. Then $\Phi(u\times v,4,2)\leq
\lfloor\frac{u}{4}(\lfloor\frac{uv-1}{3}$ $\lfloor\frac{uv-2}{2}
\rfloor\rfloor-2)\rfloor$.
\end{Lemma}

\begin{Lemma}
\label{bound-u-v-even} Let $uv\equiv 4,8$ $({\rm mod}$ $12)$ and
$v\equiv 0$ $({\rm mod}$ $2)$. Then $\Phi(u\times v,4,2)\leq
\lfloor\frac{u}{4}(\lfloor\frac{uv-1}{3}\lfloor\frac{uv-2}{2}
\rfloor\rfloor-1)\rfloor$.
\end{Lemma}

\proof For each $a\in I_u$ and each $0\leq i<v/2$, consider the
number $n$ of the base blocks containing the two points $(a,i)$,
$(a,v/2+i)$ in a strictly $v$-cyclic $3$-$(u\times v,4,1)$-packing.
Since each $3$-subset of $I_u\times Z_v$ occurs in at most one block
and each base block containing the two points $(a,i)$, $(a,v/2+i)$
generates exactly two different blocks containing the same two
points, the number $n$ is at most $\lfloor(uv-2)/4\rfloor=(uv-4)/4$.
Thus there are at least two $3$-subsets of the form
$\{(a,i),(a,v/2+i),(*,*)\}$ in the leave. Note that the above
conclusion holds for each $a\in I_u$ and each $0\leq i<v/2$. It
follows that there are at least $uv$ $3$-subsets in the leave. It
implies that $\Phi(u\times v,4,2)\leq
\lfloor({{uv}\choose{3}}-uv)/(4v)\rfloor=\lfloor\frac{1}{24}u(u^2v^2-3uv-4)\rfloor.$
It is readily checked that
$\lfloor\frac{u}{4}(\lfloor\frac{uv-1}{3}\lfloor\frac{uv-2}{2}
\rfloor\rfloor-1)\rfloor=\lfloor\frac{1}{24}u(u^2v^2-3uv-4)\rfloor.$
This completes the proof. \qed

\begin{Corollary}
\label{cor no perfect} For any $u\equiv 4,8$ $({\rm mod }$ $12)$ and
$v\equiv 2,4$ $({\rm mod }$ $6)$, there is no perfect $2$-D
$(u\times v,4,2)$-OOC.
\end{Corollary}

\proof If there were a perfect $2$-D $(u\times v,4,2)$-OOC for
$u\equiv 4,8$ $({\rm mod }$ $12)$ and $v\equiv 2,4$ $({\rm mod }$
$6)$, then it should have $u(uv-1)(uv-2)/24$ codewords. By Lemma
\ref{bound-u-v-even}, the largest possible size of the perfect $2$-D
$(u\times v,4,2)$-OOC should be $u(uv+1)(uv-4)/24$. A contradiction
occurs. \qed

\begin{Lemma}
\label{bound-u-7-11} Let $u\equiv 7,11$ $({\rm mod}$ $12)$. Then
$\Phi(u\times 2,4,2)\leq
\lfloor\frac{u}{4}\lfloor\frac{2u-1}{3}\lfloor\frac{2u-2}{2}
\rfloor\rfloor\rfloor-1$.
\end{Lemma}

\proof It is known that $\Phi(u\times 2,4,2)\leq J(u\times
2,4,2)=\lfloor\frac{u}{4}\lfloor\frac{2u-1}{3}\lfloor\frac{2u-2}{2}
\rfloor\rfloor\rfloor$. Suppose that $\Phi(u\times 2,4,2)=J(u\times
2,4,2)$. Then there were a strictly $2$-cyclic $3$-$(u\times
2,4,1)$-packing with $J(u\times 2,4,2)$ base blocks. Count the
number of $3$-subsets in the leave $\cal L$ of the strictly
$2$-cyclic $3$-$(u\times 2,4,1)$-packing. It is ${{2u}\choose{3}}-
J(u\times 2,4,2)\cdot2\cdot4=4$. Each $3$-subset in the leave is of
the form $\{(a,i),(b,j),(c,k)\}$ or $\{(a,i),(b,j),(a,i+1)\}$, where
$a,b,c$ are distinct elements in $I_u$, and $i,j,k\in Z_2$.

Assume that $\{(a,i),(b,j),(x,k)\}$ is a $3$-subset in the leave,
where $a,b,x \in I_u$, $a\neq b$ and $i,j,k\in Z_2$. Consider the
number $n$ of the blocks containing the two points $(a,i)$, $(b,j)$.
Since each $3$-subset of $I_u\times Z_2$ occurs in at most one
block, the number $n$ is at most $\lfloor(2u-3)/2\rfloor=(2u-4)/2$.
Thus there must be another $3$-subset $\{(a,i),(b,j),(y,l)\}$ in the
leave, where $(y,l)\neq (x,k)$. Due to $|{\cal L}|=4$, we have
${\cal
L}=\{\{(a,i),(b,j),(x,k)\},\{(a,i+1),(b,j+1),(x,k+1)\},\{(a,i),(b,j),(y,l)\},\{(a,i+1),(b,j+1),(y,l+1)\}\}$.

If $x\neq a$ and $x\neq b$, since each $3$-subset of $I_u\times Z_2$
occurs in at most one block, the number of blocks containing the two
points $(a,i)$, $(x,k)$ is exactly $(2u-3)/2$, which is not an
integer. A contradiction. If $x=a$, then $(x,k)=(a,i+1)$ and there
are exactly $(2u-4)/4$ base blocks containing the points $(a,i)$,
$(x,k)$. If $x=b$, then $(x,k)=(b,j+1)$ and there are also exactly
$(2u-4)/4$ base blocks containing the points $(b,j)$, $(x,k)$. The
number $(2u-4)/4$ is not an integer. A contradiction. Hence $|{\cal
L}|\neq 4$ and $\Phi(u\times 2,4,2)\leq J(u\times 2,4,2)-1$. \qed

Combine the results of Lemmas \ref{bound-6 mod
0}-\ref{bound-u-7-11}. Let $A=\{(u,v):u\equiv 0$ $({\rm mod}$ $12)$, $v\equiv 2,4$ $({\rm
mod}$ $6)\}$ and $B=\{(u,v):uv\equiv 0$ $({\rm mod}$ $12)$, $v\equiv
0$ $({\rm mod}$ $6)\}$. In the rest of this paper, we always assume
that

$$ J^*(u\times v)=\left\{
\begin{array}{lll}
\lfloor\frac{u}{4}\lfloor\frac{2u-1}{3}\lfloor\frac{2u-2}{2}
\rfloor\rfloor\rfloor-1, & \hbox{if\ $u\equiv 7,11$ $({\rm mod}$
$12)$ and $v=2$};   \\[1ex]
    \lfloor\frac{u}{4}(\lfloor\frac{uv-1}{3}\lfloor\frac{uv-2}{2}
\rfloor\rfloor-1)\rfloor, & \hbox{if\ $uv\equiv 0$ $({\rm mod}$ $6)$
and $(u,v)\not\in A\cup B$}, \\ & or\ \hbox{$uv\equiv
4,8$ $({\rm mod}$ $12)$ and $v\equiv 0$ $({\rm mod}$ $2)$};
    \\[1ex]
    \lfloor\frac{u}{4}(\lfloor\frac{uv-1}{3}\lfloor\frac{uv-2}{2}
\rfloor\rfloor-1)\rfloor-1, & \hbox{if\ $(u,v)\in A$};
  \\[1ex]
\lfloor\frac{u}{4}(\lfloor\frac{uv-1}{3}\lfloor\frac{uv-2}{2}
\rfloor\rfloor-2)\rfloor, & \hbox{if\ $(u,v)\in B$};   \\[1ex]
\lfloor\frac{u}{4}\lfloor\frac{uv-1}{3}\lfloor\frac{uv-2}{2}
\rfloor\rfloor\rfloor, & \hbox{otherwise}.
\end{array}
\right.$$ We have the following theorem.

\begin{Theorem}
\label{bound-total} $\Phi(u\times v,4,2)\leq J^*(u\times v)$.
\end{Theorem}

Now the question is whether there are optimal $2$-D
$(u\times v,4,2)$-OOCs to achieve the upper bounds established in Theorem
\ref{bound-total}. In Section $9$ we shall give many infinite families
for optimal $2$-D $(u\times v,4,2)$-OOCs, which achieve the upper bound in Theorem \ref{bound-total}.

\section{Auxiliary designs and filling constructions}

In this section and the next section, some recursive constructions for optimal $2$-D OOCs will be given, called filling constructions and weighting constructions, respectively. These constructions are the generalization of standard constructions for $3$-designs in combinatorial design theory. So far the research on combinatorial constructions for $2$-D OOCs mainly focuses on $\lambda=1$ \cite{cws,wsy}, which corresponds to the theory of $2$-designs. However, when $\lambda=2$, the research is related to the theory of $3$-designs. Compared to $2$-designs, the known results on $3$-designs are limited, and the auxiliary structures to construct $3$-designs are more complex. Thus the following auxiliary designs will be a little strange for the reader who first meets them. If one is familiar with $2$-designs, it is useful to notice that the concepts of $s$-fan designs and $H$ designs are two possible generalizations of group divisible designs. Group divisible design is one of the most basic research objects in combinatorial design theory \cite{bjl}.

\subsection{$\bm s$-fan designs}

Hartman \cite{hartman94} first introduced the concept of $s$-fan designs in $1994$. Let $s$ be a non-negative integer. An {\em $s$-fan design} is an
$(s+3)$-tuple $(X,{\cal G},{\cal B}_1,{\cal B}_2,\ldots,{\cal
B}_s,{\cal T})$ satisfying that $(X,{\cal G})$ is a $1$-design,
$(X,{\cal G}\cup {\cal B}_i)$ is a $2$-design for each $1\leq i\leq
s$ and $(X,{\cal G}\cup (\bigcup_{i=1}^s {\cal B}_i) \cup {\cal T})$
is a $3$-design. The elements of ${\cal G}$ and $(\bigcup_{i=1}^s
{\cal B}_i)\cup {\cal T}$ are called {\em groups} and {\em blocks},
respectively.

For understanding the concept of $s$-fan designs, we first consider the case of $s=0$. A $0$-fan design is a $3$-tuple $(X,{\cal G},{\cal T})$ satisfying that $(X,{\cal G})$ is a $1$-design and $(X,{\cal G}\cup {\cal T})$ is a $3$-design.

\begin{Example}
\label{0-fan(4^2)} Take $X=I_8$ and ${\cal G}=\{\{0,2,4,6\},\{1,3,5,7\}\}$. Then $(X,{\cal G})$ is a $1$-design. Let $\cal T$ consists of the following $12$ blocks
{\footnotesize
\begin{center}
\begin{tabular}{llllll}
$\{0,1,2,3\}$,&$\{0,1,4,5\}$,&$\{0,1,6,7\}$,&
$\{0,2,5,7\}$,&$\{0,3,4,7\}$,&$\{0,3,5,6\}$,\\
$\{1,2,4,7\}$,&$\{1,2,5,6\}$,&$\{1,3,4,6\}$,&
$\{2,3,4,5\}$,&$\{2,3,6,7\}$,&$\{4,5,6,7\}$.
\end{tabular}
\end{center}}
\noindent It is readily checked that each $3$-subset of $I_8$ is either
contained in exactly one block of $\cal T$ or in exactly one group of $\cal G$, but not in
both. Hence, $(X,{\cal G}\cup {\cal T})$ is a $3$-design. This is an example of $0$-fan designs.
\end{Example}

Next we give an example of $1$-fan designs. A {\em $1$-fan design} is a
$4$-tuple $(X,{\cal G},{\cal B},{\cal T})$ satisfying that $(X,{\cal G})$ is a $1$-design,
$(X,{\cal G}\cup {\cal B})$ is a $2$-design and $(X,{\cal G}\cup {\cal B} \cup {\cal T})$
is a $3$-design.

\begin{Example}
\label{1-fan(3^3)} Take $X=I_9$ and ${\cal G}=\{\{0,1,8\},\{2,3,6\},\{4,5,7\}\}$. Then $(X,{\cal G})$ is a $1$-design. Let $\cal B$ consists of the following $9$ blocks
{\footnotesize
\begin{center}
\begin{tabular}{llllll}
$\{2,4,8\}$,&$\{3,5,8\}$,&$\{0,2,7\}$,&
$\{0,3,4\}$,&$\{1,3,7\}$,&$\{1,4,6\}$,\\
$\{1,2,5\}$,&$\{0,5,6\}$,&$\{6,7,8\}$.
\end{tabular}
\end{center}}
\noindent It is readily checked that each $2$-subset of $I_9$ is either
contained in exactly one block of $\cal B$ or in exactly one group of $\cal G$, but not in
both. Hence, $(X,{\cal G}\cup {\cal B})$ is a $2$-design. Let $\cal T$ consists of the following $18$ blocks
{\footnotesize
\begin{center}
\begin{tabular}{llllll}
$\{0,1,2,3\}$,&$\{0,1,4,5\}$,&$\{0,1,6,7\}$,&
$\{0,2,4,6\}$,&$\{0,2,5,8\}$,&$\{0,3,5,7\}$,\\
$\{0,3,6,8\}$,&$\{0,4,7,8\}$,&$\{1,2,4,7\}$,&
$\{1,2,6,8\}$,&$\{1,3,4,8\}$,&$\{1,3,5,6\}$,\\
$\{1,5,7,8\}$,&$\{2,3,4,5\}$,&$\{2,3,7,8\}$,&
$\{2,5,6,7\}$,&$\{3,4,6,7\}$,&$\{4,5,6,8\}$.
\end{tabular}
\end{center}}
\noindent It is readily checked that $(X,{\cal G}\cup {\cal B}\cup {\cal T})$ is a $3$-design. This is an example of $1$-fan designs.
\end{Example}

If there are $a_i$ groups of size $g_i$ in an $s$-fan design, $1
\leq i \leq m$, then the {\em type} of the $s$-fan design is defined to be $g_1^{a_1} g_2^{a_2} \cdots g_m^{a_m}$. Let $K_1$, $K_2,\ldots,K_s$, $K_T$
be sets of positive integers. If block sizes of ${\cal
B}_i$ and ${\cal T}$ are from $K_i$ ($1\leq i\leq s$) and $K_T$,
respectively, then the $s$-fan design is denoted by
$s$-FG$(3,(K_1,K_2,\ldots,K_s,K_T),\sum_{i=1}^m g_i a_i)$ of type
$g_1^{a_1} g_2^{a_2}\cdots g_m^{a_m}$. Example \ref{0-fan(4^2)} shows a $0$-FG $(3,(\emptyset,4),8)$ of type $4^2$. Example \ref{1-fan(3^3)} presents a $1$-FG$(3,(3,4),9)$ of type $3^3$.

\begin{Lemma}
\label{nece-0-FG}{\rm (\cite{fcj})} The necessary conditions for the
existence of a $0$-FG$(3,(\emptyset,K_T),gn)$ of type $g^n$ $(n \geq
2)$ are

\begin{enumerate}
\item[{\rm (1)}] $g^2 n(n-1)(gn+g-3)\equiv 0$ $({\rm mod }$
$\alpha)$, where $\alpha=gcd\{k(k-1)(k-2):k\in K_T\}$;

\item[{\rm (2)}] $g(n-1)(gn+g-3)\equiv 0$ $({\rm mod }$ $\beta)$,
where $\beta=gcd\{(k-1)(k-2):k\in K_T\}$;

\item[{\rm (3)}] if $g=1$, then $n \equiv 2$ $({\rm mod }$ $\gamma);$ if $g> 1$,
then $gn \equiv g \equiv 2$ $({\rm mod }$ $\gamma)$, where
$\gamma=gcd\{k-2:k\in K_T\}$.
\end{enumerate}
\end{Lemma}

\begin{Theorem}
\label{G-design} {\rm (\cite{zkk})} There exists a
$0$-FG$(3,(\emptyset,4),gn)$ of type $g^n$ if and only if either
$g=1$ and $n\equiv 2,4$ $({\rm mod }$ $6)$, or $g$ is even and
$g(n-1)(n-2)\equiv 0$ $({\rm mod }$ $3)$.
\end{Theorem}

\subsubsection{The basic idea}

Since an optimal $2$-D $(u\times v,k,2)$-OOC is
equivalent to an optimal strictly $v$-cyclic $3$-$(u\times
v,k,1)$-packing, we first consider how to construct a $3$-packing without the restriction of automorphism groups.

\begin{itemize}
\item Step $1$: Start from a $0$-FG$(3,(\emptyset,k),gn)$ of type $g^n$ $(X,{\cal G},\emptyset,{\cal T})$. By the definition of $s$-fan designs, $(X,{\cal G}\cup {\cal T})$ is a $3$-design. $(X,{\cal T})$ satisfies that each $3$-subset of $X$ not contained in some group of $\cal G$ occurs in exactly one block of $\cal T$, and each $3$-subset of $X$ contained in some group of $\cal G$ never occur in any block of $\cal T$.
\item Step $2$: If a $3$-$(g,k,1)$-packing exists, then one can construct a $3$-$(g,k,1)$-packing on the set $G$ for each $G\in{\cal G}$. Denote its block set by ${\cal A}_G$.
\item Step $3$: Let ${\cal A}=\cup_{G\in{\cal G}}{\cal A}_G$. It follows that each $3$-subset of $X$ contained in some group of $\cal G$ occurs in at most one block of $\cal A$.
\item Step $4$: Let ${\cal C}={\cal A}\cup{\cal T}$. We have that $(X,{\cal C})$ is a $3$-$(gn,k,1)$-packing.
\end{itemize}

The main idea of the above construction is to fill in the groups of a $0$-fan design with a $3$-packing. So this construction is termed as ``Filling Construction''. Furthermore, if one hope to obtain an optimal $3$-$(gn,k,1)$-packing, it is necessary to input an optimal $3$-$(g,k,1)$-packing. Note that the reverse is not always correct. Now our purpose is to construct strictly $v$-cyclic $3$-$(u\times
v,k,1)$-packings. We need to modify the above ``Filling Construction" such that the initial $0$-fan design admits some special automorphisms.

An {\em automorphism} of an $s$-fan design $(X,{\cal G},{\cal
B}_1,{\cal B}_2,\ldots,{\cal B}_s,{\cal T})$ is a permutation on $X$
leaving ${\cal G}$, ${\cal B}_1$, ${\cal B}_2,\ldots,{\cal B}_s$,
${\cal T}$ invariant, respectively. All automorphisms of an $s$-fan
design form a group, called the {\em full automorphism group} of the
$s$-fan design. Any subgroup of the full automorphism group is
called an {\em automorphism group} of the $s$-fan design.

Let $G$ be an automorphism group of an $s$-fan design. All blocks of
the $s$-fan design can be partitioned into some orbits under $G$.
Choose any fixed block from each orbit and then call it a {\em base
block} of this $s$-fan design. For any block $B$ of the $s$-fan
design, the subgroup $\{\pi\in G: B^{\pi}=B\}$ is called the {\em
stabilizer} of $B$ in $G$.

\begin{Example}
\label{0-fan(4^2)-automorphism} Observe the $0$-FG $(3,(\emptyset,4),8)$ of type $4^2$ from Example \ref{0-fan(4^2)}. Consider the permutation $\alpha=(0\ 1\ 2\ 3)(4\ 5\ 6\ 7)$ on $I_8$. It is easy to checked that $\alpha$ is an automorphism of this $0$-FG. All blocks are partitioned into $5$ orbits under the action of $\alpha$. The $5$ base blocks are $\{0,1,2,3\}^*$, $\{4,5,6,7\}^*$, $\{1,2,5,6\}$, $\{0,1,6,7\}$ and $\{0,2,5,7\}$, where the stabilizer of each base block marked with a $^*$ is trivial, i.e., it contains only the identity permutation.
\end{Example}

In the following we introduce two kinds of $s$-fan designs with
special automorphism groups.

\subsubsection{$\bm h$-cyclic $\bm s$-fan designs}

Construct an $s$-fan design of type $(hg_1)^{a_1} (hg_2)^{a_2}
\cdots (hg_m)^{a_m}$ on $(\bigcup_{i=1}^m(I_{a_i}\times
I_{g_i}))\times Z_h$ with the group set $\{\{x\}\times I_{g_i}\times
Z_h: x\in I_{a_i},1\leq i\leq m\}$. If this $s$-fan design admits an
automorphism $\pi$ mapping $(x,y,j)\longmapsto(x,y,j+1)$ (mod
$(-,-,h)$), $x\in I_{a_i}$, $y\in I_{g_i}$ and $j\in Z_h$, then the
$s$-fan design is said to be {\em $h$-cyclic}.

For each block $B$ of an $h$-cyclic $s$-fan design of type
$(hg_1)^{a_1} (hg_2)^{a_2} \cdots (hg_m)^{a_m}$, if the stabilizer
of $B$ in $Z_h$ is trivial, i.e., $\{\delta\in Z_h:
B+\delta=B\}=\{0\}$, where $B+\delta=\{(x,y,j+\delta):(x,y,j)\in
B\}$, then the $s$-fan design is called {\em strictly $h$-cyclic}. A
(strictly) $h$-cyclic $s$-fan design of type $h^n$ is often referred
to as a {\em $($strictly$)$ semi-cyclic $s$-fan design} of type
$h^n$ (cf. \cite{fcj}).

The following construction is straightforward.

\begin{Construction}
\label{filling II} $(${\rm Filling Construction-I}$)$ Suppose that the following exist:

\begin{enumerate}
\item[$(1)$] a strictly $h$-cyclic
$0$-FG$(3,(\emptyset,k),\Sigma_{i=1}^m g_i a_ih)$ of type
$(hg_1)^{a_1} (hg_2)^{a_2} \cdots (hg_m)^{a_m}$ with $b_0$ base
blocks;
\item[$(2)$] a strictly $h$-cyclic $3$-$(g_i\times h,k,1)$ packing with $b_i$ base blocks for each $1\leq i\leq m$.
\end{enumerate}

\noindent Then there exists a strictly $h$-cyclic
 $3$-$((\Sigma_{i=1}^m g_i a_i)\times h, k,1)$ packing with $b_0+\sum_{i=1}^ma_ib_i$ base blocks, which is a $2$-D $((\Sigma_{i=1}^m g_i a_i)\times h,k,1)$-OOC.

Furthermore, if the given strictly $h$-cyclic $3$-$(g_i\times
h,k,1)$ packing is a strictly $h$-cyclic $S(3,k,g_i\times h)$ for
each $1\leq i\leq m$, then we obtain a strictly $h$-cyclic
$S(3,k,(\Sigma_{i=1}^m g_i a_i)\times h)$, which is a perfect $2$-D
$((\Sigma_{i=1}^m g_i a_i)\times h,k,2)$-OOC.
\end{Construction}

\begin{Example}
\label{example-2-cyclic-4^2} In this example, we construct an optimal
 $2$-D $(4\times 2,4,2)$-OOC.
\begin{itemize}
\item Step $1$: First we construct a strictly $2$-cyclic
$0$-FG$(3,(\emptyset,4),8)$ of type $4^2$ on $I_2\times I_2\times
Z_2$ with the group set $\{\{x\}\times I_2\times Z_2: x\in I_2\}$. All the $6$ base blocks are listed below.

{\footnotesize
\begin{center}
\begin{tabular}{ll}
$\{(0,0,0),(0,0,1),(1,0,0),(1,1,0)\}$, &
$\{(0,0,0),(1,0,0),(1,0,1),(0,1,0)\}$,\\
$\{(0,0,0),(1,0,0),(1,1,1),(0,1,1)\}$,&
$\{(0,0,0),(1,0,1),(1,1,0),(0,1,1)\}$,\\
$\{(0,0,0),(1,1,0),(1,1,1),(0,1,0)\}$,&
$\{(1,0,0),(1,1,0),(0,1,0),(0,1,1)\}$.
\end{tabular}
\end{center}}
\item Step $2$: Take an optimal strictly $2$-cyclic $3$-$(2\times 2,4,1)$
packing, which is trivial without base blocks.
\item Step $3$: Apply Construction $\ref{filling II}$ to obtain a strictly $2$-cyclic
 $3$-$(4\times 2, 4,1)$ packing with $6$ base blocks, which achieves the
  upper bound in Theorem $\ref{bound-total}$ and is an optimal
 $2$-D $(4\times 2,4,2)$-OOC with $6$ codewords. Note that $I_2\times I_2\times Z_2\cong I_4\times Z_2$. Hence $\Phi(4\times 2,4,2)=J^*(4\times
 2)=6$.
\end{itemize}
\end{Example}

\begin{Example}
\label{example-3-cyclic-6^2} In this example, we construct an optimal $2$-D $(4\times
  3,4,2)$-OOC.

\begin{itemize}
\item Step $1$: First we construct a strictly $3$-cyclic
$0$-FG$(3,(\emptyset,4),12)$ of type $6^2$ on $I_2\times I_2\times
Z_3$ with the group set $\{\{x\}\times I_2\times Z_3:x\in I_2\}$.
All the $15$ base blocks are listed below:

{\footnotesize
\begin{center}
\begin{tabular}{ll}
$\{(0,0,0),(0,0,1),(1,0,0),(1,0,1)\}$, &
$\{(0,0,0),(0,0,1),(1,0,2),(1,1,0)\}$,\\
$\{(0,0,0),(0,0,1),(1,1,1),(1,1,2)\}$,&
$\{(0,0,0),(0,1,0),(1,0,0),(1,1,0)\}$,\\
$\{(0,0,0),(0,1,0),(1,0,1),(1,1,1)\}$,&
$\{(0,0,0),(0,1,0),(1,0,2),(1,1,2)\}$,\\
$\{(0,0,0),(0,1,1),(1,0,0),(1,1,1)\}$,&
$\{(0,0,0),(0,1,1),(1,0,1),(1,0,2)\}$,\\
$\{(0,0,0),(0,1,1),(1,1,0),(1,1,2)\}$,&
$\{(0,0,0),(0,1,2),(1,0,0),(1,1,2)\}$,\\
$\{(0,0,0),(0,1,2),(1,0,1),(1,1,0)\}$,&
$\{(0,0,0),(0,1,2),(1,0,2),(1,1,1)\}$,\\
$\{(0,1,0),(0,1,1),(1,0,0),(1,0,2)\}$,&
$\{(0,1,0),(0,1,1),(1,0,1),(1,1,2)\}$,\\
$\{(0,1,0),(0,1,1),(1,1,0),(1,1,1)\}$.\\
\end{tabular}
\end{center}}
\item Step $2$: Construct an optimal strictly $3$-cyclic $3$-$(2\times 3,4,1)$
packing on $\{x\}\times I_2\times Z_3$ for each $x\in I_2$, which has $1$ base block and exists by Example $\ref{example-6}$. Then this step contributes $2$ base blocks as follows
{\footnotesize
\begin{center}
\begin{tabular}{ll}
$\{(0,0,0),(0,1,0),(0,0,1),(0,1,1)\}$, &
$\{(1,0,0),(1,1,0),(1,0,1),(1,1,1)\}$.
\end{tabular}
\end{center}}
\item Step $3$: Apply Construction $\ref{filling II}$ to obtain a strictly $3$-cyclic
 $(4\times 3, 4,1)$ packing with $17$ base blocks, which achieves the
  upper bound in Theorem $\ref{bound-total}$ and is an optimal $2$-D $(4\times
  3,4,2)$-OOC with $17$ codewords. Note that $I_2\times I_2\times Z_3\cong I_4\times Z_3$. Hence $\Phi(4\times 3,4,2)=J^*(4\times
 3)=17$.
\end{itemize}
\end{Example}

\begin{Lemma}
\label{(2v)^2} For any $v\equiv 1$ $({\rm mod}$ $2)$, there exists a
strictly $v$-cyclic $0$-FG$(3,(\emptyset,4),4v)$ of type $(2v)^2$.
\end{Lemma}
\proof By Lemma $2.11$ in \cite{fcj2}, there is a semi-cyclic
$0$-FG$(3,(\emptyset,4),4v)$ of type $(2v)^2$ on the point set
$X=I_2\times Z_{2v}$ and the group set ${\cal G}=\{\{x\}\times
Z_{2v}:x\in I_2\}$. Denote the family of its blocks by $\cal T$. For
each $(x,i)\in X$, define a mapping
$$\tau: \hspace{0.5cm} (x,i) \longmapsto (x,i-2\lfloor i/2\rfloor,\lfloor i/2\rfloor).$$
Let $X'=I_2\times I_2\times Z_v$ and ${\cal G}'=\{\{x\}\times
I_2\times Z_v:x\in I_2\}$. Let ${\cal T}'=\bigcup_{T\in{\cal T}}
\tau(T)$, where $\tau(T)=\{\tau(r):r\in T\}$. Since $v\equiv 1$
$({\rm mod}$ $2)$, it is readily checked that $(X',{\cal
G}',\emptyset,{\cal T}')$ is a strictly $v$-cyclic
$0$-FG$(3,(\emptyset,4),4v)$ of type $(2v)^2$. \qed

\begin{Remark}
\label{remark filling II} In Construction $\ref{filling II}$, even
if the given strictly $h$-cyclic $3$-$(g_i\times h,k,1)$ packing is
optimal for each $1\leq i\leq m$, the resulting strictly $h$-cyclic
 $3$-$((\sum_{i=1}^m g_i a_i)\times h, k,1)$ packing may not be
 optimal.
\end{Remark}

\subsubsection{$(u,h)$-regular $s$-fan designs}

Let $h$ divide $v$ and $H$ be a subgroup of order $h$ in $Z_v$,
i.e., $H=\{0,v/h,\ldots,(h-1)v/h\}$. Let $H_i=H+i$ be a coset of $H$
in $Z_v$, $0\leq i<v/h$. Construct an $s$-fan design of type
$(uh)^{v/h}$ on $I_u\times Z_v$ with the group set $\{I_u\times H_i:
0\leq i< v/h\}$. If this $s$-fan design admits an automorphism $\pi$
mapping $(x,j)\longmapsto(x,j+1)$ (mod $(-,v)$), $x\in I_u$ and
$j\in Z_v$, then the $s$-fan design is said to be {\em
$(u,h)$-regular}.

For each block $B$ of a $(u,h)$-regular $s$-fan design of type
$(uh)^{v/h}$, if the stabilizer of $B$ in $Z_v$ is trivial, i.e.,
$\{\delta\in Z_v: B+\delta=B\}=\{0\}$, where
$B+\delta=\{(x,j+\delta):(x,j)\in B\}$, then the $s$-fan design is
called {\em strictly $(u,h)$-regular}.

\begin{Example}
\label{0-fan(4^2)-excise} By Example \ref{0-fan(4^2)-automorphism}, the $0$-FG $(3,(\emptyset,4),8)$ of type $4^2$ from Example \ref{0-fan(4^2)} admits an automorphism $(0\ 1\ 2\ 3)(4\ 5\ 6\ 7)$. Actually the reader may check that this $0$-FG is isomorphic to a $(2,2)$-regular $0$-FG under the mapping $\tau:v\rightarrow(\lfloor v/4\rfloor,v\ (mod\ 4))$ from $I_8$ to $I_2\times Z_4$. But it is not strictly $(2,2)$-regular.
\end{Example}

When $u=1$, a (strictly)
$(1,h)$-regular $s$-fan design of type $h^{v/h}$ is often referred
to as a {\em $($strictly$)$ cyclic $s$-fan design} of type $h^{v/h}$
(cf. \cite{fcj}). We quote the following results for later use.

\begin{Lemma}
\label{strictly cyclic g^n}  {\rm (\cite{fcj})} \item{$(1)$} There
exists a strictly cyclic $0$-FG$(3,(\emptyset,4),2h)$ of type $h^2$
for any $h\equiv 0\ ({\rm mod }\ 8)$. \item{$(2)$} There exists a
strictly cyclic $0$-FG$(3,(\emptyset,4),3h)$ of type $h^3$ for any
$h\equiv 0\ ({\rm mod }\ 12)$.
\item{$(3)$} There exists a strictly cyclic
$0$-FG$(3,(\emptyset,4),5h)$ of type $h^5$ for any $h\equiv 0\ ({\rm
mod }\ 2)$.
\end{Lemma}

\begin{Lemma}
\label{equiv-RoSQS} {\rm (\cite{fcj1})} An $RoSQS(v+1)$ for $v\equiv
1$ $({\rm mod}$ $6)$ is equivalent to a strictly cyclic
$1$-$FG(3,(3,4),v)$ of type $1^v$. An $RoSQS(v+1)$ for $v\equiv 3$
$({\rm mod}$ $6)$ is equivalent to a strictly cyclic
$1$-$FG(3,(3,4),v)$ of type $3^{v/3}$.
\end{Lemma}

\begin{Construction}
\label{filling III} $(${\rm Filling Construction-II}$)$ Let $uh\geq k\geq 3$. Suppose that the following exist.

\begin{enumerate}
\item[$(1)$] a strictly $(u,h)$-regular $0$-FG$(3,(\emptyset,k),uv)$ of type
$(uh)^{v/h}$ with $b_0$ base blocks;
\item[$(2)$] a
 strictly $h$-cyclic $3$-$(u\times h,k,1)$ packing with $b_1$ base blocks.
\end{enumerate}

\noindent Then there exists a strictly $v$-cyclic
 $3$-$(u\times v, k,1)$ packing with $b_0+b_1$ base blocks.

Furthermore, if the given strictly $h$-cyclic $3$-$(u\times h,k,1)$
 packing is optimal with $J(u\times h,k,2)$ base blocks, then the
 derived strictly $v$-cyclic $3$-$(u\times v, k,1)$ packing is also
 optimal with $J(u\times v,k,2)$ base blocks, which is an optimal $2$-D $(u\times v,k,2)$-OOC with $J(u\times v,k,2)$ codewords.
\end{Construction}

\proof First we prove the first part of this construction.
\begin{itemize}
\item Step $1$: Start from a strictly
$(u,h)$-regular $0$-FG$(3,(\emptyset,k),uv)$ of type $(uh)^{v/h}$.
Denote the family of base blocks of this design by ${\cal F}$.
\item Step $2$: Let ${\cal E}$ be the family of base blocks of a strictly $h$-cyclic
$3$-$(u\times h,k,1)$ packing. For each
$B=\{(x_1,j_1),(x_2,j_2),\ldots,(x_k,j_k)\}\in {\cal E}$ we take
$${\frac{v}{h}}B=\{(x_1,{\frac{v}{h}}j_1),(x_2,{\frac{v}{h}}j_2),\ldots,(x_k,{\frac{v}{h}}j_k)\}.$$
\item Step $3$: Then ${\cal F}\cup \{\frac{v}{h}B:B\in {\cal E}\}$ forms the family
of base blocks of the desired strictly $v$-cyclic $(u\times v, k,1)$
packing.
\end{itemize}

For checking optimality of the required design in the second part,
it suffices to show that
\begin{eqnarray}
\label{p1}  &&\frac{u((uv-1)(uv-2)-(uh-1)(uh-2))}{k(k-1)(k-2)} +
\lfloor\frac{u}{k}\lfloor\frac{uh-1}{k-1}
\lfloor\frac{uh-2}{k-2}\rfloor\rfloor\rfloor \nonumber\\
&&=\lfloor\frac{u}{k}\lfloor\frac{uv-1}{k-1}
\lfloor\frac{uv-2}{k-2}\rfloor\rfloor\rfloor.
\end{eqnarray}

\noindent By Lemma \ref{nece-0-FG} $(3)$, since $uh>1$, the existence of a
strictly $(u,h)$-regular $0$-FG$(3,(\emptyset,k)$, $uv)$ of type
$(uh)^{v/h}$ implies that $uv-2 \equiv uh-2 \equiv 0$ $({\rm mod }$
$k-2)$. By Lemma \ref{nece-0-FG} $(2)$, one can verify that
$(uv-1)(uv-2) \equiv (uh-1)(uh-2)$ $({\rm mod }$ $(k-1)(k-2))$. Let
$(uv-1)(uv-2)=a_1(k-1)(k-2)+r$ and $(uh-1)(uh-2)=a_2(k-1)(k-2)+r$,
where $0\leq r<(k-1)(k-2)$. Thus for obtaining the equation
(\ref{p1}), it suffices to prove that
\begin{equation}
\label{p2}\frac{u(a_1-a_2)}{k}+
\lfloor\frac{ua_2}{k}\rfloor=\lfloor\frac{ua_1}{k}\rfloor.
\end{equation}
Note that $u(a_1-a_2) \equiv 0$ $({\rm mod }$ $k)$. Let
$ua_1=b_1k+r_1$ and $ua_2=b_2k+r_1$, where $0\leq r_1<k$. It is
readily   checked that the equation (\ref{p2}) holds. \qed

\begin{Example}
\label{example-4-cyclic-4^2} In this example, we construct an optimal
 $2$-D $(2\times 4,4,2)$-OOC.

\begin{itemize}
\item Step $1$: First we construct a strictly
$(2,2)$-regular $0$-FG$(3,(\emptyset,4),8)$ of type $4^2$ on
$I_2\times Z_4$ with the group set $\{I_2\times H_i:0\leq i\leq
1\}$, where $H_0=\{0,2\}$ is a subgroup of
order $2$ in $Z_4$ and $H_1=\{1,3\}$. All the $3$ base blocks are listed below:

{\footnotesize
\begin{center}
\begin{tabular}{lll}
$\{(0,0),(0,1),(0,2),(1,1)\}$, & $\{(0,0),(0,1),(1,2),(1,3)\}$,&
$\{(0,0),(1,0),(1,1),(1,3)\}$.
\end{tabular}
\end{center}}

\item Step $2$: Take an optimal strictly $2$-cyclic $3$-$(2\times 2,4,1)$
packing, which is trivial without base blocks.

\item Step $3$: Apply Construction $\ref{filling III}$ to obtain a strictly $4$-cyclic
 $3$-$(2\times 4, 4,1)$ packing with $3$ base blocks, which achieves the
  upper bound in Theorem $\ref{bound-total}$ and is an optimal
 $2$-D $(2\times 4,4,2)$-OOC with $3$ codewords. Hence $\Phi(2\times 4,4,2)=J^*(2\times
 4)=3$.
\end{itemize}
\end{Example}

\begin{Example}
\label{example 8^2}  In this example, we construct an optimal
 $2$-D $(2\times 8,4,2)$-OOC.
\begin{itemize}
\item Step $1$: First we construct a strictly $(2,4)$-regular
$0$-FG$(3,(\emptyset,4),16)$ of type $8^2$ on $I_2\times Z_8$ with
the group set $\{I_2\times H_i:0\leq i\leq 1\}$, where $H_0=\{0,2,4,6\}$ is a subgroup of order $4$ in $Z_8$ and $H_1=\{1,3,5,7\}$. All the $14$ base
blocks are listed below:

{\footnotesize
\begin{center}
\begin{tabular}{lll}
$\{(0,0),(0,1),(0,2),(0,5)\}$, & $\{(0,0),(0,1),(0,3),(1,0)\}$,&
$\{(0,0),(0,1),(0,6),(1,1)\}$,\\
$\{(0,0),(0,1),(1,2),(1,3)\}$,& $\{(0,0),(0,1),(1,4),(1,5)\}$,&
$\{(0,0),(0,1),(1,6),(1,7)\}$,\\
$\{(0,0),(0,2),(1,1),(1,5)\}$,& $\{(0,0),(0,3),(1,1),(1,6)\}$,&
$\{(0,0),(0,3),(1,2),(1,5)\}$,\\
$\{(0,0),(0,3),(1,4),(1,7)\}$,& $\{(0,0),(0,4),(1,1),(1,7)\}$,&
$\{(0,0),(1,0),(1,1),(1,3)\}$,\\
$\{(0,0),(1,0),(1,5),(1,7)\}$,& $\{(1,0),(1,1),(1,2),(1,5)\}$.
\end{tabular}
\end{center}}

\item Step $2$: Construct an optimal strictly $4$-cyclic $3$-$(2\times 4,4,1)$
packing with $3$ base blocks, which exists by Example $\ref{example-4-cyclic-4^2}$. Then this step contributes $3$ base blocks as follows

{\footnotesize
\begin{center}
\begin{tabular}{lll}
$\{(0,0),(0,2),(0,4),(1,2)\}$, & $\{(0,0),(0,2),(1,4),(1,6)\}$,&
$\{(0,0),(1,0),(1,2),(1,6)\}$.
\end{tabular}
\end{center}}

\item Step $3$: Apply Construction $\ref{filling III}$ to obtain a strictly $8$-cyclic
 $(2\times 8,4,1)$ packing with $17$ base blocks, which achieves the
  upper bound in Theorem $\ref{bound-total}$ and is an optimal
 $2$-D $(2\times 8,4,2)$-OOC with $17$ codewords. Hence $\Phi(2\times 8,4,2)=J^*(2\times
 8)=17$.
\end{itemize}
\end{Example}

The following result is simple but very useful.

\begin{Lemma}
\label{relation between regular and v-cyclic} If there exists a
strictly $(u,h)$-regular $0$-FG$(3,(\emptyset,k),uv)$ of type
$(uh)^{v/h}$, then for any integer divisor $h_1$ of $h$, there
exists a strictly $h_1$-cyclic $0$-FG$(3,(\emptyset,k),uv)$ of type
$(uh)^{v/h}$.
\end{Lemma}

\begin{Corollary}
\label{(2v)^2-1} For any $v\not\equiv 2$ $({\rm mod}$ $4)$, there
exists a strictly $v$-cyclic $0$-FG$(3,(\emptyset,4),4v)$ of type
$(2v)^2$.
\end{Corollary}

\proof When $v\equiv 0$ $({\rm mod}$ $4)$, by Lemma \ref{strictly
cyclic g^n} there is a strictly cyclic $0$-FG$(3,(\emptyset,4),4v)$
of type $(2v)^2$. Apply Lemma \ref{relation between regular and
v-cyclic} to obtain a strictly $v$-cyclic
$0$-FG$(3,(\emptyset,4),4v)$ of type $(2v)^2$. When $v\equiv 1$
$({\rm mod}$ $2)$, the conclusion follows from Lemma \ref{(2v)^2}.
\qed

\begin{Lemma}
\label{strictly $(2,1)$-regular $1$-FG} If there is a perfect $2$-D
$(2\times v,4,2)$-OOC with $v\equiv 1,5$ $({\rm mod}$ $6)$, then
there is a strictly $(2,1)$-regular $1$-FG$(3,(2,4),2v)$ of type
$2^v$.
\end{Lemma}

\proof By Lemma \ref{nece-perfect}, the necessary condition for the
existence of a perfect $2$-D $(2\times v,4,2)$-OOC is $v\equiv 1,5$
$({\rm mod}$ $6)$. Suppose that $(X,{\cal T})$ is a strictly
$v$-cyclic SQS$(2\times v)$ with $X=I_2\times Z_v$, which is
equivalent to a perfect $2$-D $(2\times v,4,2)$-OOC. Let ${\cal
G}=\{I_2\times \{x\}:x\in Z_v\}$. Then $(X,{\cal G},\emptyset,{\cal
T})$ is a strictly $(2,1)$-regular $0$-FG$(3,(\emptyset,4),2v)$ of
type $2^v$. Collect all $2$-subsets of $X$ from distinct groups of
$\cal G$ into a set $\cal B$. Since $v$ is odd, $(X,{\cal G},{\cal
B},{\cal T})$ is a strictly $(2,1)$-regular $1$-FG$(3,(2,4),2v)$ of
type $2^v$. \qed

\subsection{$H$ designs}

Mills first used the terminology of $H$ designs in \cite{mills74}. Let $n$, $g$, $t$ be positive integers and $K$ be a set of positive
integers. An {\em $H$ design} is a triple $(X,{\cal G},{\cal B})$,
where ${\cal G}$ is a partition of a set of points $X$ into $n$
subsets (called {\em groups}), each of cardinality $g$, and ${\cal
B}$ is a collection of subsets of $X$ (called {\em blocks}), each of
cardinality from $K$, such that each block intersects any given
group in at most one point, and each $t$-subset of $X$ from $t$
distinct groups is contained in a unique block. Such a design is
denoted by $H(n,g,K,t)$.

\begin{Example}
\label{H(4,2,4,3)} Take $X=I_8$ and ${\cal G}=\{\{0,4\},\{1,5\},\{2,6\},\{3,7\}\}$. Let $\cal B$ consists of the following $8$ blocks
{\footnotesize
\begin{center}
\begin{tabular}{llllll}
$\{0,1,2,3\}$,&$\{4,5,6,7\}$,&
$\{0,1,6,7\}$,&$\{2,3,4,5\}$,&
$\{0,2,5,7\}$,&$\{1,3,4,6\}$,\\
$\{0,3,5,6\}$,&$\{1,2,4,7\}$.
\end{tabular}
\end{center}}
\noindent It is easy to see that each $3$-subset of $I_8$ from three
distinct groups of $\cal G$ is contained in a unique block of $\cal B$. Then $(X,{\cal G},{\cal B})$ is an $H(4,2,4,3)$.
\end{Example}

\begin{Example}
\label{subdesigns in an s-fan} Let $(X,{\cal
G},{\cal B}_1,{\cal B}_2,\ldots,{\cal B}_s,{\cal T})$ be an
$s$-$FG(3,(K_1,K_2,\ldots,K_s,K_T),gn)$ of type $g^n$. Then for each
$1\leq i\leq s$, $(X,{\cal G},{\cal B}_i)$ is an $H(n,g,K_i,2)$
$($called the $i$-th subdesign of the $s$-fan design$)$.
\end{Example}

\begin{Lemma}
\label{H design} $(${\rm \cite{ji1,mills}}$)$ For any $n\geq 4$,
$n\neq 5$, an $H(n,g,4,3)$ exists if and only if $gn$ is even and
$g(n-1)(n-2)$ is divisible by $3$. For $n=5$, an $H(5,g,4,3)$ exists
if $g$ is even, $g\neq 2$ and $g\not\equiv 10,26$ $({\rm mod}$
$48)$.
\end{Lemma}

An {\em automorphism} of an $H$ design $(X,{\cal G},{\cal B})$ is a
permutation on $X$ leaving ${\cal G}$, ${\cal B}$ invariant,
respectively. All automorphisms of an $H$ design form a group,
called the {\em full automorphism group} of the $H$ design. Any
subgroup of the full automorphism group is called an {\em
automorphism group} of the $H$ design.

Let $G$ be an automorphism group of an $H$ design. All blocks of the
$H$ design can be partitioned into some orbits under $G$. Choose any
fixed block from each orbit and then call it a {\em base block} of
this $H$ design. For any block $B$ of the $H$ design, the subgroup
$\{\pi\in G: B^{\pi}=B\}$ is called the {\em stabilizer} of $B$ in
$G$.

\begin{Example}
\label{H(4,2,4,3)-automorphism} Observe the $H(4,2,4,3)$ from Example \ref{H(4,2,4,3)}. Consider the permutation $(0\ 4)(1\ 5)(2\ 6)(3\ 7)$ on $I_8$. It is readily checked that $\alpha$ is an automorphism of this $H$ design. All blocks are partitioned into $4$ orbits under the action of $\alpha$. The $4$ base blocks are $\{0,1,2,3\}$, $\{0,1,6,7\}$, $\{0,2,5,7\}$, $\{0,3,5,6\}$.
\end{Example}

Construct an $H(n,lh,K,t)$ on $I_n\times I_l\times Z_h$ with the
group set $\{\{x\}\times I_l\times Z_h: x\in I_n\}$. If this $H$
design admits an automorphism $\pi$ mapping
$(x,y,j)\longmapsto(x,y,j+1)$ (mod $(-,-,h)$), $x\in I_n$, $y\in
I_l$ and $j\in Z_h$, then the $H$ design is said to be {\em
$h$-cyclic}. If the stabilizer of each block of an $h$-cyclic
$H(n,lh,K,t)$ in $Z_h$ is trivial, i.e., for any block $B$,
$\{\delta\in Z_h: B+\delta=B\}=\{0\}$, where
$B+\delta=\{(x,y,j+\delta):(x,y,j)\in B\}$, then the $H$ design is
called {\em strictly $h$-cyclic}. Note that one can verify
that an $h$-cyclic $H$ design is always strictly $h$-cyclic.

\begin{Example}
\label{semi-cyclic H(4,2,4,3)-excise} By Example \ref{H(4,2,4,3)-automorphism}, the $H(4,2,4,3)$ from Example \ref{H(4,2,4,3)} admits an automorphism $(0\ 4)(1\ 5)(2\ 6)(3\ 7)$. Actually the reader may check that this $H$ design is isomorphic to a $2$-cyclic $H(4,2,4,3)$ under the mapping $\tau:v\rightarrow(v\ (mod\ 4),0,\lfloor v/4\rfloor)$ from $I_8$ to $I_4\times I_1\times Z_2$.
\end{Example}

When $l=1$, an $h$-cyclic $H(n,h,K,t)$ is often referred to as a {\em
semi-cyclic $H(n,h,K,t)$}.

\begin{Lemma}
\label{H(4,g,4,3)} $(${\rm \cite{fcj}}$)$ For any $h\geq 1$, there
exists a semi-cyclic $H(4,h,4,3)$.
\end{Lemma}

\section{Weighting constructions}

For applying Constructions \ref{filling II} and \ref{filling III}, we need some strictly
$h$-cyclic $0$-FGs and strictly $(u,h)$-regular $s$-FGs. Construction \ref{Weighting-strictly $h$-cyclic $s$-FG} shows that if one has a strictly $h_1$-cyclic $1$-FG of type $(g_1h_1)^n$, and gives each point of the $1$-FG a weight $g_2h_2$, then a strictly $h_1h_2$-cyclic $0$-FG of type $(g_1g_2h_1h_2)^{n}$ can be obtained; Construction \ref{Weighting-strictly regular $s$-FG} shows that if one has a strictly $(g_1,h_1)$-regular $1$-FG of type $(g_1h_1)^n$, and gives each point of the $1$-FG a weight $g_2h_2$, then a strictly $(g_1g_2,h_1h_2)$-regular $0$-FG of type
$(g_1g_2h_1h_2)^n$ can be obtained. So Constructions \ref{Weighting-strictly $h$-cyclic $s$-FG} and \ref{Weighting-strictly regular $s$-FG} give an approach to find some infinite families of strictly $h$-cyclic $0$-FGs and strictly $(u,h)$-regular $s$-FGs. Then apply Constructions \ref{filling II} and \ref{filling III} to fill in the groups of these infinite families. We can obtain many optimal $2$-D OOCs, which will be presented in Sections $8$ and $9$.

Condition $(3)$ in Constructions \ref{Weighting-strictly $h$-cyclic $s$-FG} and \ref{Weighting-strictly regular $s$-FG} implies that $h$-cyclic $H$ designs are important. Thus a recursive construction for $h$-cyclic $H$ designs is presented in Construction \ref{Weighting for $h$-cyclic $H$ design}. The proofs of all constructions in this section are of design theory. Here we only focus on how these constructions work. The detailed proofs of Constructions \ref{Weighting-strictly $h$-cyclic $s$-FG} and \ref{Weighting for $h$-cyclic $H$ design} have been moved to Appendix III. The detailed proof of Construction \ref{Weighting-strictly regular $s$-FG} is omitted, which is similar to that of Construction \ref{Weighting-strictly $h$-cyclic $s$-FG}.

\begin{Construction}
\label{Weighting-strictly $h$-cyclic $s$-FG} $(${\rm Weighting Construction-I}$)$ Let $K$ and $L_i$ for each $1\leq i\leq s$ be all sets of positive integers greater than
$1$. Let $K_T$ and $L_T$ be both sets of positive integers greater
than $2$. Suppose that the following exist:
\begin{enumerate}
\item[$(1)$] a strictly $h_1$-cyclic $1$-FG$(3,(K,K_T),$ $ng_1h_1)$ of type $(g_1h_1)^n$ $($called the master design$)$;
\item[$(2)$] a strictly $h_2$-cyclic $s$-FG$(3,(L_1,L_2,\ldots,L_s,L_T),kg_2h_2)$ of type
$(g_2h_2)^{k}$ for each $k\in K$;
\item[$(3)$] an $h_2$-cyclic $H(k,g_2h_2,L_T,3)$ for each $k\in K_T$.
\end{enumerate}
\noindent Then there exists a strictly
$h_1h_2$-cyclic $s$-FG$(3,(L_1,L_2,$ $\ldots,L_s,L_T)$,
$ng_1g_2h_1h_2)$ of type $(g_1g_2h_1h_2)^{n}$.
\end{Construction}

 {\footnotesize
\begin{algorithm}[H]
\caption{Algorithm in Construction \ref{Weighting-strictly $h$-cyclic $s$-FG}}

{\rm \textbf{Step $1$}}: Start from
\begin{center}
a strictly
$h_1$-cyclic $1$-FG$(3,(K,K_T),$ $ng_1h_1)$ of type $(g_1h_1)^n$ $(X,{\cal G},{\cal B},{\cal T})$,
\end{center}
where $X=I_n\times I_{g_1}\times Z_{h_1}$ and ${\cal G}=\{\{x\}\times
I_{g_1}\times Z_{h_1}:x\in I_{n}\}$.
\begin{itemize}
\item Denote the family of base
blocks of this design by ${\cal F}={\cal F}_1\cup {\cal F}_2$, where
${\cal F}_1$ and ${\cal F}_2$ \\generate all the blocks of ${\cal B}$
and ${\cal T}$ respectively.
\end{itemize}

{\rm \textbf{Step $2$ (input)}}: For any base block $B\in {\cal F}_1$, construct
\begin{center}
a strictly
$h_2$-cyclic $s$-FG$(3,(L_1,L_2$ $\ldots,L_s,L_T),|B|g_2h_2)$ of
type $(g_2h_2)^{|B|}$
\end{center}
on $B\times I_{g_2}\times Z_{h_2}$ with the
group set $\{\{x\}\times I_{g_2}\times Z_{h_2}: x\in B\}$.
\begin{itemize}
\item Denote
the family of base blocks of the $j$-th subdesign
$H(|B|,g_2h_2,L_j,2)$ by ${\cal A}_B^j$ for \\$1\leq j\leq s$. Denote the family of all the other base blocks by ${\cal D}_B$.
\end{itemize}

{\rm \textbf{Step $3$ (input)}}: For any base block $B\in {\cal F}_2$, construct
\begin{center}
an $h_2$-cyclic $H(|B|,g_2h_2,L_T,3)$
\end{center} on $B\times I_{g_2}\times Z_{h_2}$ with the
group set $\{\{x\}\times I_{g_2}\times Z_{h_2}: x\in B\}$.
\begin{itemize}
\item Denote
the family of base blocks of this design by ${{\cal D}'_B}$.
\end{itemize}

{\rm \textbf{Step $4$ (mapping)}}: Let $${\cal A}_j=\bigcup_{B\in {\cal F}_1} {\cal A}_B^j\ for\ 1\leq j
\leq s,\hspace{0.5cm} {\cal D} = (\bigcup_{B\in {\cal F}_1} {\cal
D}_B)\bigcup (\bigcup_{B\in {\cal F}_2} {{\cal D}'_B}).$$ For each
$C\in(\bigcup_{1\leq j\leq s}{\cal A}_j)\bigcup {\cal D}$ and each
$(x,y,z,u,v)\in C$, define a mapping
$$\tau: \hspace{0.5cm} (x,y,z,u,v) \longmapsto (x,y+ug_1,z+vh_1).$$ Define $\tau(C)=\{\tau(c):c\in C\}$.
Let
$$\hspace{1cm}  {{\cal A}_j^*}=\bigcup_{_{C\in {\cal A}_j}} \tau(C),
\ 1\leq j \leq s,\hspace{0.5cm}  {\cal D}^*=\bigcup_{_{C\in {\cal
D}}} \tau(C).$$

{\rm \textbf{Step $5$ (final)}}: Take $${{\cal A}'_j}=\{A+\delta:A\in {\cal A}_j^*,\delta\in
Z_{h_1h_2}\},\hspace{0.5cm} {{\cal D}}'=\{A+\delta:A\in {\cal
D}^*,\delta\in Z_{h_1h_2}\},$$ where $A+\delta=\{(x,y,z+\delta)\
({\rm mod}\ (-,-,h_1h_2)):(x,y,z)\in A\}$. Take $$X'=I_n\times
I_{g_1g_2}\times Z_{h_1h_2},\hspace{0.5cm} {\cal G}'=\{\{x\}\times
I_{g_1g_2}\times Z_{h_1h_2}: x\in I_n\}.$$ Then $(X',{\cal G}',{{\cal A}_1}',\ldots,{{\cal A}_s}',{\cal D}')$
is the required strictly $h_1h_2$-cyclic $s$-FG$(3,(L_1,L_2,$
$\ldots,L_s,L_T),$ $ng_1g_2h_1h_2)$ of type $(g_1g_2h_1h_2)^{n}$.
\end{algorithm}}

\vspace{0.5cm}

The following example illustrates the algorithm presented in Figure $2$.

\begin{Example}
\label{example-2-cyclic-8^2} In this example, we construct an optimal
 $2$-D $(8\times 2,4,2)$-OOC.
\begin{itemize}
\item Step $1$: First construct a strictly
$1$-cyclic $1$-FG$(3,(2,4),$ $4)$ of type $1^4$ $(X,{\cal G},{\cal B},{\cal T})$ on $X=I_4\times I_1\times Z_1$ with the group set ${\cal G}=\{\{x\}\times
I_1\times Z_1:x\in I_4\}$, which is trivial. Take
{\footnotesize$${\cal F}_1=\{\{(i,0,0),(j,0,0)\}:\{i,j\}\in\{\{0,1\},\{0,2\},\{0,3\},\{1,2\},\{1,3\},\{2,3\}\}\},$$} which generates $6$ blocks of $\cal B$ under $Z_1$ such that $(X,{\cal G}\cup{\cal B})$ is a $2$-design. Take {\footnotesize$${\cal F}_2=\{\{(0,0,0),(1,0,0),(2,0,0),(3,0,0)\}\},$$} which generates the unique block of $\cal T$ such that $(X,{\cal G}\cup{\cal B}\cup{\cal T})$ is a $3$-design.

\item Step $2$: For each $B=\{(i,0,0),(j,0,0)\}\in {\cal F}_1$, construct a strictly
$2$-cyclic $0$-FG$(3,(\emptyset,4),8)$ of
type $4^2$ on $B\times I_2\times Z_2$ with the
group set $\{\{x\}\times I_2\times Z_2: x\in B\}$, which exists by Example \ref{example-2-cyclic-4^2}.
All the $6$ base blocks of ${\cal D}_B$ are listed below.

{\footnotesize
\begin{center}
\begin{tabular}{l}
$\{(i,0,0,0,0),(i,0,0,0,1),(j,0,0,0,0),(j,0,0,1,0)\}$,\\
$\{(i,0,0,0,0),(j,0,0,0,0),(j,0,0,0,1),(i,0,0,1,0)\}$,\\
$\{(i,0,0,0,0),(j,0,0,0,0),(j,0,0,1,1),(i,0,0,1,1)\}$,\\
$\{(i,0,0,0,0),(j,0,0,0,1),(j,0,0,1,0),(i,0,0,1,1)\}$,\\
$\{(i,0,0,0,0),(j,0,0,1,0),(j,0,0,1,1),(i,0,0,1,0)\}$,\\
$\{(j,0,0,0,0),(j,0,0,1,0),(i,0,0,1,0),(i,0,0,1,1)\}$.
\end{tabular}
\end{center}}

\item Step $3$: For the unique $B\in {\cal F}_2$, construct a $2$-cyclic $H(4,4,4,3)$
on $B\times I_2\times Z_2$ with the group set $\{\{x\}\times I_2\times Z_2: x\in B\}$, which exists by Corollary \ref{cor for $h$-cyclic $H$ design}.
Denote the family of base blocks of this design by ${\cal D}'_B$, and $|{\cal D}'_B|=32$.

\item Step $4$: Let ${\cal D} = (\bigcup_{B\in {\cal F}_1} {\cal
D}_B)\bigcup (\bigcup_{B\in {\cal F}_2} {{\cal D}'_B}).$ For each
$C\in{\cal D}$ and each
$(x,y,z,u,v)$ $\in C$, define a mapping
$\tau: (x,y,z,u,v) \longmapsto (x,y+u,z+v).$ Define $\tau(C)=\{\tau(c):c\in C\}$.
Let ${\cal D}^*=\bigcup_{_{C\in {\cal D}}} \tau(C).$ Then $|{\cal D}^*|=68$, which is just the number of base blocks in a strictly $2$-cyclic $0$-FG$(3,(\emptyset,4),16)$ of type $4^4$.

\item Step $5$: Let ${{\cal D}}'=\{A+\delta:A\in {\cal D}^*,\delta\in Z_2\},$ where $A+\delta=\{(x,y,z+\delta)\ ({\rm mod}\ (-,-,2)):(x,y,z)\in A\}$. Take $X'=I_4\times I_2\times Z_2$ and ${\cal G}'=\{\{x\}\times I_2\times Z_2: x\in I_4\}.$ Then $(X',{\cal G}',\emptyset,{\cal D}')$ is the required strictly $2$-cyclic $0$-FG$(3,(\emptyset,4),16)$ of type $4^4$.

\item Step $6$: Apply Construction $\ref{filling II}$. Fill in the groups of the resulting strictly $2$-cyclic $0$-FG$(3,(\emptyset,4),16)$ of type $4^4$ with a trivial optimal strictly $2$-cyclic $3$-$(2\times 2, 4,1)$ packing without base blocks. We have an optimal strictly $2$-cyclic $3$-$(8\times 2, 4,1)$ packing with $68$ base blocks, which achieves the upper bound in Theorem $\ref{bound-total}$ and is an optimal $2$-D $(8\times 2,4,2)$-OOC. Hence $\Phi(8\times 2,4,2)=J^*(8\times 2)=68$.
\end{itemize}
\end{Example}


\begin{Construction}
\label{Weighting-strictly regular $s$-FG} $(${\rm Weighting Construction-II}$)$ Let $K$ and $L_i$ for each $1\leq i\leq s$ be all sets of positive integers greater than $1$.
Let $K_T$ and $L_T$ be both sets of positive integers greater than
$2$. Suppose that the following exist:

\begin{enumerate}
\item[$(1)$] a strictly $(g_1,h_1)$-regular
$1$-FG$(3,(K,K_T)$, $g_1h_1n)$ of type $(g_1h_1)^n$;

\item[$(2)$] a strictly $h_2$-cyclic $s$-FG$(3,(L_1,L_2,\ldots,L_s,L_T)$,
$kg_2h_2)$ of type $(g_2h_2)^k$ for each $k\in K$;

\item[$(3)$] an $h_2$-cyclic $H(k,g_2h_2,L_T,3)$ for each $k\in K_T$.
\end{enumerate}

\noindent Then there exists a strictly $(g_1g_2,h_1h_2)$-regular
$s$-FG$(3,(L_1,L_2,\ldots,L_s,L_T),g_1g_2h_1h_2n)$ of type
$(g_1g_2h_1h_2)^n$.
\end{Construction}

 {\footnotesize
\begin{algorithm}[H]
\caption{Algorithm in Construction \ref{Weighting-strictly regular $s$-FG}}

{\rm \textbf{Step $1$}}: Start from
\begin{center}
a strictly $(g_1,h_1)$-regular $1$-FG$(3,(K,K_T),g_1h_1n)$ of type
$(g_1h_1)^n$ $(X,{\cal G},{\cal B},{\cal T})$,
\end{center}
on $X=I_{g_1}\times Z_{h_1n}$ with the group set ${\cal
G}=\{I_{g_1}\times H_i:0\leq i< n\}$, where $H=\{0,n,\ldots,(h_1-1)n\}$ is a subgroup of order $h_1$ in $Z_{h_1n}$, and $H_i=H+i$ be a coset of $H$ in
$Z_{h_1n}$, $0\leq i<n$.
\begin{itemize}
\item Denote the family of base
blocks of this design by ${\cal F}={\cal F}_1\cup {\cal F}_2$, where
${\cal F}_1$ and ${\cal F}_2$ generate \\all the blocks of ${\cal B}$
and ${\cal T}$ respectively.
\end{itemize}

{\rm \textbf{Step $2$ (input)}}: For any base block $B\in {\cal F}_1$, construct
\begin{center}
a strictly
$h_2$-cyclic $s$-FG$(3,(L_1,L_2\ldots,L_s,L_T)$, $|B|g_2h_2)$ of
type $(g_2h_2)^{|B|}$
\end{center}
on $B\times I_{g_2}\times Z_{h_2}$ with the
group set $\{\{x\}\times I_{g_2}\times Z_{h_2}: x\in B\}$.
\begin{itemize}
\item Denote
the family of base blocks of the $j$-th subdesign
$H(|B|,g_2h_2,L_j,2)$ by ${\cal A}_B^j$ for $1\leq j$ \\$\leq s$, and
denote the family of all the other base blocks by ${\cal D}_B$.
\end{itemize}

{\rm \textbf{Step $3$ (input)}}: For any base block $B\in {\cal F}_2$, construct
\begin{center}
an $h_2$-cyclic $H(|B|,g_2h_2,L_T,3)$
\end{center} on $B\times I_{g_2}\times Z_{h_2}$ with the
group set $\{\{x\}\times I_{g_2}\times Z_{h_2}: x\in B\}$.
\begin{itemize}
\item Denote
the family of base blocks of this design by ${{\cal D}'_B}$.
\end{itemize}

{\rm \textbf{Step $4$ (mapping)}}: Let $${\cal A}_j=\bigcup_{B\in {\cal F}_1} {\cal A}_B^j\ for\ 1\leq j
\leq s,\hspace{0.5cm} {\cal D} = (\bigcup_{B\in {\cal F}_1} {\cal
D}_B)\bigcup (\bigcup_{B\in {\cal F}_2} {{\cal D}'_B}).$$ For each
$C\in(\bigcup_{1\leq j\leq s}{\cal A}_j)\bigcup {\cal D}$ and each
$(x,y,z,u)\in C$, define a mapping
$$\tau: \hspace{0.5cm} (x,y,z,u) \longmapsto (x+zg_1,y+uh_1n).$$ Define $\tau(C)=\{\tau(c):c\in C\}$.
Let
$$\hspace{1cm}  {{\cal A}_j^*}=\bigcup_{_{C\in {\cal A}_j}} \tau(C),
\ 1\leq j \leq s,\hspace{0.5cm}  {\cal D}^*=\bigcup_{_{C\in {\cal
D}}} \tau(C).$$

{\rm \textbf{Step $5$ (final)}}: Take $${{\cal A}'_j}=\{A+\delta:A\in {\cal A}_j^*,\delta\in
Z_{h_1h_2n}\},\hspace{0.5cm} {{\cal D}}'=\{A+\delta:A\in {\cal
D}^*,\delta\in Z_{h_1h_2n}\},$$ where $A+\delta=\{(x,y+\delta)\
({\rm mod}\ (-,h_1h_2n)):(x,y)\in A\}$. Let $H'=\{0,n,\ldots,(h_1h_2-1)n\}$ be a subgroup of order $h_1h_2$ in $Z_{h_1h_2n}$, and $H'_i=H'+i$ be a coset of $H'$
in $Z_{h_1h_2n}$, $0\leq i<n$. Take $$X'=I_{g_1g_2}\times Z_{h_1h_2n},\hspace{0.5cm} {\cal G}'=\{I_{g_1g_2}\times H'_i: 0\leq i<n\}.$$ Then $(X',{\cal G}',{{\cal
A}_1}',\ldots,{{\cal A}_s}',{\cal D}')$ is the required strictly
$(g_1g_2,h_1h_2)$-regular
$s$-FG$(3,(L_1,L_2,\ldots,L_s,$\\ $L_T),g_1g_2h_1h_2n)$ of type
$(g_1g_2h_1h_2)^n$.
\end{algorithm}}

\vspace{0.5cm}

The following example illustrates the algorithm presented in Figure $3$.

\begin{Example}
\label{example-(4,2)-regular-16^2} In this example, we construct an optimal
 $2$-D $(8\times 4,4,2)$-OOC.

\begin{itemize}
\item Step $1$: First we construct a strictly $(2,2)$-regular $1$-FG$(3,(2,4),8)$ of type $4^2$ as follows.

$(1)$ Take a strictly $(2,2)$-regular $0$-FG$(3,(\emptyset,4),8)$ of type $4^2$ $(X,{\cal G},\emptyset,{\cal T})$ on $X=I_2\times Z_4$ with the group set ${\cal G}=\{I_2\times H_i:0\leq i<2\}$, where $H_0=\{0,2\}$ is a subgroup of order $2$ in $Z_4$ and $H_1=\{1,3\}$. It exists by Example \ref{example-4-cyclic-4^2}. Denote the family of base blocks of this design by ${\cal F}_2$. It follows that ${\cal F}_2$ generates all the blocks of ${\cal T}$ and $|{\cal F}_2|=3$.

$(2)$ Collect all $2$-subsets from distinct groups of $\cal G$ into a set $\cal B$. Then $(X,{\cal G}\cup{\cal B})$ is a $2$-design. Hence, $(X,{\cal G},{\cal B},{\cal T})$ is a strictly $(2,2)$-regular $1$-FG$(3,(2,4),8)$ of type $4^2$. Take {\footnotesize$${\cal F}_1=\{\{(0,0),(0,1)\},\{(0,0),(1,1)\},\{(1,0),(1,1)\},\{(0,0),(1,3)\}\}.$$}
${\cal F}_1$ generates all the blocks of ${\cal B}$.

\item Step $2$: For each $B\in {\cal F}_1$, construct a strictly $1$-cyclic $0$-FG$(3,(\emptyset,4),8)$ of type $4^2$ on $B\times I_4\times Z_1$ with the group set $\{\{x\}\times I_4\times Z_1: x\in B\}$, which can be taken from Example \ref{0-fan(4^2)}. Denote the family of base blocks of this design by ${\cal D}_B$, and $|{\cal D}_B|=12$.

\item Step $3$: For each $B\in {\cal F}_2$, construct a $1$-cyclic $H(4,4,4,3)$
on $B\times I_4\times Z_1$ with the
group set $\{\{x\}\times I_4\times Z_1: x\in B\}$, which exists by Corollary \ref{cor for $h$-cyclic $H$ design}.
Denote
the family of base blocks of this design by ${{\cal D}'_B}$, and $|{\cal D}'_B|=64$.

\item Step $4$: Let ${\cal D} = (\bigcup_{B\in {\cal F}_1} {\cal D}_B)\bigcup (\bigcup_{B\in {\cal F}_2} {{\cal D}'_B}).$ For each $C\in{\cal D}$ and each $(x,y,z,u)\in C$, define a mapping $\tau: (x,y,z,u) \longmapsto (x+2z,y+4u).$ Define $\tau(C)=\{\tau(c):c\in C\}$. Let ${\cal D}^*=\bigcup_{_{C\in {\cal D}}} \tau(C)$. Then $|{\cal D}^*|=240$, which is just the number of base blocks in a strictly $(8,2)$-regular $0$-FG$(3,(\emptyset,4),32)$ of type $16^2$.

\item Step $5$: Let ${{\cal D}}'=\{A+\delta:A\in {\cal D}^*,\delta\in Z_4\},$ where $A+\delta=\{(x,y+\delta)\ ({\rm mod}\ (-,4)):(x,y)\in A\}$. Let $H'=\{0,2\}$ be a subgroup of order $2$ in $Z_4$, and $H'_1=\{1,3\}$. Take $X'=I_8\times Z_4$ and ${\cal G}'=\{I_8\times H'_i: 0\leq i<2\}.$ Then $(X',{\cal G}',\emptyset,{\cal D}')$ is the required strictly $(8,2)$-regular $0$-FG$(3,(\emptyset,4),32)$ of type $16^2$.

\item Step $6$: Apply Construction $\ref{filling III}$. Fill in the groups of the resulting strictly $(8,2)$-regular $0$-FG$(3,(\emptyset,4),32)$ of type $16^2$ with an optimal strictly $2$-cyclic $3$-$(8\times 2, 4,1)$ packing with $68$ base blocks, which exists by Example \ref{example-2-cyclic-8^2}. We have an optimal strictly $4$-cyclic $3$-$(8\times 4, 4,1)$ packing with $308$ base blocks, which achieves the upper bound in Theorem $\ref{bound-total}$ and is an optimal $2$-D $(8\times 4,4,2)$-OOC with $308$ codewords. Hence $\Phi(8\times 4,4,2)=J^*(8\times 4)=308$.
\end{itemize}
\end{Example}

 {\footnotesize
\begin{algorithm}[H]
\caption{Algorithm in Construction \ref{Weighting for $h$-cyclic $H$ design}}

{\rm \textbf{Step $1$}}: Start from an $h_1$-cyclic $H(n,g_1h_1,K,t)$ $(X,{\cal G},{\cal B})$, where $X=I_n\times I_{g_1}\times Z_{h_1}$ and
${\cal G}=\{\{x\}\times I_{g_1}\times Z_{h_1} : x\in I_n\}$.
\begin{itemize}
\item Denote
the family of base blocks of this design by ${\cal F}$.
\end{itemize}

{\rm \textbf{Step $2$ (input)}}: For any base block $B\in {\cal F}$, construct an $h_2$-cyclic $H(|B|,g_2h_2,L,t)$ on $B\times I_{g_2}\times Z_{h_2}$ with the
group set $\{\{x\}\times I_{g_2}\times Z_{h_2}: x\in B\}$.
\begin{itemize}
\item Denote
the family of base blocks of this design by ${\cal D}_B$.
\end{itemize}

{\rm \textbf{Step $3$ (mapping)}}: Let ${\cal
D}=\bigcup_{B\in {\cal F}} {\cal D}_B.$ For each $C\in {\cal D}$ and
each $(x,y,z,u,v)\in C$, define a mapping
$$\tau: \hspace{0.5cm} (x,y,z,u,v) \longmapsto (x,y+ug_1,z+vh_1).$$ Define $\tau(C)=\{\tau(c):c\in C\}$.
Let ${\cal D}^*=\bigcup_{_{C\in {\cal D}}}
\tau(C).$

{\rm \textbf{Step $4$ (final)}}: Take $${\cal D}'=\{D+\delta:D\in {\cal
D}^*,\delta\in Z_{h_1h_2}\},$$ where $D+\delta=\{(x,y,z+\delta)\
({\rm mod}\ (-,-,h_1h_2)):(x,y,z)\in D\}$. Take $$X'=I_n\times
I_{g_1g_2}\times Z_{h_1h_2},\hspace{0.5cm} {\cal G}'=\{\{x\}\times
I_{g_1g_2}\times Z_{h_1h_2}: x\in I_n\}.$$ Then $(X',{\cal G}',{\cal D}')$ is the required $h_1h_2$-cyclic $H(n,g_1g_2h_1h_2,L,t)$.
\end{algorithm}}

\begin{Construction}
\label{Weighting for $h$-cyclic $H$ design} $(${\rm Weighting Construction-III}$)$ Suppose that the following exist:

\begin{enumerate}
\item[$(1)$] an $h_1$-cyclic $H(n,g_1h_1,K,t)$;

\item[$(2)$] an $h_2$-cyclic $H(k,g_2h_2,L,t)$ for each $k\in K$.
\end{enumerate}

\noindent Then there exists an $h_1h_2$-cyclic $H(n,g_1g_2h_1h_2,$ $L,t)$.
\end{Construction}

\begin{Corollary}
\label{cor for $h$-cyclic $H$ design} For any $h\geq 1$ and $n\geq
4$, $n\neq 5$, if $gn$ is even and $g(n-1)(n-2)$ is divisible by
$3$, then there is an $h$-cyclic $H(n,gh,4,3)$. For any $h\geq 1$
and $n=5$, an $h$-cyclic $H(5,gh,4,3)$ exists if $g$ is even, $g\neq
2$ and $g\not\equiv 10,26$ $({\rm mod}$ $48)$.
\end{Corollary}

\proof By Lemma \ref{H(4,g,4,3)}, for any $h\geq 1$, there exists a
semi-cyclic $H(4,h,4,3)$ (i.e., an $h$-cyclic $H(4,h,4,3)$). Apply Construction \ref{Weighting for $h$-cyclic $H$ design} with $h_1=g_2=1$, $g_1=g$ and $h_2=h$. Combine the results of Lemma \ref{H design} to complete the proof. \qed

\section{Small orders of optimal $2$-D $(u\times v,4,2)$-OOCs}

In this section, we obtain some small orders of optimal $2$-D OOCs. Some of them are obtained by computer search, and some of them are obtained by applying filling constructions in Section $6$.

\begin{Lemma}
\label{small order} There exists an optimal $2$-D $(u\times
v,4,2)$-OOC with $J^*(u\times v)$ codewords for each
$(u,v)\in\{(3,3)$, $(2,6)$, $(3,4)$, $(6,2)$, $(7,2)$, $(2,11)\}$.
\end{Lemma}

\proof We here give a construction of a $3$-$(uv,4,1)$-packing on $I_{uv}$. Let
$\alpha=(0\ 1\ \cdots \ v-1)(v\ v+1\ \cdots \ 2v-1)\cdots((u-1)v\ \
\cdots\ uv-1)$ be a permutation on $I_{uv}$, which consists of $u$
cycles of length $v$. Let $G$ be the group generated by $\alpha$.
Only base blocks are listed below. All other blocks are obtained by
developing these base blocks under the action of $G$. Obviously this
design is isomorphic to a strictly $v$-cyclic $3$-$(u\times
v,4,1)$-packing, which achieves the upper bound in
 Theorem \ref{bound-total}, and is an optimal $2$-D $(u\times
v,4,2)$-OOC.

\begin{center}
{\footnotesize \tabcolsep 0.03in \scriptsize
\begin{tabular}{lllllll}
$(u,v)=(3,3)$:& $\{0,1,3,4\}$& $\{0,1,5,6\}$& $\{0,1,7,8\}$&
$\{0,3,6,8\}$& $\{0,4,5,7\}$& $\{3,4,7,8\}$
\\$(u,v)=(2,6)$:&$\{0,1,2,6\}$& $\{0,1,3,8\}$& $\{0,1,4,7\}$& $\{0,1,9,10\}$&
$\{0,2,8,10\}$& $\{0,6,7,9\}$\\ & $\{0,6,10,11\}$& $\{0,7,8,11\}$
\\$(u,v)=(3,4)$:&$\{0,2,4,5\}$& $\{0,1,2,8\}$& $\{0,1,4,9\}$& $\{0,1,5,7\}$&
$\{0,1,6,10\}$& $\{0,4,7,10\}$\\ & $\{0,4,8,11\}$& $\{0,5,6,11\}$&
$\{0,5,8,10\}$& $\{0,6,8,9\}$& $\{0,7,9,11\}$& $\{4,5,6,9\}$
\\$(u,v)=(6,2)$:&$\{0,1,2,4\}$& $\{0,1,6,8\}$& $\{0,2,5,9\}$& $\{0,2,3,6\}$&
$\{0,2,7,8\}$& $\{0,2,10,11\}$\\ &$\{0,3,4,7\}$& $\{0,3,8,9\}$&
$\{0,4,5,6\}$& $\{0,4,8,10\}$& $\{0,4,9,11\}$& $\{0,5,7,10\}$\\&
$\{0,5,8,11\}$& $\{0,6,7,11\}$& $\{0,6,9,10\}$& $\{2,3,4,9\}$&
$\{2,4,5,10\}$& $\{2,4,6,7\}$\\& $\{2,4,8,11\}$& $\{2,5,7,11\}$&
$\{2,6,8,10\}$& $\{2,6,9,11\}$& $\{2,7,9,10\}$& $\{4,6,8,9\}$\\&
$\{4,7,10,11\}$
\\$(u,v)=(7,2)$:& $\{0,1,2,4\}$&
$\{0,1,6,8\}$& $\{0,1,10,12\}$& $\{0,2,3,7\}$& $\{0,2,5,10\}$&
$\{0,2,6,13\}$\\& $\{0,2,8,11\}$& $\{0,2,9,12\}$& $\{0,3,4,6\}$&
$\{0,3,8,10\}$& $\{0,3,9,11\}$& $\{0,3,12,13\}$\\& $\{0,4,5,13\}$&
$\{0,4,7,10\}$& $\{0,4,8,9\}$& $\{0,4,11,12\}$& $\{0,5,6,9\}$&
$\{0,5,7,11\}$\\& $\{0,5,8,12\}$& $\{0,6,7,12\}$& $\{0,6,10,11\}$&
$\{0,7,8,13\}$& $\{0,9,10,13\}$& $\{2,3,4,9\}$\\& $\{2,3,10,12\}$&
$\{2,4,5,6\}$& $\{2,4,7,12\}$& $\{2,4,8,13\}$& $\{2,4,10,11\}$&
$\{2,5,9,13\}$\\& $\{2,5,11,12\}$& $\{2,6,7,11\}$& $\{2,6,8,12\}$&
$\{2,6,9,10\}$& $\{2,7,8,9\}$& $\{2,7,10,13\}$\\& $\{4,5,8,10\}$&
$\{4,6,7,9\}$& $\{4,6,8,11\}$& $\{4,6,12,13\}$& $\{4,7,11,13\}$&
$\{4,9,10,12\}$\\& $\{6,8,10,13\}$& $\{6,9,11,13\}$
\\$(u,v)=(2,11)$:&$\{0,1,2,4\}$&
$\{0,1,5,7\}$& $\{0,1,6,9\}$& $\{0,1,8,11\}$& $\{0,1,12,13\}$&
$\{0,1,14,15\}$\\& $\{0,1,16,17\}$& $\{0,1,18,19\}$&
$\{0,1,20,21\}$& $\{0,2,5,11\}$& $\{0,2,12,14\}$& $\{0,2,13,15\}$\\&
$\{0,2,16,19\}$& $\{0,2,17,21\}$& $\{0,2,18,20\}$& $\{0,3,7,12\}$&
$\{0,3,11,15\}$& $\{0,3,13,16\}$\\& $\{0,3,17,19\}$&
$\{0,3,18,21\}$& $\{0,4,11,18\}$& $\{0,4,12,15\}$& $\{0,4,13,17\}$&
$\{0,4,19,21\}$\\& $\{0,5,12,17\}$& $\{0,5,13,18\}$&
$\{0,5,14,19\}$& $\{0,5,15,20\}$& $\{0,5,16,21\}$&
$\{0,11,14,21\}$\\& $\{0,11,17,20\}$& $\{0,12,16,20\}$&
$\{11,12,13,17\}$& $\{11,12,14,19\}$& $\{11,12,18,20\}$
\end{tabular}}
\end{center}

\begin{Lemma}
\label{strictly $(2,6)$-regular of type $12^2$} There exists a
strictly $(2,6)$-regular $0$-FG$(3,(\emptyset,4),24)$ of type
$12^2$.
\end{Lemma}

\proof We here give a construction of a $0$-FG$(3,(\emptyset,4),24)$ of type $12^2$ on
$I_{24}$ with the group set $\{\{2i+j: 0\leq i\leq 11\} :$ $0\leq
j\leq 1\}$. Let $\alpha=(0\ 1\ \cdots \ 11)(12\ 13\ \cdots \ 23)$ be
a permutation on $I_{24}$, which consists of $2$ cycles of length
$12$. Let $G$ be the group generated by $\alpha$. Only base blocks
are listed below. All other blocks are obtained by developing these
base blocks under the action of $G$. Obviously this design is
isomorphic to a strictly $(2,6)$-regular
$0$-FG$(3,(\emptyset,4),24)$ of type $12^2$.

\begin{center}
{\footnotesize \tabcolsep 0.03in \scriptsize
\begin{tabular}{lllllll}
$\{0,1,2,5\}$& $\{0,1,3,8\}$& $\{0,1,6,9\}$& $\{0,1,7,12\}$&
$\{0,1,10,13\}$& $\{0,1,14,15\}$& $\{0,1,16,17\}$\\
$\{0,1,18,19\}$& $\{0,1,20,21\}$& $\{0,1,22,23\}$& $\{0,2,9,13\}$&
$\{0,2,17,19\}$& $\{0,2,21,23\}$& $\{0,3,12,19\}$\\
$\{0,3,13,18\}$& $\{0,3,14,21\}$& $\{0,3,17,22\}$& $\{0,3,20,23\}$&
$\{0,4,13,19\}$& $\{0,4,15,23\}$& $\{0,4,17,21\}$\\ $\{0,5,12,21\}$&
$\{0,5,13,20\}$& $\{0,5,14,23\}$& $\{0,5,15,18\}$& $\{0,5,19,22\}$&
$\{0,6,13,21\}$& $\{0,12,13,15\}$\\ $\{0,12,17,23\}$&
$\{0,13,16,23\}$& $\{12,13,14,21\}$& $\{12,13,16,19\}$&
$\{12,13,17,22\}$
\end{tabular}}
\end{center}

\begin{Lemma}
\label{2*12} There exists an optimal $2$-D $(2\times 12,4,2)$-OOC
with $J^*(2\times 12)=41$ codewords.
\end{Lemma}

\proof Start from a strictly $(2,6)$-regular
$0$-FG$(3,(\emptyset,4),24)$ of type $12^2$, which exists by Lemma
\ref{strictly $(2,6)$-regular of type $12^2$}. Apply Construction
$\ref{filling III}$ with an optimal strictly $6$-cyclic
 $3$-$(2\times 6, 4,1)$ packing from Lemma \ref{small order} to obtain a strictly $12$-cyclic
 $3$-$(2\times 12, 4,1)$ packing with $41$ base blocks, which achieves the upper bound in
 Theorem \ref{bound-total}, and is an optimal
 $2$-D $(2\times 12,4,2)$-OOC with $41$ codewords. \qed

\begin{Lemma}
\label{example-2-cyclic-12^2} There exists a strictly $2$-cyclic
$0$-FG$(3,(\emptyset,4),24)$ of type $12^2$.
\end{Lemma}

\proof We here give a construction of a $0$-FG$(3,(\emptyset,4),24)$ of type $12^2$ on
$I_{24}$ with the group set $\{\{0,1,\ldots,11\}$ $+i:i\in
\{0,12\}\}$. Let $\alpha=(0\ 1)(2\ 3)\ \cdots\ (22\ 23)$ be a
permu-tation on $I_{24}$ and $G$ be the group generated by $\alpha$.
Only base blocks are listed below. All other blocks are obtained by
developing these base blocks under the action of $G$. Obviously this
design is isomorphic to a strictly $2$-cyclic
$0$-FG$(3,(\emptyset,4),24)$ of type $12^2$.

\begin{center}
{\footnotesize \tabcolsep 0.03in \scriptsize
\begin{tabular}{lllllll}
$\{0,1,12,14\}$& $\{0,1,16,18\}$& $\{0,1,20,22\}$& $\{0,2,12,13\}$&
$\{0,2,14,15\}$& $\{0,2,16,17\}$& $\{0,2,18,19\}$\\ $\{0,2,20,21\}$&
$\{0,2,22,23\}$& $\{0,3,12,15\}$& $\{0,3,13,14\}$& $\{0,3,16,19\}$&
$\{0,3,17,18\}$& $\{0,3,20,23\}$\\ $\{0,3,21,22\}$& $\{0,4,12,16\}$&
$\{0,4,13,17\}$& $\{0,4,14,20\}$&
$\{0,4,15,21\}$& $\{0,4,18,22\}$& $\{0,4,19,23\}$\\
$\{0,5,12,17\}$& $\{0,5,13,20\}$& $\{0,5,14,22\}$&
$\{0,5,15,18\}$& $\{0,5,16,23\}$& $\{0,5,19,21\}$& $\{0,6,12,19\}$\\
$\{0,6,13,16\}$& $\{0,6,14,21\}$& $\{0,6,15,20\}$& $\{0,6,17,22\}$&
$\{0,6,18,23\}$& $\{0,7,12,20\}$& $\{0,7,13,21\}$\\ $\{0,7,14,18\}$&
$\{0,7,15,19\}$& $\{0,7,16,22\}$& $\{0,7,17,23\}$& $\{0,8,12,18\}$&
$\{0,8,13,19\}$& $\{0,8,14,23\}$\\ $\{0,8,15,22\}$& $\{0,8,16,20\}$&
$\{0,8,17,21\}$& $\{0,9,12,21\}$&
$\{0,9,13,18\}$& $\{0,9,14,16\}$& $\{0,9,15,23\}$\\
$\{0,9,17,20\}$& $\{0,9,19,22\}$& $\{0,10,12,22\}$&
$\{0,10,13,23\}$& $\{0,10,14,17\}$& $\{0,10,15,16\}$&
$\{0,10,18,21\}$\\ $\{0,10,19,20\}$& $\{0,11,12,23\}$&
$\{0,11,13,22\}$& $\{0,11,14,19\}$&
$\{0,11,15,17\}$& $\{0,11,16,21\}$& $\{0,11,18,20\}$\\
$\{2,3,12,14\}$& $\{2,3,16,18\}$& $\{2,3,20,22\}$& $\{2,4,12,17\}$&
$\{2,4,13,18\}$& $\{2,4,14,21\}$& $\{2,4,15,20\}$\\ $\{2,4,16,23\}$&
$\{2,4,19,22\}$& $\{2,5,12,18\}$& $\{2,5,13,19\}$& $\{2,5,14,23\}$&
$\{2,5,15,22\}$& $\{2,5,16,20\}$\\ $\{2,5,17,21\}$& $\{2,6,12,16\}$&
$\{2,6,13,17\}$& $\{2,6,14,20\}$&
$\{2,6,15,21\}$& $\{2,6,18,22\}$& $\{2,6,19,23\}$\\
$\{2,7,12,19\}$& $\{2,7,13,20\}$& $\{2,7,14,22\}$&
$\{2,7,15,17\}$& $\{2,7,16,21\}$& $\{2,7,18,23\}$& $\{2,8,12,22\}$\\
$\{2,8,13,23\}$& $\{2,8,14,16\}$& $\{2,8,15,18\}$& $\{2,8,17,20\}$&
$\{2,8,19,21\}$& $\{2,9,12,23\}$& $\{2,9,13,22\}$\\ $\{2,9,14,17\}$&
$\{2,9,15,16\}$& $\{2,9,18,21\}$& $\{2,9,19,20\}$& $\{2,10,12,21\}$&
$\{2,10,13,16\}$& $\{2,10,14,19\}$\\ $\{2,10,15,23\}$&
$\{2,10,17,22\}$& $\{2,10,18,20\}$& $\{2,11,12,20\}$&
$\{2,11,13,21\}$& $\{2,11,14,18\}$& $\{2,11,15,19\}$\\
$\{2,11,16,22\}$& $\{2,11,17,23\}$& $\{4,5,12,14\}$&
$\{4,5,16,18\}$& $\{4,5,20,22\}$& $\{4,6,12,22\}$& $\{4,6,13,23\}$\\
$\{4,6,14,17\}$& $\{4,6,15,16\}$& $\{4,6,18,21\}$& $\{4,6,19,20\}$&
$\{4,7,12,23\}$& $\{4,7,13,22\}$& $\{4,7,14,16\}$\\ $\{4,7,15,18\}$&
$\{4,7,17,20\}$& $\{4,7,19,21\}$& $\{4,8,12,13\}$& $\{4,8,14,15\}$&
$\{4,8,16,17\}$& $\{4,8,18,19\}$\\ $\{4,8,20,21\}$& $\{4,8,22,23\}$&
$\{4,9,12,20\}$& $\{4,9,13,21\}$&
$\{4,9,14,18\}$& $\{4,9,15,19\}$& $\{4,9,16,22\}$\\
$\{4,9,17,23\}$& $\{4,10,12,19\}$& $\{4,10,13,20\}$&
$\{4,10,14,22\}$& $\{4,10,15,17\}$& $\{4,10,16,21\}$&
$\{4,10,18,23\}$\\ $\{4,11,12,15\}$& $\{4,11,13,14\}$&
$\{4,11,16,19\}$& $\{4,11,17,18\}$&
$\{4,11,20,23\}$& $\{4,11,21,22\}$& $\{6,7,12,14\}$\\
$\{6,7,16,18\}$& $\{6,7,20,22\}$& $\{6,8,12,15\}$& $\{6,8,13,14\}$&
$\{6,8,16,19\}$& $\{6,8,17,18\}$& $\{6,8,20,23\}$\\ $\{6,8,21,22\}$&
$\{6,9,12,17\}$& $\{6,9,13,20\}$& $\{6,9,14,22\}$& $\{6,9,15,18\}$&
$\{6,9,16,23\}$& $\{6,9,19,21\}$\\ $\{6,10,12,13\}$&
$\{6,10,14,15\}$& $\{6,10,16,17\}$& $\{6,10,18,19\}$&
$\{6,10,20,21\}$& $\{6,10,22,23\}$& $\{6,11,12,18\}$\\
$\{6,11,13,19\}$& $\{6,11,14,23\}$& $\{6,11,15,22\}$&
$\{6,11,16,20\}$& $\{6,11,17,21\}$& $\{8,9,12,14\}$&
$\{8,9,16,18\}$\\ $\{8,9,20,22\}$& $\{8,10,12,17\}$&
$\{8,10,13,18\}$& $\{8,10,14,20\}$& $\{8,10,15,21\}$&
$\{8,10,16,23\}$& $\{8,10,19,22\}$\\ $\{8,11,12,16\}$&
$\{8,11,13,17\}$& $\{8,11,14,21\}$& $\{8,11,15,20\}$&
$\{8,11,18,22\}$& $\{8,11,19,23\}$& $\{10,11,12,14\}$\\
$\{10,11,16,18\}$& $\{10,11,20,22\}$
\end{tabular}}
\end{center}

\begin{Lemma}
\label{12*2} There exists an optimal $2$-D $(12\times 2,4,2)$-OOC
with $J^*(12\times 2)=248$ codewords.
\end{Lemma}

\proof Start from a strictly $2$-cyclic $0$-FG$(3,(\emptyset,4),24)$
of type $12^2$, which exists by Lemma \ref{example-2-cyclic-12^2}.
Applying Construction \ref{filling II} with an optimal strictly
$2$-cyclic $3$-$(6\times 2,4,1)$ packing from Lemma \ref{small
order}, we have a strictly $2$-cyclic $3$-$(12\times 2,4,1)$ packing
with $248$ base blocks. This number achieves the upper bound in
Theorem \ref{bound-total}. Thus an optimal $2$-D $(12\times
2,4,2)$-OOC with $248$ codewords exists. \qed

\begin{Lemma}
\label{strictly $(2,3)$-regular of type $6^5$} There exists a
strictly $(2,3)$-regular $0$-FG$(3,(\emptyset,4),30)$ of type $6^5$.
\end{Lemma}

\proof We here give a construction of a $0$-FG$(3,(\emptyset,4),30)$ of type $6^5$ on
$I_{30}$ with the group set $\{\{5i+j: 0\leq i\leq 5\} :$ $0\leq
j\leq 4\}$. Let $\alpha=(0\ 1\ \cdots \ 14)(15\ 16\ \cdots \ 29)$ be
a permutation on $I_{30}$, which consists of $2$ cycles of length
$15$. Let $G$ be the group generated by $\alpha$. Only base blocks
are listed below. All other blocks are obtained by developing these
base blocks under the action of $G$. Obviously this design is
isomorphic to a strictly $(2,3)$-regular
$0$-FG$(3,(\emptyset,4),30)$ of type $6^5$.

\begin{center}
{\footnotesize \tabcolsep 0.03in \scriptsize
\begin{tabular}{lllllll}
$\{0,1,2,4\}$& $\{0,1,5,6\}$& $\{0,1,7,9\}$& $\{0,1,8,12\}$&
$\{0,1,13,15\}$& $\{0,1,16,17\}$& $\{0,1,18,19\}$\\ $\{0,1,20,21\}$&
$\{0,1,22,23\}$& $\{0,1,24,25\}$& $\{0,1,26,27\}$& $\{0,1,28,29\}$&
$\{0,2,5,12\}$& $\{0,2,6,11\}$\\ $\{0,2,7,15\}$& $\{0,2,10,16\}$&
$\{0,2,18,20\}$& $\{0,2,19,21\}$& $\{0,2,22,24\}$& $\{0,2,23,25\}$&
$\{0,2,26,28\}$\\ $\{0,2,27,29\}$& $\{0,3,6,15\}$& $\{0,3,7,16\}$&
$\{0,3,9,19\}$& $\{0,3,18,21\}$& $\{0,3,20,23\}$& $\{0,3,22,26\}$\\
$\{0,3,24,29\}$& $\{0,3,25,28\}$& $\{0,4,15,19\}$& $\{0,4,16,20\}$&
$\{0,4,17,24\}$& $\{0,4,18,22\}$& $\{0,4,21,27\}$\\ $\{0,4,23,26\}$&
$\{0,4,25,29\}$& $\{0,5,16,22\}$& $\{0,5,17,27\}$& $\{0,5,18,23\}$&
$\{0,5,19,24\}$& $\{0,5,26,29\}$\\ $\{0,6,17,22\}$& $\{0,6,18,24\}$&
$\{0,6,19,28\}$& $\{0,6,20,26\}$& $\{0,6,21,29\}$& $\{0,6,23,27\}$&
$\{0,7,17,25\}$\\ $\{0,7,18,26\}$& $\{0,7,19,27\}$& $\{0,7,20,24\}$&
$\{0,7,22,29\}$& $\{0,7,23,28\}$& $\{0,15,17,29\}$&
$\{0,15,21,28\}$\\
$\{0,15,24,27\}$& $\{0,16,19,25\}$& $\{0,16,23,29\}$&
$\{0,16,24,28\}$& $\{15,16,17,28\}$& $\{15,16,19,23\}$&
$\{15,16,20,22\}$\\ $\{15,16,21,25\}$& $\{15,16,24,26\}$&
$\{15,17,20,27\}$
\end{tabular}}
\end{center}

\begin{Lemma}
\label{2*15} There exists an optimal $2$-D $(2\times 15,4,2)$-OOC
with $J^*(2\times 15)=67$ codewords.
\end{Lemma}

\proof Start from a strictly $(2,3)$-regular
$0$-FG$(3,(\emptyset,4),30)$ of type $6^5$, which exists by Lemma
\ref{strictly $(2,3)$-regular of type $6^5$}. Apply Construction
$\ref{filling III}$ with an optimal strictly $3$-cyclic
 $3$-$(2\times 3, 4,1)$ packing from Example \ref{example-6} to obtain a strictly $15$-cyclic
 $3$-$(2\times 15, 4,1)$ packing with $67$ base blocks, which  achieves the upper bound in
 Theorem \ref{bound-total}, and is an optimal
 $2$-D $(2\times 15,4,2)$-OOC with $67$ codewords. \qed

 \begin{Lemma}
\label{strictly $(3,2)$-regular of type $6^5$} There exists a
strictly $(3,2)$-regular $0$-FG$(3,(\emptyset,4),30)$ of type $6^5$.
\end{Lemma}

\proof We here give a construction of a $0$-FG$(3,(\emptyset,4),30)$ of type $6^5$ on
$I_{30}$ with the group set $\{\{5i+j: 0\leq i\leq 5\} :$ $0\leq
j\leq 4\}$. Let $\alpha=(0\ 1\ \cdots \ 9)(10\ 11\ \cdots \ 19)(20\
21\ \cdots \ 29)$ and $\beta=(0\ 10\ 20)(1\ 11\ 21) \cdots (9\ 19\
29)$ be two permutations on $I_{30}$ and $G$ be the group generated
by $\alpha$ and $\beta$. Only base blocks are listed below. All
other blocks are obtained by developing these base blocks under the
action of $G$. Obviously this design is isomorphic to a
$(3,2)$-regular $0$-FG$(3,(\emptyset,4),30)$ of type $6^5$.

\begin{center}
{\footnotesize \tabcolsep 0.03in \scriptsize
\begin{tabular}{llllllll}
$\{0,1,2,4\}$& $\{0,1,5,7\}$& $\{0,1,6,10\}$& $\{0,1,8,11\}$&
$\{0,1,12,13\}$& $\{0,1,14,15\}$& $\{0,1,16,17\}$& $\{0,1,18,20\}$\\
$\{0,1,19,22\}$& $\{0,1,21,23\}$& $\{0,2,7,26\}$& $\{0,2,11,22\}$&
$\{0,2,12,18\}$& $\{0,2,14,20\}$& $\{0,2,15,17\}$& $\{0,2,16,21\}$\\
$\{0,2,23,29\}$& $\{0,2,24,28\}$& $\{0,3,6,24\}$& $\{0,3,10,14\}$&
$\{0,3,11,18\}$& $\{0,3,12,27\}$& $\{0,3,15,23\}$& $\{0,3,16,20\}$\\
$\{0,3,17,28\}$& $\{0,3,19,26\}$& $\{0,4,10,21\}$& $\{0,4,12,26\}$&
$\{0,4,15,19\}$& $\{0,4,17,27\}$& $\{0,5,11,27\}$& $\{0,5,12,21\}$\\
$\{0,5,13,28\}$
\end{tabular}}
\end{center}

\begin{Lemma}
\label{3*10} There exists an optimal $2$-D $(3\times 10,4,2)$-OOC
with $J^*(3\times 10)=100$ codewords.
\end{Lemma}

\proof Start from a strictly $(3,2)$-regular
$0$-FG$(3,(\emptyset,4),30)$ of type $6^5$, which exists by Lemma
\ref{strictly $(3,2)$-regular of type $6^5$}. Apply Construction
$\ref{filling III}$ with an optimal strictly $2$-cyclic
 $3$-$(3\times 2, 4,1)$ packing from Example \ref{example-6} to obtain a strictly $10$-cyclic
 $3$-$(3\times 10, 4,1)$ packing with $100$ base blocks, which  achieves the upper bound in
 Theorem \ref{bound-total}, and is an optimal
 $2$-D $(3\times 10,4,2)$-OOC with $100$ codewords. \qed

\begin{Lemma}
\label{2D from 1D 18} Let $n\equiv 18\ ({\rm mod }\ 24)$. If there
is an optimal $1$-D $(n,4,2)$-OOC, which achieves the Johnson bound
$J(1\times
n,4,2)=\lfloor\frac{1}{4}\lfloor\frac{n-1}{3}\lfloor\frac{n-2}{2}
\rfloor\rfloor\rfloor$, then for any integer factorization $n=uv$,
there is an optimal $2$-D $(u\times v,4,2)$-OOC with $J^*(u\times
v)=\lfloor\frac{u}{4}(\lfloor\frac{n-1}{3}\lfloor\frac{n-2}{2}
\rfloor\rfloor-1)\rfloor$ codewords.
\end{Lemma}

\proof By Corollary \ref{relation1}, if there exists an optimal $1$-D
$(uv,4,2)$-OOC with $J(1\times uv,4,2)$ codewords, then there exists
a $2$-D $(u\times v,4,2)$-OOC with $uJ(1\times uv,4,2)$ codewords.
It is readily checked that $uJ(1\times
uv,4,2)=u(u^2v^2-3uv-6)/24=\lfloor\frac{u}{4}(\lfloor\frac{uv-1}{3}\lfloor\frac{uv-2}{2}
\rfloor\rfloor-1)\rfloor=J^*(u\times v)$. This number achieves the
upper bound in Theorem \ref{bound-total}. This completes the proof.
\qed

Note that when $n\equiv 18\ ({\rm mod }\ 24)$, $J(1\times
n,4,2)=\lfloor\frac{1}{4}\lfloor\frac{n-1}{3}\lfloor\frac{n-2}{2}
\rfloor\rfloor\rfloor=\lfloor\frac{1}{4}(\lfloor\frac{n-1}{3}\lfloor\frac{n-2}{2}
\rfloor\rfloor-1)\rfloor$. Hence no confusion occurs in Lemma
\ref{2D from 1D 18}. By Lemma \ref{1D k=4}$(3)$, there is an optimal
$1$-D $(n,4,2)$-OOC with $J(1\times n,4,2)$ codewords for each
$n\in\{18,42,90\}$. Then we have

\begin{Corollary}
\label{18,42,90} Let $n\in\{18,42,90\}$. For any integer
factorization $n=n_1n_2$, there is an optimal $2$-D $(n_1\times
n_2,4,2)$-OOC with $J^*(n_1\times n_2)$ codewords.
\end{Corollary}

\begin{Lemma}
\label{new small order-u*v derived} There exists an optimal $2$-D
$(u\times v,4,2)$-OOC with $J^*(u\times v)$ codewords for each
$(u,v)\in\{(5,4)$, $(7,4)$, $(6,5)\}$.
\end{Lemma}

\proof Apply Theorem \ref{relation2} with some known optimal $2$-D
$(u_1\times v_1,4,2)$-OOCs. One can have all the required optimal
$2$-D $(u\times v,4,2)$-OOCs. For illustrating the details, we give
the following table.

{\footnotesize \tabcolsep 0.05in
\begin{center}
\begin{tabular}{|ccc|l|cc|c|}\hline
$(u_1,v_1)$ & Source & number of codewords &$\Rightarrow$ & $(u,v)$
& number of codewords & $J^*(u\times v)$
\\\hline $(1,20)$& Lemma \ref{1D k=4} & $14$ & & $(5,4)$ & $70$& $70$
\\\hline $(1,28)$& Lemma \ref{1D k=4} & $29$ & & $(7,4)$ & $203$& $203$
\\\hline $(2,15)$& Lemma \ref{2*15} & $67$ & & $(6,5)$ & $201$& $201$
\\\hline
\end{tabular} \end{center} }

\section{Infinite families of optimal $2$-D $(u\times v,4,2)$-OOCs}

In this section, on one hand we shall give some infinite families of optimal $2$-D $(u\times v,4,2)$-OOCs, which will be presented as Theorems. On the other hand, although we can not complete the existence of optimal $2$-D $(u\times v,4,2)$-OOCs, we hope to present some possible approaches to complete it, which will be presented as Propositions.

\begin{Lemma}
\label{infinite family-u-2-1} There exists an optimal $2$-D
$(u\times 2,4,2)$-OOC with $J^*(u\times 2)$ codewords for any
$u\equiv 2,4$ $({\rm mod}$ $6)$.
\end{Lemma}

\proof Let $n=u/2$. Then $n\equiv 1,2$ $({\rm mod}$ $3)$. When
$n=1$, an optimal $2$-D $(2\times 2,4,2)$-OOC is trivial without
base blocks. In the following consider $n\geq 2$. First we shall
show that there is a strictly $2$-cyclic
$0$-FG$(3,(\emptyset,4),4n)$ of type $4^n$ for any $n\equiv 1,2$
$({\rm mod}$ $3)$ and $n\geq 2$. When $n=2$, a strictly $2$-cyclic
$0$-FG$(3,(\emptyset,4),8)$ of type $4^2$ exists by Example
\ref{example-2-cyclic-4^2}. When $n\equiv 1,2$ $({\rm mod}$ $3)$,
$n\geq 4$ and $n\neq 5$, start from a $1$-FG$(3,(2,n),n)$ of type
$1^n$, which contains one block of size $n$ and all $2$-subsets of
$n$ points. Apply Construction \ref{Weighting-strictly $h$-cyclic
$s$-FG} with $h_1=1$ and $h_2=2$ to obtain a strictly $2$-cyclic
$0$-FG$(3,(\emptyset,4),4n)$ of type $4^n$, where the needed
$2$-cyclic $H(n,4,4,3)$ is from Corollary \ref{cor for $h$-cyclic
$H$ design}. When $n=5$, there is a strictly cyclic
$0$-FG$(3,(\emptyset,4),20)$ of type $4^5$ from Lemma \ref{strictly
cyclic g^n}. By Lemma \ref{relation between regular and v-cyclic},
it implies a strictly $2$-cyclic $0$-FG$(3,(\emptyset,4),20)$ of
type $4^5$.

Next applying Construction \ref{filling II} with an optimal strictly
$2$-cyclic $3$-$(2\times 2,4,1)$ packing, which is trivial without
base blocks, we have a strictly $2$-cyclic $3$-$(2n\times 2,4,1)$
packing, which contains
$\lfloor\frac{2n}{4}(\lfloor\frac{4n-1}{3}\lfloor\frac{4n-2}{2}
\rfloor\rfloor-1)\rfloor=n(n-1)(4n+1)/3$ base blocks. This number
achieves the upper bound in Theorem \ref{bound-total}. Thus an
optimal $2$-D $(2n\times 2,4,2)$-OOC with $J^*(2n\times 2)$
codewords exists. It is an optimal $2$-D $(u\times 2,4,2)$-OOC with
$J^*(u\times 2)$ codewords. \qed

\begin{Proposition}
\label{infinite family-u-v-1} Let $v\equiv 1$ $({\rm mod}$ $2)$ or
$v\equiv 0$ $({\rm mod}$ $4)$. Suppose that there is an optimal
$2$-D $(2\times v,4,2)$-OOC with $J^*(2\times v)$ codewords. Then
there is an optimal $2$-D $(u\times v,4,2)$-OOC with $J^*(u\times
v)$ codewords for any $u\equiv 2,4$ $({\rm mod}$ $6)$. Especially
when $v\equiv 1,5$ $({\rm mod}$ $6)$, the resulting optimal $2$-D
$(u\times v,4,2)$-OOC is perfect.
\end{Proposition}

\proof Let $n=u/2$. Then $n\equiv 1,2$ $({\rm mod}$ $3)$. When
$n=1$, the conclusion follows from the assumption. In the following
consider $n\geq 2$. First we shall show that there is a strictly
$v$-cyclic $0$-FG$(3,(\emptyset,4),2vn)$ of type $(2v)^n$ for any
$n\equiv 1,2$ $({\rm mod}$ $3)$ and $n\geq 2$. When $n=2$, a
strictly $v$-cyclic $0$-FG$(3,(\emptyset,4),4v)$ of type $(2v)^2$ is
from Corollary \ref{(2v)^2-1}. When $n\equiv 1,2$ $({\rm mod}$ $3)$,
$n\geq 4$ and $n\neq 5$, start from a $1$-FG$(3,(2,n),n)$ of type
$1^n$, which contains a block of size $n$ and all $2$-subsets of $n$
points. Apply Construction \ref{Weighting-strictly $h$-cyclic
$s$-FG} with $h_1=1$ and $h_2=v$ to obtain a strictly $v$-cyclic
$0$-FG$(3,(\emptyset,4),2vn)$ of type $(2v)^n$, where the needed
$v$-cyclic $H(n,2v,4,3)$ is from Corollary \ref{cor for $h$-cyclic
$H$ design}. When $n=5$, there is a strictly cyclic
$0$-FG$(3,(\emptyset,4),10v)$ of type $(2v)^5$ from Lemma
\ref{strictly cyclic g^n}. By Lemma \ref{relation between regular
and v-cyclic}, it implies a strictly $v$-cyclic
$0$-FG$(3,(\emptyset,4),10v)$ of type $(2v)^5$.

Next apply Construction \ref{filling II} with an optimal strictly
$v$-cyclic $3$-$(2\times v,4,1)$ packing with $J^*(2\times v)$ base
blocks, which exists by assumption. Note that by Theorem
\ref{bound-total},
$$ J^*(2\times v)=\left\{
\begin{array}{lll}
    \lfloor\frac{2}{4}\lfloor\frac{2v-1}{3}\lfloor\frac{2v-2}{2}
\rfloor\rfloor\rfloor, & \hbox{if $v\equiv 1,5$ $({\rm mod}$
$6)$,} \\
    \\
    \lfloor\frac{2}{4}(\lfloor\frac{2v-1}{3}\lfloor\frac{2v-2}{2}
\rfloor\rfloor-1)\rfloor, & \hbox{if $v\equiv 3$ $({\rm mod}$ $6)$
or $v\equiv 4,8$ $({\rm mod}$ $12)$,}
\\\\
\lfloor\frac{2}{4}(\lfloor\frac{2v-1}{3}\lfloor\frac{2v-2}{2}
\rfloor\rfloor-2)\rfloor, & \hbox{if $v\equiv 0$ $({\rm mod}$ $12)$.} \\
\end{array}
\right.$$ Then we have a strictly $v$-cyclic $3$-$(2n\times v,4,1)$
packing, which contains
$$ \left\{
\begin{array}{lll}
    n(nv-1)(2nv-1)/6=\lfloor\frac{2n}{4}\lfloor\frac{2nv-1}{3}\lfloor\frac{2nv-2}{2}
\rfloor\rfloor\rfloor, & \hbox{if $v\equiv 1,5$ $({\rm mod}$
$6)$,} \\
    \\
    n(2n^2v^2-3nv-3)/6=\lfloor\frac{2n}{4}(\lfloor\frac{2nv-1}{3}\lfloor\frac{2nv-2}{2}
\rfloor\rfloor-1)\rfloor, & \hbox{if $v\equiv 3$ $({\rm mod}$ $6)$,} \\
    \\
    n(nv-2)(2nv+1)/6=\lfloor\frac{2n}{4}(\lfloor\frac{2nv-1}{3}\lfloor\frac{2nv-2}{2}
\rfloor\rfloor-1)\rfloor, & \hbox{if $v\equiv 4,8$ $({\rm mod}$ $12)$,} \\
    \\
    n(2n^2v^2-3nv-6)/6=\lfloor\frac{2n}{4}(\lfloor\frac{2nv-1}{3}\lfloor\frac{2nv-2}{2}
\rfloor\rfloor-2)\rfloor, & \hbox{if $v\equiv 0$ $({\rm mod}$
$12)$,}
\end{array}
\right.$$ base blocks. This number achieves the upper bound in
Theorem \ref{bound-total}. Thus an optimal $2$-D $(2n\times
v,4,2)$-OOC with $J^*(2n\times v)$ codewords exists. It is an
optimal $2$-D $(u\times v,4,2)$-OOC with $J^*(u\times v)$ codewords.
Especially by Lemma \ref{nece-perfect}, when $v\equiv 1,5$ $({\rm
mod}$ $6)$, the resulting optimal $2$-D $(u\times v,4,2)$-OOC is
perfect. \qed

\begin{Lemma}
\label{infinite family-2^n} There is an optimal $2$-D $(2\times
2^n,4,2)$-OOC with $J^*(2\times 2^n)$ codewords for any positive
integer $n$.
\end{Lemma}

\proof When $n=1$, an optimal $2$-D $(2\times 2,4,2)$-OOC is trivial
without codewords. When $n=2$, the conclusion follows from Example
\ref{example-4-cyclic-4^2}. When $n\geq 3$, by Lemma \ref{strictly
cyclic g^n} there exists a strictly cyclic
$0$-FG$(3,(\emptyset,4),2^{n+1})$ of type $(2^n)^2$, denoted by
$(X,{\cal G},\emptyset,{\cal T})$, which is also a strictly
$(1,2^n)$-regular $0$-FG$(3,(\emptyset,4),2^{n+1})$ of type
$(2^n)^2$. Collect all $2$-subsets from distinct groups of $\cal G$
into a set $\cal B$. Then $(X,{\cal G},{\cal B},{\cal T})$ is a
strictly $(1,2^n)$-regular $1$-FG$(3,(2,4),2^{n+1})$ of type
$(2^n)^2$. Start from this $1$-FG and apply Construction
\ref{Weighting-strictly regular $s$-FG} with $h_1=2^n$ and $h_2=1$
to obtain a strictly $(2,2^n)$-regular
$0$-FG$(3,(\emptyset,4),2^{n+2})$ of type $(2^{n+1})^2$, where the
needed strictly $1$-cyclic $0$-FG$(3,(\emptyset,4),4)$ of type $2^2$
is from Theorem \ref{G-design}, and the needed $1$-cyclic
$H(4,2,4,3)$ is from Corollary \ref{cor for $h$-cyclic $H$ design}.
Now use induction on $n$. When $n=3$, there is an optimal $2$-D
$(2\times 2^3,4,2)$-OOC with $J^*(2\times 2^3)$ codewords from
Example \ref{example 8^2}. Assume that an optimal $2$-D $(2\times
2^n,4,2)$-OOC with $J^*(2\times 2^n)$ codewords exists for some
$n\geq 3$. Then start from a strictly $(2,2^n)$-regular
$0$-FG$(3,(\emptyset,4),2^{n+2})$ of type $(2^{n+1})^2$, and apply
Construction \ref{filling III} with an optimal $2$-D $(2\times
2^n,4,2)$-OOC with $J^*(2\times 2^n)$ codewords to obtain an $2$-D
$(2\times 2^{n+1},4,2)$-OOC, which contains
$(2^n-1)(2^{n+2}+1)/3=\lfloor\frac{2}{4}(\lfloor\frac{2^{n+2}-1}{3}\lfloor\frac{2^{n+2}-2}{2}
\rfloor\rfloor-1)\rfloor$ codewords. This number achieves the upper
bound in Theorem \ref{bound-total} (note that for any integer $n\geq
2$, $2^n\equiv 4,8$ $({\rm mod}$ $12)$). Thus an optimal $2$-D
$(2\times 2^{n+1},4,2)$-OOC with $J^*(2\times 2^{n+1})$ codewords
exists.  \qed

Combining the results of Proposition \ref{infinite family-u-v-1} and Lemmas \ref{infinite family-u-2-1}, \ref{infinite family-2^n}, we
have

\begin{Theorem}
\label{cor-2*2^n} There is an optimal $2$-D $(u\times 2^n,4,2)$-OOC
with $J^*(u\times 2^n)$ codewords for any $u\equiv 2,4$ $({\rm mod}$
$6)$ and any positive integer $n$.
\end{Theorem}

Proposition \ref{infinite family-u-v-1} can only deal with the case of
$u\equiv 2,4$ $({\rm mod}$ $6)$ and $v\not\equiv 2$ $({\rm mod}$
$4)$. When $v\equiv 2$ $({\rm mod}$ $4)$, we have the following
proposition.

\begin{Proposition}
\label{infinite family-u-v-b1} Let $u\equiv 2,4$ $({\rm mod}$ $6)$
and $v\equiv 2$ $({\rm mod}$ $4)$. Suppose that there is an optimal
$2$-D $(u/2\times 2v,4,2)$-OOC with $J^*(u/2\times 2v)$ codewords.
Then there is an optimal $2$-D $(u\times v,4,2)$-OOC with
$J^*(u\times v)$ codewords.
\end{Proposition}

\proof By Theorem \ref{relation2}, if there exists an optimal $2$-D
$(u/2\times 2v,4,2)$-OOC with $J^*(u/2\times 2v)$ codewords, then
there exits a $2$-D $(u\times v,4,2)$-OOC with $2J^*(u/2\times 2v)$
codewords. Note that by Theorem \ref{bound-total},
$$ J^*(u/2\times 2v)=\left\{
\begin{array}{lll}
    \lfloor\frac{u}{8}(\lfloor\frac{uv-1}{3}\lfloor\frac{uv-2}{2}
\rfloor\rfloor-1)\rfloor, & \hbox{if $v\equiv 2,10$ $({\rm mod}$
$12)$,}
\\\\
\lfloor\frac{u}{8}(\lfloor\frac{uv-1}{3}\lfloor\frac{uv-2}{2}
\rfloor\rfloor-2)\rfloor, & \hbox{if $v\equiv 6$ $({\rm mod}$ $12)$.} \\
\end{array}
\right.$$ It is readily checked that $2J^*(u/2\times 2v)=$
$$ \\\left\{
\begin{array}{lll}
    u(uv+1)(uv-4)/24=\lfloor\frac{u}{4}(\lfloor\frac{uv-1}{3}\lfloor\frac{uv-2}{2}
\rfloor\rfloor-1)\rfloor, & \hbox{if $v\equiv 2,10$ $({\rm mod}$
$12)$,}
\\\\
u(u^2v^2-3uv-12)/24=\lfloor\frac{u}{4}(\lfloor\frac{uv-1}{3}\lfloor\frac{uv-2}{2}
\rfloor\rfloor-2)\rfloor, & \hbox{if $v\equiv 6$ $({\rm mod}$ $12)$.} \\
\end{array}
\right.$$ This number achieves the upper bound in Theorem
\ref{bound-total}. Thus an optimal $2$-D $(u\times v,4,2)$-OOC with
$J^*(u\times v)$ codewords exists. \qed

The following proposition shows another approach to obtain some optimal
$2$-D $(u\times v,4,2)$-OOCs with $v\equiv 2$ $({\rm mod}$ $4)$.

\begin{Proposition}
\label{infinite family-u-v-5} If there is a perfect $2$-D $(2\times
v,4,2)$-OOC with $v\equiv 1,5$ $({\rm mod}$ $6)$, then there is an
optimal $2$-D $(u\times 2v,4,2)$-OOC  for any $u\equiv 8,16$ $({\rm
mod}$ $24)$ with $J^*(u\times 2v)$ codewords.
\end{Proposition}

\proof  Let $n=u/8$. Then $n\equiv 1,2$ $({\rm mod}$ $3)$. There is
a strictly cyclic $0$-FG$(3,(\emptyset,4)$, $16n)$ of type $(8n)^2$,
which exists by Lemma \ref{strictly cyclic g^n}. By Lemma
\ref{relation between regular and v-cyclic}, it implies a strictly
$2$-cyclic $0$-FG$(3,(\emptyset,4),16n)$ of type $(8n)^2$. By Lemma
\ref{strictly $(2,1)$-regular $1$-FG}, if there is a perfect $2$-D
$(2\times v,4,2)$-OOC with $v\equiv 1,5$ $({\rm mod}$ $6)$, then
there is a strictly $(2,1)$-regular $1$-FG$(3,(2,4),2v)$ of type
$2^v$. Start from this $1$-FG and apply Construction
\ref{Weighting-strictly regular $s$-FG} with $h_1=1$ and $h_2=2$ to
obtain a strictly $(8n,2)$-regular $0$-FG$(3,(\emptyset,4),16nv)$ of
type $(16n)^v$, where the needed $2$-cyclic $H(4,8n,4,3)$ is from
Corollary \ref{cor for $h$-cyclic $H$ design}. Applying Construction
\ref{filling III} with an optimal strictly $2$-cyclic $3$-$(8n\times
2,4,1)$ packing with $J^*(8n\times
2)=\lfloor\frac{8n}{4}(\lfloor\frac{16n-1}{3}\lfloor\frac{16n-2}{2}
\rfloor\rfloor-1)\rfloor$ base blocks, which exists by Lemma
\ref{infinite family-u-2-1}, we have a strictly $2v$-cyclic
$3$-$(8n\times 2v,4,1)$ packing, which contains $
4n(4nv-1)(16nv+1)/3=\lfloor\frac{8n}{4}(\lfloor\frac{16nv-1}{3}\lfloor\frac{16nv-2}{2}
\rfloor\rfloor-1)\rfloor$ base blocks. This number achieves the
upper bound in Theorem \ref{bound-total}. Thus an optimal $2$-D
$(8n\times 2v,4,2)$-OOC with $J^*(8n\times 2v)$ codewords exists. It
is an optimal $2$-D $(u\times 2v,4,2)$-OOC with $J^*(u\times 2v)$
codewords. \qed

\begin{Theorem}
\label{cor-u*2p}  Let $p\equiv 7 \ (mod \ 12)$ be a prime or $p\in
\{37$, $61$, $73$, $109$, $157$, $181$, $229$, $277\}$. There exists
an optimal $(u\times 2p,4,2)$-OOC with $J^*(u\times 2p)$ codewords
for any $u\equiv 8,16$ $({\rm mod}$ $24)$.
\end{Theorem}
\proof Start from a perfect $2$-D $(2\times p,4,2)$-OOC, which
exists by Theorem \ref{RoSQS-result}. Apply Proposition \ref{infinite
family-u-v-5} to complete the proof. \qed

\begin{Lemma}
\label{infinite family-u-v-2} Let $v\equiv 1$ $({\rm mod}$ $2)$ or
$v\equiv 0$ $({\rm mod}$ $12)$. Suppose that there is an optimal
$2$-D $(12\times v,4,2)$-OOC with $J^*(12\times v)$ codewords. Then
there is an optimal $2$-D $(u\times v,4,2)$-OOC with $J^*(u\times
v)$ codewords for any $u\equiv 0$ $({\rm mod}$ $12)$.
\end{Lemma}

\proof Let $n=u/12$. When $n=1$, the conclusion follows from the
assumption. When $n\geq 2$, by Theorem \ref{G-design}, there exists
a $0$-FG$(3,(\emptyset,4),6n)$ of type $6^n$ $(X,{\cal
G},\emptyset,{\cal T})$. Collect all $2$-subsets from distinct
groups of $\cal G$ into a set $\cal B$. Then $(X,{\cal G},{\cal
B},{\cal T})$ is a $1$-FG$(3,(2,4),6n)$ of type $6^n$. Apply
Construction \ref{Weighting-strictly $h$-cyclic $s$-FG} with $h_1=1$
and $h_2=v$ to obtain a strictly $v$-cyclic
$0$-FG$(3,(\emptyset,4),12nv)$ of type $(12v)^n$, where the needed
strictly $v$-cyclic
 $0$-FG$(3,(\emptyset,4),4v)$ of type $(2v)^2$ is from Corollary \ref{(2v)^2-1}, and the
 needed $v$-cyclic $H(4,2v,4,3)$ is from Corollary \ref{cor for $h$-cyclic $H$
 design}. Apply Construction \ref{filling II} with an optimal strictly
$v$-cyclic $3$-$(12\times v,4,1)$ packing with $J^*(12\times v)$
base blocks, which exists by assumption. Note that by Theorem
\ref{bound-total},
$$ J^*(12\times v)=\left\{
\begin{array}{lll}
    \lfloor\frac{12}{4}(\lfloor\frac{12v-1}{3}\lfloor\frac{12v-2}{2}
\rfloor\rfloor-1)\rfloor, & \hbox{if $v\equiv 1$ $({\rm mod}$
$2)$,} \\
    \\
    \lfloor\frac{12}{4}(\lfloor\frac{12v-1}{3}\lfloor\frac{12v-2}{2}
\rfloor\rfloor-2)\rfloor, & \hbox{if $v\equiv 0$ $({\rm mod}$
$12)$.}
\\
\end{array}
\right.$$ Then we have a strictly $v$-cyclic $3$-$(12n\times v,4,1)$
packing, which contains
$$ \left\{
\begin{array}{lll}
    3n(24n^2v^2-6nv-1)=\lfloor\frac{12n}{4}(\lfloor\frac{12nv-1}{3}\lfloor\frac{12nv-2}{2}
\rfloor\rfloor-1)\rfloor, & \hbox{if $v\equiv 1$ $({\rm mod}$
$2)$,} \\
    \\
    6n(12n^2v^2-3nv-1)=\lfloor\frac{12n}{4}(\lfloor\frac{12nv-1}{3}\lfloor\frac{12nv-2}{2}
\rfloor\rfloor-2)\rfloor, & \hbox{if $v\equiv 0$ $({\rm mod}$ $12)$,} \\
\end{array}
\right.$$ base blocks. This number achieves the upper bound in
Theorem \ref{bound-total}.
Thus an optimal $2$-D $(12n\times v,4,2)$-OOC with $J^*(12n\times
v)$ codewords exists. It is an optimal $2$-D $(u\times v,4,2)$-OOC
with $J^*(u\times v)$ codewords. \qed

The use of Lemma \ref{infinite family-u-v-2} depends on the
existence of optimal $2$-D $(12\times v,4,2)$-OOCs with
$J^*(12\times v)$ codewords. The following lemma shows an approach
to obtain some optimal $2$-D $(12\times v,4,2)$-OOCs from perfect
$2$-D $(2\times v,4,2)$-OOCs.

\begin{Lemma}
\label{infinite family-u-v-3} Suppose that there is a perfect $2$-D
$(2\times v,4,2)$-OOC with $v\equiv 1,5$ $({\rm mod}$ $6)$. Then
there is an optimal $2$-D $(12\times v,4,2)$-OOC with $J^*(12\times
v)$ codewords.
\end{Lemma}

\proof By Lemma \ref{strictly $(2,1)$-regular $1$-FG}, if there is a
perfect $2$-D $(2\times v,4,2)$-OOC with $v\equiv 1,5$ $({\rm mod}$
$6)$, then there is a strictly $(2,1)$-regular $1$-FG$(3,(2,4),2v)$
of type $2^v$. Start from this $1$-FG and apply Construction
\ref{Weighting-strictly regular $s$-FG} with $h_1=1$ and $h_2=1$ to
obtain a strictly $(12,1)$-regular $0$-FG$(3,(\emptyset,4),12v)$ of
type $12^v$, where the needed strictly $1$-cyclic
 $0$-FG$(3,(\emptyset,4),12)$ of type $6^2$ is from Theorem \ref{G-design}, and the
 needed $1$-cyclic $H(4,6,4,3)$ is from Corollary \ref{cor for $h$-cyclic $H$
 design}. Applying Construction \ref{filling III} with an optimal strictly
$1$-cyclic $3$-$(12\times 1,4,1)$ packing with $51$ base blocks from
Theorem \ref{ji-v=1}, we have a strictly $v$-cyclic $3$-$(12\times
v,4,1)$ packing, which contains
$72v^2-18v-3=\lfloor\frac{12}{4}(\lfloor\frac{12v-1}{3}\lfloor\frac{12v-2}{2}
\rfloor\rfloor-1)\rfloor$ base blocks. This number achieves the
upper bound in Theorem \ref{bound-total}. Thus an optimal $2$-D
$(12\times v,4,2)$-OOC with $J^*(12\times v)$ codewords exists. \qed

\begin{Lemma}
\label{infinite family-u-v-a1} Let $v\equiv 3$ $({\rm mod}$ $6)$ or
$v\equiv 0$ $({\rm mod}$ $12)$. Suppose that there is an optimal
$2$-D $(2\times v,4,2)$-OOC  with $J^*(2\times v)$ codewords. Then
there is an optimal $2$-D $(12\times v,4,2)$-OOC with $J^*(12\times
v)$ codewords.
\end{Lemma}

\proof By Proposition \ref{infinite family-u-v-1}, if there is an
optimal $2$-D $(2\times v,4,2)$-OOC with $J^*(2\times v)$ codewords
for $v\equiv 3$ $({\rm mod}$ $6)$ or $v\equiv 0$ $({\rm mod}$ $12)$,
then there is an optimal $2$-D $(4\times v,4,2)$-OOC with
$J^*(4\times v)$ codewords. Note that by Theorem \ref{bound-total},
$$ J^*(4\times v)=\left\{
\begin{array}{lll}
    \lfloor\frac{4}{4}(\lfloor\frac{4v-1}{3}\lfloor\frac{4v-2}{2}
\rfloor\rfloor-1)\rfloor, & \hbox{if $v\equiv 3$ $({\rm mod}$
$6)$,} \\
    \\
    \lfloor\frac{4}{4}(\lfloor\frac{4v-1}{3}\lfloor\frac{4v-2}{2}
\rfloor\rfloor-2)\rfloor, & \hbox{if $v\equiv 0$ $({\rm mod}$
$12)$.}
\\
\end{array}
\right.$$ By Lemma \ref{strictly cyclic g^n}, there is a strictly
cyclic $0$-FG$(3,(\emptyset,4),12v)$ of type $(4v)^3$. By Lemma
\ref{relation between regular and v-cyclic}, it implies a strictly
$v$-cyclic $0$-FG$(3,(\emptyset,4),12v)$ of type $(4v)^3$. Start
from this strictly $v$-cyclic $0$-FG and apply Construction
\ref{filling II} with an optimal strictly $v$-cyclic $3$-$(4\times
v,4,1)$ packing with $J^*(4\times v)$ base blocks, which is
equivalent to an optimal $2$-D $(4\times v,4,2)$-OOC with
$J^*(4\times v)$ codewords, to obtain a strictly $v$-cyclic
$3$-$(12\times v,4,1)$ packing, which contains
$$ \left\{
\begin{array}{lll}
    72v^2-18v-3=\lfloor\frac{12}{4}(\lfloor\frac{12v-1}{3}\lfloor\frac{12v-2}{2}
\rfloor\rfloor-1)\rfloor, & \hbox{if $v\equiv 3$ $({\rm mod}$
$6)$,} \\
    \\
    72v^2-18v-6=\lfloor\frac{12}{4}(\lfloor\frac{12v-1}{3}\lfloor\frac{12v-2}{2}
\rfloor\rfloor-2)\rfloor, & \hbox{if $v\equiv 0$ $({\rm mod}$ $12)$,} \\
\end{array}
\right.$$ base blocks. This number achieves the upper bound in
Theorem \ref{bound-total}. Thus an optimal $2$-D $(12\times
v,4,2)$-OOC with $J^*(12\times v)$ codewords exists. \qed

Combining the results of Lemmas \ref{infinite family-u-v-2},
\ref{infinite family-u-v-3} and \ref{infinite family-u-v-a1}, we
have the following proposition.

\begin{Proposition}
\label{infinite family-u-v-a2} Let $v\equiv 1$ $({\rm mod}$ $2)$ or
$v\equiv 0$ $({\rm mod}$ $12)$. Suppose that there is an optimal
$2$-D $(2\times v,4,2)$-OOC  with $J^*(2\times v)$ codewords. Then
there is an optimal $2$-D $(u\times v,4,2)$-OOC with $J^*(u\times
v)$ codewords for any $u\equiv 0$ $({\rm mod}$ $12)$.
\end{Proposition}

\begin{Lemma}
\label{example-3-cyclic-6^3} There exists a strictly $3$-cyclic
$0$-FG$(3,(\emptyset,4),18)$ of type $6^3$.
\end{Lemma}

\proof We here give a construction of a $0$-FG$(3,(\emptyset,4),18)$ of type $6^3$ on
$I_{18}$ with the group set $\{\{0,1,2,3,4,5\}$ $+i:i\in
\{0,6,12\}\}$. Let $\alpha=(0\ 1\ 2)(3\ 4\ 5)(6\ 7\ 8)(9\ 10\
11)(12\ 13\ 14)(15\ $ $16\ 17)$ be a permutation on $I_{18}$ and $G$ be
the group generated by $\alpha$. Only base blocks are listed below.
All other blocks are obtained by developing these base blocks under
the action of $G$. Obviously this design is isomorphic to a strictly
$3$-cyclic $0$-FG$(3,(\emptyset,4),18)$ of type $6^3$.

\begin{center}
{\footnotesize \tabcolsep 0.03in \scriptsize
\begin{tabular}{lllllll}
$\{0,1,6,7\}$& $\{0,1,8,11\}$& $\{0,1,9,10\}$& $\{0,1,12,13\}$&
$\{0,1,14,15\}$& $\{0,1,16,17\}$& $\{0,3,6,9\}$\\ $\{0,3,7,8\}$&
$\{0,3,10,11\}$& $\{0,3,12,15\}$& $\{0,3,13,16\}$&
$\{0,3,14,17\}$& $\{0,4,6,12\}$& $\{0,4,7,14\}$\\
$\{0,4,8,16\}$& $\{0,4,9,13\}$& $\{0,4,10,15\}$& $\{0,4,11,17\}$&
$\{0,5,6,16\}$& $\{0,5,7,12\}$&
$\{0,5,8,17\}$\\
$\{0,5,9,14\}$& $\{0,5,10,13\}$& $\{0,5,11,15\}$& $\{0,6,10,17\}$&
$\{0,6,11,14\}$& $\{0,6,13,15\}$& $\{0,7,9,16\}$\\ $\{0,7,11,13\}$&
$\{0,7,15,17\}$& $\{0,8,9,15\}$& $\{0,8,10,12\}$& $\{0,8,13,14\}$&
$\{0,9,12,17\}$& $\{0,10,14,16\}$\\ $\{0,11,12,16\}$& $\{3,4,6,7\}$&
$\{3,4,8,11\}$& $\{3,4,9,10\}$& $\{3,4,12,17\}$& $\{3,4,13,14\}$&
$\{3,4,15,16\}$\\ $\{3,6,10,14\}$& $\{3,6,11,17\}$& $\{3,6,12,16\}$&
$\{3,7,9,13\}$&
$\{3,7,11,16\}$& $\{3,7,12,14\}$& $\{3,8,9,12\}$\\
$\{3,8,10,15\}$& $\{3,8,16,17\}$& $\{3,9,14,15\}$& $\{3,10,13,17\}$&
$\{3,11,13,15\}$& $\{6,7,12,17\}$& $\{6,7,13,16\}$\\
$\{6,7,14,15\}$& $\{6,9,12,13\}$& $\{6,9,14,17\}$& $\{6,9,15,16\}$&
$\{9,10,12,14\}$& $\{9,10,13,16\}$& $\{9,10,15,17\}$
\end{tabular}}
\end{center}

\begin{Lemma}
\label{infinite family-u-3-2} There exists an optimal $2$-D
$(u\times 3,4,2)$-OOC with $J^*(u\times 3)$ codewords for any
$u\equiv 6$ $({\rm mod}$ $12)$.
\end{Lemma}

\proof Let $n=(u+2)/2$. Then $n\equiv 4$ $({\rm mod}$ $6)$. There is
an SQS$(n)$ \cite{hanani60}. Delete one point to obtain a
$1$-FG$(3,(3,4),n-1)$ of type $1^{n-1}$. Start from this $1$-FG and
apply Construction \ref{Weighting-strictly $h$-cyclic $s$-FG} with
$h_1=1$ and $h_2=3$ to obtain a strictly $3$-cyclic
$0$-FG$(3,(\emptyset,4),6(n-1))$ of type $6^{n-1}$, where the needed
strictly $3$-cyclic $0$-FG$(3,(\emptyset,4),18)$ of type $6^3$ is
from Lemma \ref{example-3-cyclic-6^3}, and the needed $3$-cyclic
$H(4,6,4,3)$ is from Corollary \ref{cor for $h$-cyclic $H$ design}.
Applying Construction \ref{filling II} with an optimal strictly
$3$-cyclic $3$-$(2\times 3,4,1)$ packing with $J^*(2\times 3)=1$
base block from Example \ref{example-6}, we have a strictly
$3$-cyclic $3$-$(2(n-1)\times 3,4,1)$ packing, which contains
$(n-1)(6n^2-15n+8)/2$ base blocks. This number achieves the upper
bound in Theorem \ref{bound-total}. Thus an optimal $2$-D
$(2(n-1)\times 3,4,2)$-OOC with $J^*(2(n-1)\times 3)$ codewords
exists. It is an optimal $2$-D $(u\times 3,4,2)$-OOC with
$J^*(u\times 3)$ codewords.  \qed

\begin{Theorem}
\label{the-u*3} There exists an optimal $2$-D $(u\times 3,4,2)$-OOC
with $J^*(u\times 3)$ codewords for any $u\equiv 0$ $({\rm mod}$
$2)$.
\end{Theorem}
\proof When $u\equiv 2,4$ $({\rm mod}$ $6)$, apply Proposition
\ref{infinite family-u-v-1} with an optimal $2$-D $(2\times
3,4,2)$-OOC with $J^*(2\times 3)=1$ base block from Example
\ref{example-6} to obtain an optimal $2$-D $(u\times 3,4,2)$-OOC
with $J^*(u\times 3)$ codewords. When $u\equiv 0$ $({\rm mod}$
$12)$, apply Proposition \ref{infinite family-u-v-a2} with an optimal
$2$-D $(2\times 3,4,2)$-OOC with $J^*(2\times 3)=1$ base block from
Example \ref{example-6} to obtain an optimal $2$-D $(u\times
3,4,2)$-OOC with $J^*(u\times 3)$ codewords. When $u\equiv 6$ $({\rm
mod}$ $12)$, the conclusion follows from Lemma \ref{infinite
family-u-3-2}.  \qed

\begin{Proposition}
\label{infinite family-Ro} If there is an $RoSQS(v+1)$, then there
is an optimal $2$-D $(u\times v,4,2)$-OOC with $J^*(u\times v)$
codewords for any $u\equiv 0$ $({\rm mod}$ $6)$.
\end{Proposition}

\proof By Lemma \ref{equiv-RoSQS}, when $v\equiv 1$ $({\rm mod}$
$6)$, an RoSQS$(v+1)$ is equivalent to a strictly cyclic
$1$-$FG(3,(3,4),v)$ of type $1^v$, which is also a strictly
$(1,1)$-regular $1$-FG$(3,(3,4),v)$ of type $1^v$. Start from this
$1$-FG and apply Construction \ref{Weighting-strictly regular
$s$-FG} with $h_1=1$ and $h_2=1$ to obtain a strictly
$(u,1)$-regular $0$-FG$(3,(\emptyset,4),uv)$ of type $u^v$ for any
$u\equiv 0$ $({\rm mod}$ $6)$, where the needed strictly $1$-cyclic
$0$-FG$(3,(\emptyset,4),3u)$ of type $u^3$ exists from Theorem
\ref{G-design}, and the needed $1$-cyclic $H(4,u,4,3)$ is from
Corollary \ref{cor for $h$-cyclic $H$
 design}. Applying Construction \ref{filling III} with an optimal strictly
$1$-cyclic $3$-$(u\times 1,4,1)$ packing with $J^*(u\times 1)$ base
blocks from Theorem \ref{ji-v=1}, we have a strictly $v$-cyclic
$3$-$(u\times v,4,1)$ packing, which contains
$u(u^2v^2-3uv-6)/24=\lfloor\frac{u}{4}(\lfloor\frac{uv-1}{3}\lfloor\frac{uv-2}{2}
\rfloor\rfloor-1)\rfloor$ base blocks. This number achieves the
upper bound in Theorem \ref{bound-total}. Thus an optimal $2$-D
$(u\times v,4,2)$-OOC with $J^*(u\times v)$ codewords exists.

By Lemma \ref{equiv-RoSQS}, when $v\equiv 3$ $({\rm mod}$ $6)$, an
RoSQS$(v+1)$ is equivalent to a strictly cyclic $1$-$FG(3,(3,4),v)$
of type $3^{v/3}$, which is also a strictly $(1,3)$-regular
$1$-FG$(3,(3,4),v)$ of type $3^{v/3}$. Start from this $1$-FG and
apply Construction \ref{Weighting-strictly regular $s$-FG} with
$h_1=3$ and $h_2=1$ to obtain a strictly $(u,3)$-regular
$0$-FG$(3,(\emptyset,4),uv)$ of type $(3u)^{v/3}$ for any $u\equiv
0$ $({\rm mod}$ $6)$. Applying Construction \ref{filling III} with
an optimal strictly $3$-cyclic $3$-$(u\times 3,4,1)$ packing with
$J^*(u\times 3)$ base blocks from Theorem \ref{the-u*3}, we have a
strictly $v$-cyclic $3$-$(u\times v,4,1)$ packing, which contains
$u(u^2v^2-3uv-6)/24=\lfloor\frac{u}{4}(\lfloor\frac{uv-1}{3}\lfloor\frac{uv-2}{2}
\rfloor\rfloor-1)\rfloor$ base blocks. This number achieves the
upper bound in Theorem \ref{bound-total}. Thus an optimal $2$-D
$(u\times v,4,2)$-OOC with $J^*(u\times v)$ codewords exists. \qed

Combining Theorem \ref{RoSQS-update} and Proposition \ref{infinite
family-Ro}, many infinite families of optimal $2$-D $(u\times
v,4,2)$-OOCs with $J^*(u\times v)$ codewords will be obtained. As an
example, we have

\begin{Theorem}
\label{cor-u*p}  Let $p\equiv 7 \ (mod \ 12)$ be a prime or $p\in
\{37$, $61$, $73$, $109$, $157$, $181$, $229$, $277\}$. There exist
a perfect $2$-D $(u\times p,4,2)$-OOC for any $u\equiv 2,4$ $({\rm
mod}$ $6)$, and an optimal $2$-D $(u\times p,4,2)$-OOC with
$J^*(u\times p)$ codewords for any $u\equiv 0$ $({\rm mod}$ $6)$.
\end{Theorem}
\proof When $u\equiv 2,4$ $({\rm mod}$ $6)$, start from a perfect
$2$-D $(2\times p,4,2)$-OOC, which exists by Theorem
\ref{RoSQS-result}, and apply Proposition \ref{infinite family-u-v-1} to
obtain a perfect $2$-D $(u\times p,4,2)$-OOC. When $u\equiv 0$
$({\rm mod}$ $6)$, start from an RoSQS$(p+1)$, which exists by
Theorem \ref{RoSQS-update}, and apply Proposition \ref{infinite
family-Ro} to obtain an optimal $2$-D $(u\times p,4,2)$-OOC with
$J^*(u\times p)$ codewords.  \qed\vspace{-4mm}

\begin{Lemma}
\label{infinite family-12u-2v} If there is a perfect $2$-D $(2\times
v,4,2)$-OOC with $v\equiv 1,5$ $({\rm mod}$ $6)$, then there is an
optimal $2$-D $(12\times 2v,4,2)$-OOC with $J^*(12\times 2v)$
codewords.
\end{Lemma}

\proof By Lemma \ref{strictly $(2,1)$-regular $1$-FG}, if there is a
perfect $2$-D $(2\times v,4,2)$-OOC with $v\equiv 1,5$ $({\rm mod}$
$6)$, then there is a strictly $(2,1)$-regular $1$-FG$(3,(2,4),2v)$
of type $2^v$. Start from this $1$-FG and apply Construction
\ref{Weighting-strictly regular $s$-FG} with $h_1=1$ and $h_2=2$ to
obtain a strictly $(12,2)$-regular $0$-FG$(3,(\emptyset,4),24v)$ of
type $24^v$, where the needed strictly $2$-cyclic
 $0$-FG$(3,(\emptyset,4),24)$ of type $12^2$ is from Lemma \ref{example-2-cyclic-12^2}, and the
 needed $2$-cyclic $H(4,12,4,3)$ is from Corollary \ref{cor for $h$-cyclic $H$
 design}. Applying Construction \ref{filling III} with an optimal strictly
$2$-cyclic $3$-$(12\times 2,4,1)$ packing with $J^*(12\times 2)=248$
base blocks from Lemma \ref{12*2}, we have a strictly $2v$-cyclic
$3$-$(12\times 2v,4,1)$ packing, which contains
$4(72v^2-9v-1)=\lfloor\frac{12}{4}(\lfloor\frac{24v-1}{3}\lfloor\frac{24v-2}{2}
\rfloor\rfloor-1)\rfloor-1$ base blocks. This number achieves the
upper bound in Theorem \ref{bound-total}. Thus an optimal $2$-D
$(12\times 2v,4,2)$-OOC with $J^*(12\times 2v)$ codewords exists.
\qed

\begin{Theorem}
\label{cor-12*2p}  Let $p\equiv 7 \ (mod \ 12)$ be a prime or $p\in
\{37$, $61$, $73$, $109$, $157$, $181$, $229$, $277\}$. There exists
an optimal $(12\times 2p,4,2)$-OOC with $J^*(12\times 2p)$
codewords.
\end{Theorem}
\proof Start from a perfect $2$-D $(2\times p,4,2)$-OOC, which
exists by Theorem \ref{RoSQS-result}. Apply Lemma \ref{infinite
family-12u-2v} to complete the proof. \qed

 {\footnotesize
\tabcolsep 0.07in
\begin{center}
{\bf Table II\\ Small orders of optimal $2$-D $(u\times v,4,2)$-OOCs
with $\Phi(u\times v,4,2)=J^*(u\times v)$\\ codewords for $6\leq
uv\leq 34$} \vspace{6pt}

\begin{tabular}{|c|ccc|ccc|}\hline
$uv$ & $u\times v$ & $\Phi(u\times v,4,2)$ & Source & $u\times v$ &
$\Phi(u\times v,4,2)$ & Source
\\\hline
$6$ & $2\times3$ & 1 & Example \ref{example-6} & $3\times2$ & $1$ &
Example \ref{example-6}
\\ $8$ & $2\times4$ & $3$ & Example \ref{example-4-cyclic-4^2} &
$4\times2$ & $6$ & Example \ref{example-2-cyclic-4^2}
\\ $9$ & $3\times3$ & $6$ & Lemma \ref{small order} & & &
\\ $10$ & $2\times5$ & $6$ & Theorem \ref{2D k=4 from 1D}$(2)$ & $5\times2$ & $15$ &  Theorem \ref{2D k=4 from 1D}$(2)$
\\ $12$ & $2\times6$ & $8$ & Lemma \ref{small order} & $3\times4$ & $12$ & Lemma \ref{small order}
\\ & $4\times3$ & $17$ & Example \ref{example-3-cyclic-6^2} &$6\times2$ & $25$ & Lemma \ref{small order}
\\ $14$ & $2\times7$ & $13$ & Theorem \ref{RoSQS-result} & $7\times2$ & $44$ & Lemma \ref{small order}
\\ $15$ & $3\times5$ & $21$ & Theorem \ref{2D k=4 from 1D}$(2)$ & $5\times3$ & $35$ & Theorem \ref{2D k=4 from 1D}$(2)$
\\ $16$ & $2\times8$ & $17$ & Example \ref{example 8^2} & $4\times4$ &$34$ & Theorem \ref{cor-2*2^n}
\\ & $8\times2$ & $68$ & Lemma \ref{infinite family-u-2-1}& & &
\\ $18$ & $2\times9$ & $22$ & Corollary \ref{18,42,90} & $3\times6$ & $33$ & Corollary \ref{18,42,90}
\\ & $6\times3$ & $66$ & Corollary \ref{18,42,90} & $9\times2$ & $99$ & Corollary \ref{18,42,90}
\\ $20$ & $2\times10$ & $28$ & Theorem \ref{2D k=4 from 1D}$(4)$ & $4\times5$ & $57$ & Proposition \ref{infinite family-u-v-1}
\\ & $5\times4$ & $70$ & Lemma \ref{new small order-u*v derived} & $10\times2$ & $140$ & Lemma \ref{infinite family-u-2-1}
\\ $21$ & $3\times7$ & $45$ & Theorem \ref{2D k=4 from 1D}$(2)$ & $7\times3$ & $105$ & Theorem \ref{2D k=4 from 1D}$(2)$
\\ $22$ & $2\times11$ & $35$ & Lemma \ref{small order} & $11\times2$ & & ??
\\ $24$ & $2\times12$ & $41$ & Lemma \ref{2*12} & $3\times8$ & & ??
\\ & $4\times6$ & $82$ & Proposition \ref{infinite family-u-v-b1}&
$6\times4$ & & ??
\\& $8\times3$ & $166$ & Theorem \ref{the-u*3} & $12\times2$ & $248$ & Lemma \ref{12*2}
\\ $25$ & $5\times5$ & $110$ &Theorem \ref{2D k=4 from 1D}$(2)$& & &
\\ $26$ & $2\times13$ & $50$ & Theorem \ref{2D k=4 from 1D}$(2)$& $13\times2$ & $325$ & Theorem \ref{2D k=4 from 1D}$(2)$
\\ $27$ & $3\times9$ & $78$ & Theorem \ref{2D k=4 from 1D}$(2)$& $9\times3$ & $234$ & Theorem \ref{2D k=4 from 1D}$(2)$
\\ $28$ & $2\times14$ & $58$ & Theorem \ref{2D k=4 from 1D}$(4)$& $4\times7$ & $117$ & Proposition \ref{infinite family-u-v-1}
\\ &$7\times4$ & $203$ & Lemma \ref{new small order-u*v derived}& $14\times2$ & $406$ & Lemma \ref{infinite family-u-2-1}
\\ $30$ & $2\times15$ & $67$ & Lemma \ref{2*15}& $3\times10$ & $100$ & Lemma \ref{3*10}
\\& $5\times6$ & & ?? & $6\times5$ & $201$ & Lemma \ref{new small order-u*v
derived}\\& $10\times3$ & $335$ & Theorem \ref{the-u*3} &
$15\times2$ & & ??
\\ $32$ & $2\times16$ & $77$ & Theorem \ref{cor-2*2^n}& $4\times8$ & $154$ & Theorem \ref{cor-2*2^n}
\\& $8\times4$ & $308$ & Theorem \ref{cor-2*2^n}& $16\times2$ & $616$ & Lemma
\ref{infinite family-u-2-1}
\\ $33$ & $3\times11$ & $120$ & Theorem \ref{2D k=4 from 1D}$(2)$&
$11\times3$ & $440$ & Theorem \ref{2D k=4 from 1D}$(2)$
\\ $34$ & $2\times17$ & $88$ & Theorem \ref{2D k=4 from 1D}$(2)$&
$17\times2$ & $748$ & Theorem \ref{2D k=4 from 1D}$(2)$\\\hline
\end{tabular} \end{center} }

Finally we summarize the existence of small orders of optimal $2$-D
$(u\times v,4,2)$-OOCs with $J^*(u\times v)$ codewords as follows.

\begin{Theorem}
\label{summarize small orders} There exists an optimal $2$-D
$(u\times v,4,2)$-OOC with $J^*(u\times v)$ codewords for each
$(u,v)$ satisfying $6\leq uv\leq 34$ and $(u,v)\not\in
\{(1,9),(1,12),(1,13),(11,2),(23,1)$, $(3,8),(6,4),(5,6),(15,2)\}$.
When $(u,v)\in \{(1,9),(1,12),(1,13)\}$, there exists an optimal
$2$-D $(u\times v,4,2)$-OOC with $J(1\times uv,4,2)-1$ codewords.
\end{Theorem}

\proof By Theorem \ref{ji-v=1}, there exists an optimal $2$-D
$(u\times 1,4,2)$-OOC with $J^*(u\times 1)$ codewords for each
$6\leq u\leq 34$ and $u\neq 23$. By Lemma \ref{1D k=4}$(3)$, there
exists an optimal $2$-D $(1\times v,4,2)$-OOC with $J^*(1\times v)$
codewords for each $7\leq v\leq 34$ and $v\not\in\{9,12,13\}$. Note
that $J^*(1\times v)=J(1\times v,4,2)$ when $v\not\equiv 0$ $({\rm
mod}$ $24)$, and $J^*(1\times v)=J(1\times v,4,2)-1$ when $v\equiv
0$ $({\rm mod}$ $24)$. When $v\in\{9,12,13\}$, by Lemma \ref{1D
k=4}$(4)$, there exists an optimal $2$-D $(1\times v,4,2)$-OOC with
$J(1\times v,4,2)-1$ codewords. An optimal $2$-D $(1\times
6,4,2)$-OOC is trivial without codewords. For all other cases of
$(u,v)$ such that there is an optimal $2$-D $(u\times v,4,2)$-OOC
with $J^*(u\times v)$ codewords, we show the sources in Table II,
where the question marks "??" indicates the orders for each of which
the existence of an optimal $2$-D $(u\times v,4,2)$-OOC with
$J^*(u\times v)$ codewords is still open. \qed\vspace{-2mm}

\section{Conclusion}

In this paper, we gave some combinatorial constructions for optimal
$2$-D $(u\times v,k,2)$-OOCs. As applications, many infinite
families of optimal $2$-D $(u\times v,4,2)$-OOCs are obtained. We summarize all infinite families obtained in this paper in Table III. Although we can not complete the existence of optimal $2$-D $(u\times v,4,2)$-OOCs, we hope to present some possible approaches to reduce the existence problem. We summarize these approaches in Table IV.

 {\footnotesize
\tabcolsep 0.07in
\begin{center}
{\bf Table III\\ New infinite families of optimal $2$-D $(u\times v,4,2)$-OOCs with $J^*(u\times v)$ codewords \\ in this paper} \vspace{6pt}

\begin{tabular}{|c|c|c|}\hline
Parameters & Conditions & Source\\\hline

& $n_1n_2=uv$, & \\
&$u\in \{4^n-1: \ n\ge 1\}\cup \{1,27,33,39,51,87,123,183\}$ and &\\
$(n_1\times n_2,4,2)$ &$v\in S=\{p\equiv 7 \ (mod \ 12):
p$ is a prime$\}$ $\cup$& Theorem \ref{2D k=4 from 1D}$(1)$\\
&$ \{2^n-1: {\rm odd \ integer} \ n\geq 1\}$ $\cup\{25,37,61,73,109,157,181,229,277\}$,&\\
& or $v$ is a product of integers from $S$ &\\\hline

& $n=p_1^{r_1}p_2^{r_2}\cdots p_s^{r_s}$ & \\
$(2\times 2n,4,2)$ & $p_i=13$ or $p_i\equiv 5\ (mod\ 12)$ is a
prime and $p_i<1500000$, &Theorem \ref{2D k=4 from 1D}$(3)$\\
&$r_i\geq 1$ for
$1\leq i\leq s$&\\\hline

$(u\times 2^n,4,2)$ & $u\equiv 2,4$ $({\rm mod}$
$6)$ and $n\geq 1$& Theorem \ref{cor-2*2^n}\\\hline

$(u\times 2p,4,2)$ & $u\equiv 8,16$ $({\rm mod}$ $24)$ or $u=12$ & Theorem \ref{cor-u*2p}\\
& $p\equiv 7 \ (mod \ 12)$ a prime or $p\in
\{37,61,73,109,157,181,229,277\}$ & Theorem \ref{cor-12*2p}\\\hline

$(u\times p,4,2)$ & $u\equiv 0$ $({\rm mod}$ $2)$ & Theorem \ref{the-u*3}\\
& $p\equiv 7 \ (mod \ 12)$ a prime or $p\in
\{3,37,61,73,109,157,181,229,277\}$ & Theorem \ref{cor-u*p}\\\hline
\end{tabular} \end{center} }

\vspace{0.5cm}
 {\footnotesize
\tabcolsep 0.07in
\begin{center}
{\bf Table IV\\ Possible approaches to construct optimal $2$-D $(u\times v,4,2)$-OOCs} \vspace{2pt}

\begin{tabular}{|c|c|c|c|}\hline
Input & $\Rightarrow$ & Output & Source\\\hline

optimal $1$-D $(n,4,2)$-OOC & & optimal $2$-D
$(u\times v,4,2)$-OOC & \\
with $J(1\times n,4,2)$ codewords, && with $J(u\times v,4,2)$ codewords, & Corollary \ref{2D from 1D k=4}$(1)$ \\
$n\equiv 1,3\
({\rm mod }\ 6)$ or $n\equiv 2,10\ ({\rm mod }\ 24)$ & & for any
integer factorization $n=uv$ & \\\hline

optimal $2$-D$(2\times v,4,2)$-OOC && optimal $2$-D $(u\times v,4,2)$-OOC &  \\
with $J^*(2\times v)$ codewords, && with $J^*(u\times v)$ codewords, & Proposition \ref{infinite family-u-v-1} \\
$v\equiv 1$ $({\rm mod}$ $2)$ or $v\equiv 0$ $({\rm mod}$ $4)$ && $u\equiv 2,4$ $({\rm mod}$ $6)$ & \\\hline

optimal $2$-D $(u/2\times 2v,4,2)$-OOC && optimal
$2$-D $(u\times v,4,2)$-OOC & \\
with $J^*(u/2\times 2v)$ codewords, && with
$J^*(u\times v)$ codewords & Proposition \ref{infinite family-u-v-b1} \\
$u\equiv 2,4$ $({\rm mod}$ $6)$
and $v\equiv 2$ $({\rm mod}$ $4)$ &&&\\\hline

perfect $2$-D $(2\times v,4,2)$-OOC && optimal $2$-D $(u\times 2v,4,2)$-OOC &\\
$v\equiv 1,5$ $({\rm mod}$ $6)$ && with $J^*(u\times 2v)$ codewords, & Proposition \ref{infinite family-u-v-5}\\
&& $u\equiv 8,16$ $({\rm mod}$ $24)$ & \\\hline

optimal $2$-D $(2\times v,4,2)$-OOC && optimal $2$-D $(u\times v,4,2)$-OOC &\\
with $J^*(2\times v)$ codewords, && with $J^*(u\times v)$ codewords, & Proposition \ref{infinite family-u-v-a2}\\
$v\equiv 1$ $({\rm mod}$ $2)$ or $v\equiv 0$ $({\rm mod}$ $12)$ && $u\equiv 0$ $({\rm mod}$ $12)$ &\\\hline

$RoSQS(v+1)$ && optimal $2$-D $(u\times v,4,2)$-OOC &\\
$v\equiv 1,3$ $({\rm mod}$ $6)$&& with $J^*(u\times v)$ codewords, & Proposition \ref{infinite family-Ro} \\
&& $u\equiv 0$ $({\rm mod}$ $6)$ &\\

\hline
\end{tabular} \end{center} }

By Theorem \ref{bound-total}, we see that in many cases the Johnson bound can not be achieved. A natural question is whether the bounds established in Theorem
\ref{bound-total} is good enough to make each optimal $2$-D
$(u\times v,4,2)$-OOC achieve it. Although many infinite families
are given to achieve the upper bound in Theorem \ref{bound-total},
we still tend to think it not true. For example we conjecture that
when $u\equiv 0$ $({\rm mod}$ $6)$ and $v\equiv 2,4$ $({\rm mod}$
$6)$, the upper bound is
$\lfloor\frac{u}{4}(\lfloor\frac{uv-1}{3}\lfloor\frac{uv-2}{2}
\rfloor\rfloor-1)\rfloor-\lfloor\frac{u}{12}\rfloor$. If the
conjecture is correct, the condition in Lemma \ref{infinite
family-u-v-2} can be relaxed to $v\not\equiv 2$ $({\rm mod}$ $4)$, which implies that the condition in Proposition \ref{infinite family-u-v-a2} can also be relaxed to $v\not\equiv 2$ $({\rm mod}$ $4)$.

Another question is to find more constructions for optimal $2$-D
$(2\times v,4,2)$-OOCs, which are very useful by Propositions
\ref{infinite family-u-v-1}, \ref{infinite family-u-v-5} and \ref{infinite family-u-v-a2}. In $1991$ Phelps \cite{Phelps} constructed a class of $2$-chromatic
SQS$(22)$ using cyclic large sets of $2$-$(11,3,1)$ packings. It
seems that Phelps's method can be generalized to construct some
strictly $v$-cyclic SQS$(2\times v)$s for $v\equiv 1,5$ $({\rm mod}$
$6)$, which are also perfect $2$-D $(2\times v,4,2)$-OOCs. The
interested reader may refer to the paper \cite{Phelps}.

\begin{center}\noindent{\bf Appendix I}\end{center}
\textbf{Proof of Construction \ref{RoSQS-recur}:} For checking the correctness of the algorithm shown in Figure $1$, first count the number of base blocks in $\cal A$. It is clear that $|{\cal A}_1\cup{\cal A}'_1|=(p-1)(p-3)/12$, $|{\cal A}_2\cup{\cal
A}'_2|=(p-1)/3$, $|{\cal A}_3|\leq 3\times
(p-1)/6\times (p-1)/2=(p-1)^2/4$. Thus $|{\cal A}|\leq
(p-1)(2p-1)/6$.

Since the number $(p-1)(2p-1)/6$ is the right number of base blocks
in a strictly $p$-cyclic $SQS(2\times p)$, in the following it
suffices to show that each triple of $I_2\times Z_p$ appears in at
least one block of the resulting design. $(1)$ Each triple of
$\{0\}\times Z_p$ appears in one block of ${\cal A}_1\cup {\cal
A}_2$ and their cyclic shifts. $(2)$ Each triple of $\{1\}\times
Z_p$ appears in one block of ${\cal A}'_1\cup {\cal A}'_2$ and their
cyclic shifts. $(3)$ Each triple of the form $\{x_0,y_0,z_1\}$
appears in one block of ${\cal A}_2\cup {\cal A}_3$ and their cyclic
shifts. $(4)$ Each triple of the form $\{x_1,y_1,z_0\}$ appears in
one block of ${\cal A}'_2\cup {\cal A}_3$ and their cyclic shifts. \qed

\begin{center}\noindent{\bf Appendix II}\end{center}
\textbf{Lemma \ref{bound-u-v-12-2}}\hspace{0.2cm} {\it Let $u\equiv 0$ $({\rm mod}$ $12)$ and $v\equiv 2,4$ $({\rm mod }$ $6)$. Then $\Phi(u\times v,4,2)\leq
\lfloor\frac{u}{4}(\lfloor\frac{uv-1}{3}\lfloor\frac{uv-2}{2}
\rfloor\rfloor-1)\rfloor-1$.}

\vspace{0.4cm}
\proof First we shall show that $\Phi(u\times 2,4,2)\leq
\lfloor\frac{u}{4}(\lfloor\frac{2u-1}{3}\lfloor\frac{2u-2}{2}
\rfloor\rfloor-1)\rfloor-1$. By Lemma \ref{bound-6 mod 0},
$\Phi(u\times 2,4,2)\leq
\lfloor\frac{u}{4}(\lfloor\frac{2u-1}{3}\lfloor\frac{2u-2}{2}
\rfloor\rfloor-1)\rfloor$. Suppose that $\Phi(u\times 2,4,2)=
\lfloor\frac{u}{4}(\lfloor\frac{2u-1}{3}\lfloor\frac{2u-2}{2}
\rfloor\rfloor-1)\rfloor$. Then there were a strictly $2$-cyclic
$3$-$(u\times 2,4,1)$-packing with
$\lfloor\frac{u}{4}(\lfloor\frac{2u-1}{3}\lfloor\frac{2u-2}{2}
\rfloor\rfloor-1)\rfloor$ base blocks. Let $\cal L$ be the leave of
the strictly $2$-cyclic $3$-$(u\times 2,4,1)$-packing. Count the
number of $3$-subsets in the leave $\cal L$. It is
${{2u}\choose{3}}-
\lfloor\frac{u}{4}(\lfloor\frac{2u-1}{3}\lfloor\frac{2u-2}{2}
\rfloor\rfloor-1)\rfloor\cdot 2\cdot4=8u/3$.

For each $a\in I_u$ and each $i\in Z_2$, consider the number $n$ of
$3$-subsets containing the point $(a,i)$ in the leave $\cal L$.
Delete one point from a strictly $2$-cyclic $3$-$(u\times
2,4,1)$-packing to obtain a $2$-$(2u-1,3,1)$-packing, which contains
at most $\lfloor(2u-1)(2u-2)/6\rfloor-1$ blocks when $2u\equiv 0$
$({\rm mod}$ $6)$ \cite{hanani75}. Since each $3$-subset of
$I_u\times Z_2$ occurs in at most one block, we have $n\geq
{{2u-1}\choose{2}}-3(\lfloor(2u-1)(2u-2)/6\rfloor-1)=4$, which
implies that $|{\cal L}|\geq 4\cdot2u/3$. Due to $|{\cal L}|=8u/3$,
$n$ must be equal to $4$. Note that the above conclusion holds for
each $a\in I_u$ and each $i\in Z_2$.

For each $a\in I_u$, consider the number $m$ of the base blocks
containing the two points $(a,0)$, $(a,1)$. Since each $3$-subset of
$I_u\times Z_2$ occurs in at most one block and each base block
containing the two points $(a,0)$, $(a,1)$ generates exactly two
different blocks containing the same two points, the number $m$ is
at most $\lfloor(2u-2)/4\rfloor=(2u-4)/4$. Thus there are at least
two $3$-subsets containing the two points $(a,0)$, $(a,1)$ in the
leave, denoted by $\{(a,0),(a,1),(b_a,0)\}$ and
$\{(a,0),(a,1),(b_a,1)\}$, where $b_a\in I_u$ and $b_a\neq a$. Note
that the above conclusion holds for each $a\in I_u$. We have that
${\cal L}\supset
\{\{(a,0),(a,1),(b_a,0)\},\{(a,0),(a,1),(b_a,1)\}:a\in I_u\}$.

Given any $a\in I_u$, consider the number $r$ of the blocks
containing the two points $(a,0)$, $(b_a,0)$. Since each $3$-subset
of $I_u\times Z_2$ occurs in at most one block, the number $r$ is at
most $\lfloor(2u-3)/2\rfloor=(2u-4)/2$. Thus there is at least
another one $3$-subset in the leave containing the two points
$(a,0)$, $(b_a,0)$. Assume that $\{(a,0),(b_a,0),(x,k)\}\in{\cal
L}$, where $(x,k)\neq (a,1)$. Similarly, consider the blocks
containing the two points $(a,0)$, $(b_a,1)$ and assume that
$\{(a,0),(b_a,1),(y,l)\}\in{\cal L}$, where $(y,l)\neq (a,1)$.

If $(x,k)\neq (b_a,1)$, then
$\{(a,0),(b_a,0),(x,k)\}\neq\{(a,0),(b_a,1),(y,l)\}$. Since the
number of $3$-subsets containing the point $(a,0)$ in the leave is
exactly four, they must be $\{(a,0),(a,1),(b_a,0)\}$,
$\{(a,0),(a,1),(b_a,1)\}$, $\{(a,0),(b_a,0),(x,k)\}$,
$\{(a,0),(b_a,1),(y,l)\}$. Note that $(x,k)\neq (b_a,0)$. Consider
the number $s$ of the blocks containing the two points $(a,0)$,
$(x,k)$. Since each $3$-subset of $I_u\times Z_2$ occurs in at most
one block, the number $s$ is at most
$\lfloor(2u-3)/2\rfloor=(2u-4)/2$. Thus there is at least another
one $3$-subset in the leave containing the two points $(a,0)$,
$(x,k)$. It implies that $(y,l)=(x,k)$. Due to
$\{(a,0),(b_a,1),(y,l)\}\in {\cal L}$, i.e.,
$\{(a,0),(b_a,1),(x,k)\}\in {\cal L}$, under the action of $Z_2$ we
have $\{(a,1),(b_a,0),(x,k+1)\}\in {\cal L}$. It implies that there
are at least five $3$-subsets containing the point $(b_a,0)$, i.e.,
$\{(a,0),(a,1),(b_a,0)\}$, $\{(a,0),(b_a,0),(x,k)\}$,
$\{(a,1),(b_a,0),(x,k+1)\}$, $\{(b_a,0),(b_a,1)$, $(b_{b_a},0)\}$,
$\{(b_a,0),(b_a,1),(b_{b_a},1)\}$. A contradiction.

If $(x,k)=(b_a,1)$ and $(y,l)\neq (b_a,0)$, then
$\{(a,0),(b_a,0),(x,k)\}\neq\{(a,0),(b_a,1),(y,l)\}$. Note that
$(y,l)\neq (b_a,1)$, and hence $(y,l)\neq (x,k)$. Since the number
of $3$-subsets containing the point $(a,0)$ in the leave is exactly
four, they must be $\{(a,0),(a,1),(b_a,0)\}$,
$\{(a,0),(a,1),(b_a,1)\}$, $\{(a,0),(b_a,0),(x,k)\}$,
$\{(a,0),(b_a,1)$, $(y,l)\}$. It implies that there is only one
$3$-subset containing the two points $(a,0)$ and $(y,l)$ in the
leave. Similar arguments to those in the paragraph $4$ of this
proof, it is impossible.

If $(x,k)=(b_a,1)$ and $(y,l)=(b_a,0)$, then
$\{(a,0),(b_a,0),(x,k)\}=\{(a,0),(b_a,1),(y,l)\}$. Since the number
of $3$-subsets containing the point $(a,0)$ in the leave is exactly
four, three of them must be $\{(a,0),(a,1),(b_a,0)\}$,
$\{(a,0),(a,1),(b_a,1)\}$ and $\{(a,0),(b_a,0),(b_a,1)\}$. Assume
that the $4$th $3$-subset containing the point $(a,0)$ is
$\{(a,0),(z,i),(w,j)\}$. Similar arguments to those  in the
paragraph $4$ of this proof, we have that the number of $3$-subsets
containing the points $(a,0)$, $(z,i)$ in the leave must be even. A
contradiction. Hence $\Phi(u\times 2,4,2)\leq
\lfloor\frac{u}{4}(\lfloor\frac{2u-1}{3}\lfloor\frac{2u-2}{2}
\rfloor\rfloor-1)\rfloor-1$.

Next consider the number of $\Phi(u\times v,4,2)$. If there is an
optimal $2$-D $(u\times v,4,2)$-OOC with $\Phi(u\times v,4,2)$
codewords, then by Theorem \ref{relation2}, for integer
factorization $v=2v_1$, there exits a $2$-D $(uv_1\times 2,4,2)$-OOC
with $v_1\Phi(u\times v,4,2)$ codewords. Since $uv_1\equiv 0$ $({\rm
mod}$ $12)$, we have $v_1\Phi(u\times v,4,2)\leq
\lfloor\frac{uv_1}{4}(\lfloor\frac{uv-1}{3}\lfloor\frac{uv-2}{2}
\rfloor\rfloor-1)\rfloor-1$. It is readily checked that
$\Phi(u\times v,4,2)\leq
\lfloor\frac{1}{v_1}(\lfloor\frac{uv_1}{4}(\lfloor\frac{uv-1}{3}\lfloor\frac{uv-2}{2}
\rfloor\rfloor-1)\rfloor-1)\rfloor=u(u^2v^2-3uv-6)/24-1=\lfloor\frac{u}{4}(\lfloor\frac{uv-1}{3}\lfloor\frac{uv-2}{2}
\rfloor\rfloor-1)\rfloor-1$.  \qed

\vspace{0.4cm}
\noindent \textbf{Lemma \ref{bound-u-v-4-6}}\hspace{0.2cm} {\it Let $uv\equiv 0$ $({\rm mod}$ $12)$ and $v\equiv 0$ $({\rm mod}$ $6)$. Then $\Phi(u\times v,4,2)\leq
\lfloor\frac{u}{4}(\lfloor\frac{uv-1}{3}$ $\lfloor\frac{uv-2}{2}
\rfloor\rfloor-2)\rfloor$.}

\vspace{0.4cm}
\proof Let $\cal L$ be the leave of a strictly $v$-cyclic
$3$-$(u\times v,4,1)$-packing. Let ${\cal
L}_1=\{\{(a,i),(a,v/3+i),(a,2v/3+i)\}:a\in I_u,0\leq i<v/3\}$. Since
each orbit of $3$-subsets in ${\cal L}_1$ is of length $v/3$ under
the action of $Z_v$, each $3$-subset in ${\cal L}_1$ must be
contained in $\cal L$, i.e., ${\cal L}_1\subset {\cal L}$.

For each $a\in I_u$ and each $i\in Z_v$, consider the number $n$ of
the blocks containing the two points $(a,i)$, $(a,v/3+i)$. Since
each $3$-subset of $I_u\times Z_v$ occurs in at most one block and
$\{(a,i),(a,v/3+i),(a,2v/3+i)\}\in {\cal L}$, the number $n$ is at
most $\lfloor(uv-3)/2\rfloor=(uv-4)/2$. Thus there is at least
another one $3$-subset in the leave containing the two points
$(a,i)$, $(a,v/3+i)$, denoted by
$\{(a,i),(a,v/3+i),(b_{a,i},j_{a,i})\}$, where
$(b_{a,i},j_{a,i})\neq (a,2v/3+i)$. Note that the above conclusion
holds for each $a\in I_u$ and each $i\in Z_v$. Thus we have that
${\cal L}_2=\{\{(a,i),(a,v/3+i),(b_{a,i},j_{a,i})\}:a\in I_u,i\in
Z_v\}\subset{\cal L}\setminus{\cal L}_1$.

For each $a\in I_u$ and each $0\leq i<v/2$, consider the number $m$
of the base blocks containing the two points $(a,i)$, $(a,v/2+i)$.
Since each $3$-subset of $I_u\times Z_v$ occurs in at most one block
and each base block containing the two points $(a,i)$, $(a,v/2+i)$
generates exactly two different blocks containing the same two
points, the number $m$ is at most $\lfloor(uv-2)/4\rfloor=(uv-4)/4$.
Thus there are at least two $3$-subsets containing the two points
$(a,i)$, $(a,v/2+i)$ in the leave, denoted by
$\{(a,i),(a,v/2+i),(c_{a,i},k_{a,i})\}$ and
$\{(a,i),(a,v/2+i),(c_{a,i},v/2+k_{a,i})\}$. Note that the above
conclusion holds for each $a\in I_u$ and each $0\leq i<v/2$. Thus we
have that ${\cal
L}_3=\{\{(a,i),(a,v/2+i),(c_{a,i},k_{a,i})\},\{(a,i),(a,v/2+i),(c_{a,i},v/2+k_{a,i})\}:a\in
I_u,0\leq i<v/2\}\subset{\cal L}\setminus{\cal L}_1$. For
convenience assume that ${\cal
L}_3=\{\{(a,i),(a,v/2+i),(c_{a,i},l_{a,i})\}:a\in I_u,i\in Z_v\}$,
where $l_{a,i}=k_{a,i}$ when $0\leq i<v/2$, and
$l_{a,i}=v/2+k_{a,i}$ when $v/2\leq i<v$.

If ${\cal L}_2\cap{\cal L}_3=\emptyset$, then $|{\cal L}|\geq
7uv/3$. If ${\cal L}_2\cap{\cal L}_3\neq\emptyset$, assume that
$\{(x,i_1),(x,v/3+i_1),(b_{x,i_1},j_{x,i_1})\}=\{(x,i_2),(x,v/2+i_2),(c_{x,i_2},l_{x,i_2})\}$
for some $x\in I_u$ and some $i_1,i_2\in Z_v$. If $(x,i_1)\neq
(c_{x,i_2},l_{x,i_2})$, we have $(b_{x,i_1},j_{x,i_1})=(x,v/2+i_1)$.
If $(x,i_1)=(c_{x,i_2},l_{x,i_2})$, we have
$(b_{x,i_1},j_{x,i_1})=(x,5v/6+i_1)$. Thus each $3$-subset in ${\cal
L}_2\cap{\cal L}_3$ is of the form
$\{(x,i_1),(x,v/3+i_1),(x,v/2+i_1)\}$ or
$\{(x,i_1),(x,v/3+i_1),(x,5v/6+i_1)\}$. Let ${\cal L}_2\cap{\cal
L}_3=\{\{(x,i),(x,v/3+i),(x,v/2+i)\}:x\in A,i\in Z_v\}\cup
\{\{(x,i),(x,v/3+i),(x,5v/6+i)\}:x\in B,i\in
Z_v\}=\{\{(x,i),(x,v/3+i),(x,v/2+i)\}:x\in A,i\in Z_v\}\cup
\{\{(x,v/2+i),(x,5v/6+i),(x,v/3+i)\}:x\in B,i\in Z_v\}$, where
$A,B\subset I_u$ and $A\cap B=\emptyset$. Then ${\cal
L}_2\setminus({\cal L}_2\cap{\cal
L}_3)=\{\{(a,i),(a,v/3+i),(b_{a,i},j_{a,i})\}:a\in I_u\setminus
(A\cup B),i\in Z_v\}$ and ${\cal L}_3\setminus({\cal L}_2\cap{\cal
L}_3)=\{\{(a,i),(a,v/2+i),(c_{a,i},l_{a,i})\}:a\in I_u\setminus
(A\cup B),i\in Z_v\}$.

Let $\{(x,v/3+i),(x,v/2+i)\}\subset T\in {\cal L}_2\cap{\cal L}_3$.
Consider the number of the blocks containing the two points
$(x,v/3+i)$, $(x,v/2+i)$. It is at most
$\lfloor(uv-3)/2\rfloor=(uv-4)/2$. Thus there is at least another
one $3$-subset containing the two points $(x,v/3+i)$, $(x,v/2+i)$ in
${\cal L}\setminus({\cal L}_2\cap{\cal L}_3)$, denoted by
$\{(x,v/3+i),(x,v/2+i),(d_{x,i},r_{x,i})\}$, where
$(d_{x,i},r_{x,i})\neq(x,i)$ if $x\in A$, and
$(d_{x,i},r_{x,i})\neq(x,5v/6+i)$ if $x\in B$. Let ${\cal
L}_4=\{\{(x,v/3+i),(x,v/2+i),(d_{x,i},r_{x,i})\}:\{(x,v/3+i),(x,v/2+i)\}\subset
T\in {\cal L}_2\cap{\cal L}_3\}$. Then ${\cal L}_4\subset {\cal L}$
and ${\cal L}_4\cap({\cal L}_2\cup{\cal L}_3)=\emptyset$. Since
$|{\cal L}_4|=|{\cal L}_2\cap {\cal L}_3|$, we have $|{\cal L}|\geq
|{\cal L}_1|+|{\cal L}_2|+|{\cal L}_3\setminus({\cal L}_2\cap {\cal
L}_3)|+|{\cal L}_4|=7uv/3$.

Thus there are at least $7uv/3$ $3$-subsets in the leave. It implies
that $\Phi(u\times v,4,2)\leq
\lfloor({{uv}\choose{3}}-\frac{7}{3}uv)/(4v)\rfloor=\lfloor\frac{1}{24}u(u^2v^2-3uv-12)\rfloor.$
It is readily checked that
$\lfloor\frac{u}{4}(\lfloor\frac{uv-1}{3}\lfloor\frac{uv-2}{2}
\rfloor\rfloor-2)\rfloor=\lfloor\frac{1}{24}u(u^2v^2-3uv-12)\rfloor.$ \qed

\begin{center}\noindent{\bf Appendix III}\end{center}
\textbf{Proof of Construction \ref{Weighting-strictly $h$-cyclic $s$-FG}:} For checking the correctness of the algorithm shown in Figure $2$, it suffices to show that: ($1$) the resulting design is strictly $h_1h_2$-cyclic; ($2$) any
$3$-subset $S$ of $X'$ satisfying that $|S \cap G'|<3$ for each
$G'\in {\cal G}'$ is contained in a unique block of the resulting
design; ($3$) any $2$-subset $R$ of $X'$ satisfying that $|R \cap
G'|<2$ for each $G'\in {\cal G}'$ is contained in a unique block of
${{\cal A}'_i}$ for each $1\leq i\leq s$.

For convenience assume that ${\cal A}_B=\bigcup_{j=1}^s {\cal A}_B^j$ for each $B\in {\cal F}_1$.

($1$) Suppose that $A=\{(x_l,y_l+u_lg_1,z_l+v_lh_1):1\leq l\leq r\}$
is a base block of the resulting design, where $x_l\in I_n$, $y_l\in
I_{g_1}$, $u_l\in I_{g_2}$, $z_l\in Z_{h_1}$, $v_l\in Z_{h_2}$. We
need to show that the stabilizer of $A$ is trivial, i.e., $A+\delta=
A$ if and only if $\delta\equiv 0$ (mod $h_1h_2)$. The sufficiency
follows immediately, so we consider the necessity. Assume that
$\delta=\delta_1+\delta_2 h_1$, $\delta_1\in Z_{h_1}$, $\delta_2\in
Z_{h_2}$. If $A+\delta=A$, we have
$$\{(x_l,y_l+u_lg_1,z_l+v_lh_1):1\leq l\leq r\}
=\{(x_l,y_l+u_lg_1,z_l+\delta_1+(v_l+\delta_2)h_1):1\leq l\leq r\}
,$$ where the arithmetic is modulo $(-,-,h_1h_2)$. It follows that
$$\{(x_l,y_l,z_l):1\leq l\leq r\}
=\{(x_l,y_l,z_l+\delta_1):1\leq l\leq r\} ,$$ where the arithmetic
is modulo $(-,-,h_1)$. Let $U=\{(x_l,y_l,z_l):1\leq l\leq r\}$.

If $A\in {{\cal A}'_j}$, $1\leq j\leq s$, then $|U|=r\geq 2$. Since
the subdesign $(X,{\cal G},{\cal B})$ of the master design
$1$-FG$(3,(K,K_T),ng_1h_1)$ of type $(g_1h_1)^n$ $(X,{\cal G},{\cal
B},{\cal T})$ is strictly $h_1$-cyclic and it requires that any
$2$-subset of $X$ which intersects each group of ${\cal G}$ in at
most one point occurs in exactly one block, we have $\delta_1=0$.

If $A\in {\cal D}'$, without loss of generality assume that $A\in
{\cal D}^*$. If $A=\tau(C)$ for some $C\in \bigcup_{B\in{\cal F}_2}
{{\cal D}'_B}$, then $|U|=r\geq 3$. Since the master design
$1$-FG$(3,(K,K_T),ng_1h_1)$ of type $(g_1h_1)^n$ is strictly
$h_1$-cyclic and it requires that any $3$-subset of $X$ which
intersects each group of ${\cal G}$ in at most two points occurs in
exactly one block, we have $\delta_1=0$. If $A=\tau(C)$ for some
$C\in \bigcup_{B\in{\cal F}_1} {{\cal D}_B}$, then $|U|\geq 2$. Note
that in this case $U$ may be a multiset, i.e., $|U|$ may be not
equal to $r$. By similar arguments to those in the case of $A\in
{{\cal A}'_j}$, we have $\delta_1=0$.

Hence,
$$\{(x_l,y_l+u_lg_1,z_l+v_lh_1):1\leq l\leq r\}
=\{(x_l,y_l+u_lg_1,z_l+(v_l+\delta_2)h_1):1\leq l\leq r\} ,$$ where
the arithmetic is modulo $(-,-,h_1h_2)$. It follows that
$$\{(x_l,y_l,z_l,u_l,v_l):1\leq l\leq r\}
=\{(x_l,y_l,z_l,u_l,v_l+\delta_2):1\leq l\leq r\} ,$$ where the
arithmetic is modulo $(-,-,-,-,h_2)$. Since the input designs are
all strictly $h_2$-cyclic, we have $\delta_2=0$. Thus the resulting
design is strictly $h_1h_2$-cyclic.

($2$) Take any triple $S=\{(x_l,y_l+u_lg_1,z_l+v_lh_1):1\leq l\leq
3\}\subset X'$, where $x_l\in I_n$, $y_l\in I_{g_1}$, $u_l\in
I_{g_2}$, $z_l\in Z_{h_1}$, $v_l\in Z_{h_2}$ and $x_1,x_2,x_3$ are
not equal at the same time.

Case $1$. Suppose that $x_1$, $x_2$, $x_3$ are pairwise distinct.
Then there exist a unique base block $F$ in ${\cal F}$ and a unique
element $\delta_1\in Z_{h_1}$, such that $\{(x_l,y_l,z_l^*):1\leq
l\leq 3\}\subseteq F$ and
$(x_l,y_l,z_l^*)+\delta_1=(x_l,y_l,z_l^*+\delta_1)=(x_l,y_l,z_l)$,
$1\leq l\leq 3$, where the arithmetic is modulo $(-,-,h_1)$. It
follows that $(x_l,y_l,z_l)-\delta_1=(x_l,y_l,z_l^*)$.

If $F\in {\cal F}_1$, then there exist a unique base block $B\in
{\cal A}_F\bigcup {\cal D}_F$ and a unique element $\delta_2\in
Z_{h_2}$, such that $\{(x_l,y_l,z_l^*,u_l,v_l^*):1\leq l\leq
3\}\subseteq B$ and
$(x_l,y_l,z_l^*,u_l,v_l^*)+\delta_2=(x_l,y_l,z_l^*,u_l,v_l^*+\delta_2)
=(x_l,y_l,z_l^*,u_l,v_l)$, $1\leq l\leq 3$, where the arithmetic is
modulo $(-,-,-,-,h_2)$. By the mapping $\tau$, we have that
$(x_l,y_l+u_lg_1,z_l^*+(v_l^*+\delta_2)h_1)
=(x_l,y_l+u_lg_1,z_l^*+v_lh_1)
  =(x_l,y_l+u_lg_1,z_l-\delta_1+v_lh_1)$, where the arithmetic is modulo $(-,-,h_1h_2)$.
   Let $\delta=\delta_1+\delta_2 h_1$. It
follows that $(x_l,y_l+u_lg_1,z_l^*+v_l^*h_1+\delta)=
(x_l,y_l+u_lg_1,z_l+v_lh_1)$. By $(1)$ the resulting design is
strictly $h_1h_2$-cyclic, so $\{(x_l,y_l+u_lg_1,z_l+v_lh_1):1\leq
l\leq 3\}$ is contained in the unique block $\tau(B)+\delta$, which
is generated by $\tau(B)$. Similar arguments hold for $F\in {\cal
F}_2$, where $B\in {{\cal D}'_F}$.

Case $2$. Suppose that $x_1=x_2$, $x_1\neq x_3$, and $(y_1,z_1)\neq
(y_2,z_2)$. Then there exist a unique base block $F$ in ${\cal F}_2$
and a unique element $\delta_1\in Z_{h_1}$, such that
$\{(x_l,y_l,z_l^*):1\leq l\leq 3\}\subseteq F$ and
$(x_l,y_l,z_l^*+\delta_1)=(x_l,y_l,z_l)$, $1\leq l\leq 3$, where the
arithmetic is modulo $(-,-,h_1)$. There exist a unique base block
$B\in {\cal D}'_F$ and a unique element $\delta_2\in Z_{h_2}$, such
that $\{(x_l,y_l,z_l^*,u_l,v_l^*):1\leq l\leq 3\}\subseteq B$ and
$(x_l,y_l,z_l^*,u_l,v_l^*+\delta_2) =(x_l,y_l,z_l^*,u_l,v_l)$,
$1\leq l\leq 3$, where the arithmetic is modulo $(-,-,-,-,h_2)$. By
similar arguments to those in Case $1$,
$\{(x_l,y_l+u_lg_1,z_l+v_lh_1):1\leq l\leq 3\}$ is contained in the
unique block $\tau(B)+\delta$, where $\delta=\delta_1+\delta_2 h_1$.

Case $3$. Suppose that $(x_1,y_1,z_1)=(x_2,y_2,z_2)$, $(u_1,v_1)\neq
(u_2,v_2)$ and $x_1\neq x_3$. Then there exist a unique base block
$F$ in ${\cal F}_1$ and a unique element $\delta_1\in Z_{h_1}$, such
that $\{(x_l,y_l,z_l^*):1\leq l\leq 3\}\subseteq F$ and
$(x_l,y_l,z_l^*+\delta_1)=(x_l,y_l,z_l)$, $1\leq l\leq 3$, where the
arithmetic is modulo $(-,-,h_1)$. There exist a unique base block
$B\in {\cal D}_F$ and a unique element $\delta_2\in Z_{h_2}$, such
that $\{(x_l,y_l,z_l^*,u_l,v_l^*):1\leq l\leq 3\}\subseteq B$ and
$(x_l,y_l,z_l^*,u_l,v_l^*+\delta_2) =(x_l,y_l,z_l^*,u_l,v_l)$,
$1\leq l\leq 3$, where the arithmetic is modulo $(-,-,-,-,h_2)$. By
similar arguments to those in Case $1$,
$\{(x_l,y_l+u_lg_1,z_l+v_lh_1):1\leq l\leq 3\}$ is contained in the
unique block $\tau(B)+\delta$, where $\delta=\delta_1+\delta_2 h_1$.

($3$) Take any $2$-subset $R=\{(x_l,y_l+u_lg_1,z_l+v_lh_1):1\leq
l\leq 2\}\subset X'$, where $x_l\in I_n$, $y_l\in I_{g_1}$, $u_l\in
I_{g_2}$, $z_l\in Z_{h_1}$, $v_l\in Z_{h_2}$ and
 $x_1\neq x_2$. Then there exist a unique
base block $F$ in ${\cal F}_1$ and a unique element $\delta_1 \in
Z_{h_1}$, such that $\{(x_l,y_l,z_l^*):1\leq l\leq 2\}\subseteq F$
and $(x_l,y_l,z_l^*+\delta_1)=(x_l,y_l,z_l)$, $1\leq l\leq 2$, where
the arithmetic is modulo $(-,-,h_1)$.

Given any $1\leq j\leq s$. There exist a unique base block $B$ in
${\cal A}_j$ and a unique element $\delta_{2}\in Z_{h_2}$, such that
$\{(x_l,y_l,z_l^*,u_l,v_l^*):1\leq l\leq 2\}\subseteq B$ and
$(x_l,y_l,z_l^*,u_l,v_l^*+\delta_2) =(x_l,y_l,z_l^*,u_l,v_l)$,
$1\leq l\leq 2$, where the arithmetic is modulo $(-,-,-,-,h_2)$. By
the mapping $\tau$, we have that
$(x_l,y_l+u_lg_1,z_l^*+(v_l^*+\delta_2)h_1)
=(x_l,y_l+u_lg_1,z_l^*+v_lh_1)
  =(x_l,y_l+u_lg_1,z_l-\delta_1+v_lh_1)$, where the arithmetic is modulo $(-,-,h_1h_2)$.
   Let $\delta=\delta_1+\delta_2 h_1$. It
follows that $(x_l,y_l+u_lg_1,z_l^*+v_l^*h_1+\delta)=
(x_l,y_l+u_lg_1,z_l+v_lh_1)$. By $(1)$ the resulting design is
strictly $h_1h_2$-cyclic, so $\{(x_l,y_l+u_lg_1,z_l+v_lh_1):1\leq
l\leq 2\}$ is contained in the unique block $\tau(B)+\delta$, which
is generated by $\tau(B)$. \qed

\vspace{0.4cm} \noindent \textbf{Proof of Construction \ref{Weighting for
$h$-cyclic $H$ design}:} For checking the correctness of the algorithm shown in Figure $4$, take any $t$-subset $T=\{(x_l,y_l+u_lg_1,z_l+v_lh_1):1\leq l\leq
t\}\subset X'$, where $x_l\in I_n$, $y_l\in I_{g_1}$, $u_l\in
I_{g_2}$, $z_l\in Z_{h_1}$, $v_l\in Z_{h_2}$ and
 $|\{x_l:1\leq l\leq t\}|=t$. Then there exist a unique
base block $F$ in ${\cal F}$ and a unique element $\delta_1 \in
Z_{h_1}$, such that $\{(x_l,y_l,z_l^*):1\leq l\leq t\}\subseteq F$
and $(x_l,y_l,z_l^*+\delta_1)=(x_l,y_l,z_l)$, $1\leq l\leq t$, where
the arithmetic is modulo $(-,-,h_1)$. There exist a unique base
block $B$ in ${\cal D}_F$ and a unique element $\delta_{2}\in
Z_{h_2}$, such that $\{(x_l,y_l,z_l^*,u_l,v_l^*):1\leq l\leq
t\}\subseteq B$ and $(x_l,y_l,z_l^*,u_l,v_l^*+\delta_2)
=(x_l,y_l,z_l^*,u_l,v_l)$, $1\leq l\leq t$, where the arithmetic is
modulo $(-,-,-,-,h_2)$. By the mapping $\tau$, we have that
$(x_l,y_l+u_lg_1,z_l^*+(v_l^*+\delta_2)h_1)
=(x_l,y_l+u_lg_1,z_l^*+v_lh_1)
  =(x_l,y_l+u_lg_1,z_l-\delta_1+v_lh_1)$, where the arithmetic is modulo $(-,-,h_1h_2)$.
   Let $\delta=\delta_1+\delta_2 h_1$. It
follows that $(x_l,y_l+u_lg_1,z_l^*+v_l^*h_1+\delta)=
(x_l,y_l+u_lg_1,z_l+v_lh_1)$. Thus
$\{(x_l,y_l+u_lg_1,z_l+v_lh_1):1\leq l\leq t\}$ is contained in the
unique block $\tau(B)+\delta$, which is generated by $\tau(B)$. \qed

\vspace{0.6cm} \noindent{\bf Acknowledgements} The authors would
like to thank the anonymous referees and Prof. Navin Kashyap, the Associate Editor for Coding Theory, for their
helpful comments and valuable suggestions to improve the readability
of this paper.

\end{document}